\documentclass[a4paper,10pt]{article}
\usepackage{authblk}

\usepackage[ngerman,english]{babel}
\usepackage[latin1]{inputenc}
\usepackage{csquotes}
\usepackage{amsfonts,amsmath,amsthm}
\usepackage{empheq}
\usepackage[titletoc,title]{appendix}
\usepackage[backend=bibtex8,doi=false,eprint=false,firstinits=true,isbn=false,style=numeric-comp,maxnames=99]{biblatex}
\makeatletter
\def\blx@maxline{77}
\makeatother

\bibliography{bibl_SW.bib}

\usepackage{cases}
\usepackage{mathabx}
\usepackage{bbm}
\usepackage{xfrac}
\usepackage{fancyhdr}
\usepackage{color}
\usepackage[colorinlistoftodos,textsize=small]{todonotes}
\usepackage[colorlinks]{hyperref}
\definecolor{blue75}{rgb}{0,0,.75}
\definecolor{green75}{rgb}{0,.75,0}
\hypersetup{colorlinks=true, urlcolor=blue75,linkcolor=blue75,citecolor=green75,pdfstartview=FitB,bookmarksopen=true,bookmarksopenlevel=1}
\usepackage[a4paper, left=2.5cm, right=2.0cm, top=2.5cm,bottom=2cm]{geometry}
\usepackage{constants}
\newcommand{\parenthezises}[1]{\arabic{#1}}
\newconstantfamily{C}{
	symbol=C,
	format=\parenthezises,
	reset={section}
}
\newconstantfamily{M}{
	symbol=M,
	format=\parenthezises,
	reset={section}
}
\newconstantfamily{B}{
	symbol=B,
	format=\parenthezises,
	reset={section}
}
\newconstantfamily{xi}{
	symbol=\xi,
	format=\parenthezises,
	reset={section}
}
\usepackage{enumerate}
\usepackage{graphicx}
\graphicspath{ {images/} }
\usepackage{wrapfig}
\usepackage{figbib}
\allowdisplaybreaks
\usepackage[capitalise]{cleveref}

\crefdefaultlabelformat{{\it #2#1#3}}

\crefname{equation}{}{}

\crefname{enumi}{}{}
\creflabelformat{enumi}{{(#2#1#3)}}

\crefname{section}{{\it Section}}{{\it Sections}}
\crefname{subsection}{{\it Subsection}}{{\it Subsections}}
\crefname{subsubsection}{{\it Paragraph}}{{\it Paragraphs}}

\newtheorem{Theorem}{Theorem}[section]
\crefname{Theorem}{{\it Theorem}}{{\it Theorems}}

\crefname{Definition}{{\it Definition}}{{\it Definitions}}
\newtheorem{Lemma}[Theorem]{Lemma}
\crefname{Lemma}{{\it Lemma}}{{\it Lemmas}}

\crefname{Proposition}{{\it Proposition}}{{\it Propositions}}

\crefname{Assumption}{{\it Assumption}}{{\it Assumptions}}
\newtheorem{Corollary}[Theorem]{Corollary}
\crefname{Corollary}{{\it Corollary}}{{\it Corollaries}}

\theoremstyle{definition}

\crefname{Remark}{{\it Remark}}{{\it Remarks}}

\crefname{Notation}{{\it Notation}}{{\it Notations}}

\crefname{Example}{{\it Example}}{{\it Examples}}

\newcommand{\proofc}{{\sc Proof} \ }
\newcommand{\be}{\begin{equation} \label}
\newcommand{\ee}{\end{equation}}
\newcommand{\bea}{\begin{eqnarray}\label}
\newcommand{\eea}{\end{eqnarray}}
\newcommand{\bas}{\begin{eqnarray*}}
	\newcommand{\eas}{\end{eqnarray*}}
\newcommand{\bit}{\begin{itemize}}
	\newcommand{\eit}{\end{itemize}}

\newcommand{\nn}{\nonumber}
\newcommand{\R}{\mathbb{R}}
\newcommand{\N}{\mathbb{N}}
\newcommand{\pO}{\partial\Omega}
\newcommand{\bom}{\overline{\Omega}}
\newcommand{\na}{\nabla}
\newcommand{\eps}{\varepsilon}

\newcommand{\io}{\int_\Omega}

\newcommand{\abs}{\\[5pt]}
\newcommand{\hs}{\hspace*}

\newcommand{\ov}{\overline}

\newcommand{\tm}{T_{max}}
\newcommand{\Cchi}{C_\chi}
\newcommand{\Cxi}{C_\xi}
\newcommand{\Cphi}{C_\phi}
\newcommand{\Cpsi}{C_\psi}
\newcommand{\tCphi}{\widetilde{C}_\phi}
\newcommand{\tCpsi}{\widetilde{C}_\psi}

\newcommand{\vthet}{\vartheta}
\newcommand{\om}{\omega }
\newcommand{\Sbb}{\mathbb S^{N-1}}
\begin{document}
\enlargethispage{10mm}
\title{Mathematical modeling and analysis of a tumor invasion problem with angiogenesis and taxis cascade}
\author{
Christina Surulescu\footnote{surulescu@mathematik.uni-kl.de}\\
{\small Department of Mathematics, RPTU Kaiserslautern-Landau, } \\
{\small 67663 Kaiserslautern, Germany}
\and
Michael Winkler\footnote{michael.winkler@math.uni-paderborn.de}\\
{\small Institut f\"ur Mathematik, Universit\"at Paderborn,}\\
{\small 33098 Paderborn, Germany} 
}

\date{}
\maketitle
\begin{abstract}
\noindent 
We propose a mathematical model for tumor invasion supported by angiogenesis and interactions with the surrounding tissue. For the model deduction we employ a multiscale approach starting from lower scales and obtaining by an informal parabolic upscaling a system of reaction-diffusion-taxis equations with a so-called 'taxis cascade', where one species is performing taxis towards a signal  whose production/decay is controled by the other, for which it also serves as a tactic cue. We prove global existence and uniqueness of solutions to the obtained PDE-ODE system and perform numerical simulations to illustrate the behavior of solutions. 
\abs
\end{abstract}
%
%
%
%
%
%
%
%
\section{Introduction and model set up}\label{intro}
\subsection{Biological motivation and previous continuous  models}\label{subsec:bio-motiv}
Cancer metastasis is crucially influenced by tumor cell migration and interactions with the peritumoral environment. Angiogenesis is one of the cancer hallmarks \cite{Hanahan2011}. It involves the growth and development of new blood capillaries by chemotactically  attracting endothelial cells (ECs) from neighboring vessels. The environment of a growing tumor becomes increasingly hypoxic, due to sustained glycolytic activity of the cancer cells \cite{Abakarova1995,Gatenby2004}. To further develop and metastasize, tumors require blood supply: they release pro-angiogenic factors (e.g., vascular endothelial growth factor = VEGF) that serve as chemoattractants for ECs \cite{Wang2019}.  Hence, they stimulate the formation of new blood vessels, eventually forming a network of capillaries that infiltrate the neoplasm, thus providing it with nutrients
\cite{Hardee2012,Xu2002}. Tumor angiogenesis is often targeted in cancer therapy, in order to inhibit tumor expansion. \\[-2ex]

\noindent
Migrating cancer cells are also guided by gradients of surrounding tissue, see e.g., \cite{Oudin2016}. For instance, this is particularly prominent in glioma cells, which are known to follow white matter tracts and blood vessels to infiltrate the brain at distant sites from the original tumor \cite{Giese1996,Liu2019,EHKS}. Such behavior is typically described by haptotaxis in the framework of continuous reaction-diffusion-drift models.  \\[-2ex]

\noindent
Previous continuous models describing tumor angiogenesis and haptotaxis typically consider PDE-ODE systems set up in a more or less heuristic manner and coupling the dynamics of endothelial cells and extracellular matrix (ECM) or its components \cite{Xu2020}, possibly with a further PDE for some chemotactic signal \cite{Anderson1998,MORALESRODRIGO2013,Chen2017,Pang2019}. Some of these hence belong to the class of models with multiple taxis reviewed in \cite{Kolbe2021}. For a recent review of mathematical models for tumor angiogenesis we refer to \cite{Hormuth2021}. In this note we propose a model for tumor invasion guided by angiogenesis (here described by endothelial cell dynamics) and by ECM. Our multiscale approach connects meso- and macroscopic dynamics of tumor and endothelial cells, ECM, and VEGF and leads to a PDE-ODE system of reaction-diffusion-taxis equations where the ECs perform taxis towards VEGF gradients, thus following a signal produced by the tumor cells. The latter, in turn, follow gradients of EC density. Such tactic scenario is sometimes called a 'taxis cascade', see e.g. \cite{Tao2019,Winkler2019,Tao2023} for its use in another context.

\subsection{Model set up}\label{subsec:model}

We aim to deduce a macroscopic description for the dynamics of a population of tumor cells migrating in a tissue (mainly ECM) and interacting with endothelial cells. The cancer cells preferentially bias their motion towards gradients of tissue and EC densities, while the ECs follow the gradient of VEGF expressed by the tumor cells, which also degrade the ECM. Hence, the obtained mathematical model should involve partial differential equations with diffusion and taxis terms describing the mentioned biases. We follow the informal derivation performed e.g., in \cite{CS20,Conte2023,Kumar2021} within the kinetic theory of active particles (KTAP) framework \cite{bellom3} and adapt it here to the biological problem sketched above. 

\subsubsection{Micro-meso description}\label{ssec:micro-meso}

Cells perceive signals from their surroundings by way of their receptors binding to the respective ligands and of subsequent processing the therewith transmitted biochemical and/or biophysical information. Thus, we start by describing simple mass action kinetics for receptor binding (subcellular level) of the involved cell types (cancer cells, ECs) to their respective environmental cues (ECs and tissue for tumor cells, VEGF for ECs). \\[-2ex]

\noindent
Let $y_1$ and $y_2$ denote the amounts of cancer cell receptors bound to tissue and ECs, respectively. We will also use $u(t,x)$, $w(t,x)$, and $v(t,x)$ to denote the macroscopic densities (volume fractions) of tumor cells, tissue, and ECs, respectively, as well as $z(t,x)$ for the concentration of VEGF. Further, the mesoscopic quantities $p(t,x,\vthet_c,y)$ and $\om (t,x,\vthet_e)$ represent the distribution functions of tumor cells and of ECs, respectively, each of them depending on a kinetic variable $\vthet_{\cdot}$ representing cell velocity. We take $t\ge 0$, $x\in \R^N$, $\vthet_c\in s\Sbb =:\Theta _c$, $\vthet_e\in \sigma \Sbb =:\Theta _e$, where $\Sbb $ is the unit sphere in $\R^N$ and we assume constant speeds for the tumor cells and for the ECs:  $s=\frac{\vthet_c}{|\vthet_c|}>0$, $\sigma =\frac{\vthet_e}{|\vthet_e|}>0$; hence, only the directions of the velocity vectors change.
\begin{align*}
&{\frac{v}{v_M}+(R_0-y_1-y_2)} \xrightleftharpoons[k_1^{-}]{k_1^{+}} {y_1},\\
&{\frac{w}{w_M}+(R_0-y_1-y_2)} \xrightleftharpoons[k_2^{-}]{k_2^{+}} {y_2},
\end{align*}
where $v_M$ and $w_M$ denote reference quantities for $v$ and $w$, respectively; they are used for purposes of nondimensionalization occurring later on. Here $R_0$ denotes the total amount of receptors on a cell membrane, which for simplicity we assume to be constant. If we rescale $y_1$ and $y_2$ by $R_0$, while keeping the same notation, and consider $y:=y_1+y_2$ to be the total amount of bound receptors, then we can describe the receptor binding dynamics on a cancer cell by 
\begin{equation}
\dot y=k_1^+\frac{v}{v_M}+k_2^+\frac{w}{w_M}-(k_1^+\frac{v}{v_M}+k_2^+\frac{w}{w_M}+k^-)y,
\end{equation}
where we also assumed that the cell detachment rates from tissue and ECs are equal: $k_1^-=k_2^-=:k^-$. The so-called activity variable $y$ is supposed to belong to the open interval $(0,1)$; for more details on this we refer to \cite{Lorenz2014}. As in \cite{CS20,CEKNSSW,CKSS,EHKS,EHS,EKS,Kumar2021}, we consider deviations $\zeta :=y_*-y$ from the steady-state 
\begin{equation}
y_*=\frac{k_1^+\frac{v}{v_M}+k_2^+\frac{w}{w_M}}{B(v,w)},
\end{equation}
where we denoted $B(v,w):=k_1^+\frac{v}{v_M}+k_2^+\frac{w}{w_M}+k^-$. They satisfy the ODE
\begin{equation}
\dot \zeta =-B(v,w)\zeta +\frac{k^-}{(B(v,w))^2}\left (\frac{k_1^+}{v_M}D_tv+\frac{k_2^+}{w_M}D_tw\right )=:G(\zeta, v,w),
\end{equation}
with $D_tv=v_t+\vthet_c\cdot \nabla v$ denoting the pathwise gradient of $v$ (analogously $D_tw$). Since the receptor binding dynamics is very fast when compared to the other biological processes described here, we will assume in the following that the deviations $\zeta $ are very small. We denote by $Z\subseteq [y_*-1,y_*]$ the set to which the shifted activity variable $\zeta $ belongs. \\[-2ex]

\noindent
On the mesoscopic level the following kinetic transport equations (KTEs) hold for the cell distribution functions $p(t,x,\vthet_c,y)$ and $\om (t,x,\vthet_e)$:
\begin{align}
&p_t+\na _x\cdot (\vthet_cp)+\partial_{\zeta}(G(\zeta,v,w)p)=\mathcal L[\lambda(\zeta )]p\label{KTE-p}\\
&\om_t+\na _x\cdot (\vthet_e\om )=\mathcal L[\eta (\vthet_e,z)]\om , \label{KTE-om}
\end{align}
where the turning operators on the right hand sides characterize cell reorientations in response to extracellular influences:
\begin{align}\label{eq:turning-op-cancer}
L[\lambda(\zeta )]p(t,x,\vthet_c,\zeta ):=-\lambda (\zeta )p(t,x,\vthet_c,\zeta )+\lambda (\zeta )\int _{\Theta_c}K(\vthet_c,\vthet_c')p(t,x,\vthet_c',\zeta )\ d\vthet_c'=\lambda (\zeta )\left (\frac{\bar p}{|\Theta_c|}-p\right )
\end{align}
for the cancer cells, where we chose a uniform turning kernel $K(\vthet_c,\vthet_c')=\frac{1}{|\Theta_c|}$ and denoted $\bar p(t,x,\zeta )=\int _{\Theta_c}p(t,x,\vthet_c,\zeta )\ d\vthet_c$. In \cref{eq:turning-op-cancer} the coefficient $\lambda (\zeta )$ denotes the turning rate and, following \cite{CS20,EHKS,EHS,EKS,Kumar2021}, to which we also refer for more details about calculations in this Subsection, we take it in the affine form $\lambda (\zeta)=\lambda _0-\lambda_1\zeta$, with $\lambda_0,\ \lambda_1>0$ constants. \\[-2ex]

\noindent
On the other hand, we let the VEGF concentration act upon reorientations of ECs by way of selecting the turning rate in \cref{KTE-om} of the form
\begin{equation}
\eta (\vthet_e,z):=\eta _0e^{-D_tz/z_M},\qquad D_t(z/z_M)=z_t/z_M+\vthet_e\cdot \na_xz/z_M,\qquad \eta _0>0 \text{ constant}.
\end{equation}
We therewith have the turning operator for ECs (still with a uniform turning kernel)
\begin{align}\label{eq:turning-op-ECs}
L[\eta(\vthet_e,z )]\om (t,x,\vthet_e)=-\eta (\vthet_e,z)\om (t,x,\vthet_e)+\frac{1}{|\Theta_e|}\int _{\Theta_e}\eta (\vthet_e',z)\om (t,x,\vthet_e')\ d\vthet_e'.
\end{align}
\noindent
This approach to including effects of environmental signals into KTE descriptions of mesoscopic cell density dynamics was proposed in \cite{OtHi} for bacteria swimming and adapted in \cite{CS20,Kumar2021} to the case of cancer cell migration. It provides an alternative way\footnote{yet another one can be found in \cite{CEKNSSW,DKSS20} or -more rigorously- in \cite{ZS22}} to the previous receptor binding dynamics for obtaining chemotaxis terms on the macrolevel. An informal relationship between these approaches was studied in \cite{Kumar2021}, a rigorous one (in the case of bacteria) in \cite{PeTaWa}. 

\subsubsection{Parabolic scaling}

\noindent
We perform a parabolic scaling of the mesoscopic system \cref{KTE-p}, \cref{KTE-om}, i.e. we rescale $t\leadsto \eps^2t$, $x\leadsto \eps x$. At his stage we also want to include some proliferation terms for the two types of cells. Since cells need a much longer time to proliferate than perform other functions (in particular migration), such terms will be correspondingly scaled (by $\eps^2$), as done e.g. in  \cite{CS20,CEKNSSW,EHS,Kumar2021}. Thus, \cref{KTE-p}, \cref{KTE-om} become
\begin{align}
\eps^2p_t+\eps \na _x\cdot (\vthet_cp)-&\partial_{\zeta}\left(\left (\zeta B(v,w)-\frac{k^-}{(B(v,w))^2}\left (\frac{k_1^+}{v_M}(\eps ^2\frac{v_t}{v_M}+\eps \vthet_c\cdot \na_xv)+\frac{k_2^+}{w_M}(\eps ^2\frac{w_t}{w_M}+\eps \vthet_c\cdot \na_xw)\right )\right )p\right )\notag\\
=&\mathcal L[\lambda(\zeta )]p+\eps^2\mu_c(u,v,w)\int _Z\Gamma (x,\zeta, \zeta')p(t,x,\vthet_c,\zeta')\ d\zeta ',\label{KTE-p-parab}\\
\eps^2\om_t+\eps \na _x\cdot (\vthet_e\om )=&\mathcal L[\eta ^\eps(\vthet_e,z)]\om +\eps^2\mu_e(u,v,z)\om , \label{KTE-om-parab}
\end{align}
where $u(t,x)=\int_{\Theta_c}\int _Zp(t,x,\vthet_c,\zeta)\ d\zeta d\vthet_c$,  $v(t,x)=\int _{\Theta_e}\om (t,x,\vthet_e)\ d\vthet_e$, and $\Gamma (x,\zeta, \zeta')$ is a kernel satisfying $\int _Z\Gamma (x,\zeta, \zeta')d\zeta =1$ and characterizing the transition from the (deviation of) cell receptor binding state $\zeta'$ to $\zeta $ during interaction with signals (tissue, inter- and intrapopulation exchange); thereby only the volume fractions (and not local directionality) are relevant, hence only macroscopic quantities are considered. In \cref{KTE-om-parab} we denoted 
\begin{align*}
\eta ^\eps(\vthet_e,z)=\eta _0\exp(-(\eps ^2z_t/z_M+\eps \vthet_e\cdot \na_xz/z_M)). 
\end{align*}
Using the moments w.r.t. $\zeta $:
\begin{align*}
m(t,x,\vthet_c)=\int  \limits_Zp(t,x,\vthet_c,\zeta )\ d\zeta,\quad m^\zeta (t,x,\vthet_c)=\int  \limits_Z\zeta p(t,x,\vthet_c,\zeta )\ d\zeta,\quad u^\zeta(t,x)=\int \limits _{\Theta_c}\int \limits _Z\zeta p(t,x,\vthet_c,\zeta )\ d\zeta d\vthet_c,
\end{align*}
along with their Hilbert expansions $m=\sum \limits _{i=0}^{\infty}\eps ^im_i$, $m^\zeta=\sum \limits _{i=0}^{\infty}\eps ^im^\zeta_i$, and $u^\zeta=\sum\limits  _{i=0}^{\infty}\eps ^iu^\zeta_i$, a linearization of $\eta ^\eps$, neglecting the higher order moments w.r.t. $\zeta $ in virtue of the smallness assumption about these deviations, and equating the powers of $\eps$ in \cref{KTE-p-parab}, \cref{KTE-om-parab} we obtain in the usual way (for details see e.g. \cite{CS20,CEKNSSW,Kumar2021}) the following macroscopic equations for $u$ and $v$:
\begin{align}
&u_t=\na \cdot (\mathbb D_T\na u)-\na \cdot \left (\left (g(v,w)\frac{k_1^+}{v_M}\mathbb D_T\na v+g(v,w)\frac{k_2^+}{w_M}\mathbb D_T\na w\right)u\right )+\mu_c(u,v,w)u\label{q-macro-u}\\
&v_t=\na \cdot (\mathbb D_E\na v)-\na \cdot \left (\eta_0v\ \mathbb D_E\na z/z_M\right )+\mu_e(u,v,z)v,\label{q-macro-v}
\end{align}
where:
\begin{align}
&\mathbb D_T=\frac{1}{\lambda_0}\int _{\Theta_c}\vthet_c\otimes \vthet_c \frac{1}{|\Theta_c|}\ d\vthet_c=\frac{s^2}{N\lambda_0}\mathbb I_N\quad \text{(tumor diffusion tensor)}\\
&\mathbb D_E=\frac{1}{\eta_0}\int _{\Theta_e}\vthet_e\otimes \vthet_e \frac{1}{|\Theta_e|}\ d\vthet_e=\frac{\sigma ^2}{N\eta_0}\mathbb I_N\quad \text{(EC diffusion tensor)}\\
&g(v,w)=\frac{\lambda_1k^-}{(B(v,w))^2(\lambda_0+B(v,w))}\quad \text{(tactic sensitivity of tumor cells)}.
\end{align}
The concrete forms of coefficient functions $\mu_c,\mu_e$ will be provided in \cref{sec:numerics}.\\[-2ex]

\noindent
Performing an appropriate nondimensionalization leads to the system\footnote{we kept the previous notations for the scaled $\mu_c$, $\mu_e$, and $B$ functions}
\be{-1}
\left\{ \begin{array}{l}
u_t=\Delta u - \nabla \cdot \Big( u \chi(v,w)\nabla v\Big) - \nabla \cdot \Big( u\xi(v,w)\nabla w \Big)+\mu_c(u,v,w)u,\\
v_t = \Delta v - \nabla \cdot (v\nabla z)+\mu_e(u,v,z)v, 
\end{array} \right.
\ee
with
\begin{align}\label{eq:chi-xi}
\chi(v,w)=\frac{\kappa _1}{(B(v,w))^2(1+B(v,w))}, \qquad \xi(v,w)=\frac{\kappa _2}{(B(v,w))^2(1+B(v,w))},
\end{align}
$\kappa_1, \kappa_2>0$ constants.\\[-2ex]

\noindent
Proceeding e.g., as in \cite{CEKNSSW}, no-flux boundary conditions can be obtained (still in a formal way) for the formulation on a bounded space domain $\Omega \subset \R^N$.  This system is to be supplemented with equations for the evolution of the macroscopic quantities $w$ and $z$, for which no deduction from lower scales is needed. We assume that the tissue (of density $w$) is degraded by the tumor cells\footnote{Unlike other works with stronger modeling focus (e.g., \cite{Lorenz2014,CEKNSSW}) here we do not take into account the fine, possibly anisotropic structure of tissue fibers and their degradation by tumor cells which depends on directionality and orientation of tissue and cells, but  consider a depletion due to the mere cell-fibers interaction.} and that VEGF (of concentration $z$) is uptaken by ECs and produced by cancer cells. This leads to the (nondimensional) equations
\be{-1p2}
\left\{ \begin{array}{l}
w_t=-\psi(uw)\\
z_t=D_z\Delta z-\mu_zvz+\phi (u),	
\end{array} \right.	
\ee
where $D_z,\mu_z>0$ are constants, and with functions $\psi, \phi$ yet to be specified. 
\section{Global existence and uniqueness of solutions}
\subsection{Main results}\label{subsec:announce_results}
In the sequel we consider a version of system \cref{-1} together with the equations for $w$ and $z$, and with more general tactic sensitivity functions $\chi ,\xi$. However, since we are mainly interested in the dynamics due to taxis and self-diffusion, we will not take into account any source terms for the two types of cells (tumor and ECs). For simplicity of writing we will omit the constants in \eqref{-1p2}, but reconsider them when performing numerical simulations in \cref{sec:numerics}.\\[-2ex]

\noindent
Concretely, we study here the following system:
\be{0}
	\left\{ \begin{array}{ll}
	u_t=\Delta u - \nabla \cdot \Big( u \chi(v,w)\nabla v\Big) - \nabla \cdot \Big( u\xi(v,w)\nabla w \Big),
	\qquad & x\in\Omega, \ t>0, \\[1mm]
	v_t = \Delta v - \nabla \cdot (v\nabla z), 
	\qquad & x\in\Omega, \ t>0, \\[1mm]
	w_t = - \psi(uw),
	\qquad & x\in\Omega, \ t>0, \\[1mm]
	z_t = \Delta z - vz + \phi(u),
	\qquad & x\in\Omega, \ t>0, \\[1mm]
	\frac{\partial u}{\partial\nu} - u \chi(v,w) \frac{\partial v}{\partial\nu} - u\xi(v,w)\frac{\partial w}{\partial\nu}
	= \frac{\partial v}{\partial\nu} = \frac{\partial z}{\partial\nu} =0,
	\qquad & x\in\pO, \ t>0, \\[1mm]
	u(x,0)=u_0(x), \quad v(x,0)=v_0(x), \quad w(x,0)=w_0(x), \quad z(x,0)=z_0(x),
	\quad & x\in\Omega,
	\end{array} \right.
\ee
in a bounded convex domain $\Omega\subset \R^2$ with smooth boundary.
Our assumptions on the initial data are that
\be{init}
	\left\{ \begin{array}{l}
	u_0 \in C^{2+\vartheta}(\bom)
	\quad \mbox{is nonnegative with $\frac{\partial u_0}{\partial\nu}=0$ on $\pO$ and } u_0\not\equiv 0, \\[1mm]
	v_0 \in C^{2+\vartheta}(\bom)
	\quad \mbox{is nonnegative with $\frac{\partial v_0}{\partial\nu}=0$ on $\pO$ and } v_0\not\equiv 0, \\[1mm]
	w_0 \in C^{2+\vartheta}(\bom)
	\quad \mbox{is nonnegative with $\frac{\partial w_0}{\partial\nu}=0$ on $\pO$, and that} \\[1mm]
	z_0 \in C^{2+\vartheta}(\bom)
	\quad \mbox{is nonnegative with $\frac{\partial z_0}{\partial\nu}=0$ on $\pO$ and } z_0\not\equiv 0, \\[1mm]
	\end{array} \right.
\ee
with some $\vartheta\in (0,1)$, and concerning the coefficient functions in \cref{0} we shall suppose that
\be{reg}
	\phi\in C^1([0,\infty)),
	\quad 
	\psi\in C^3([0,\infty)),
	\quad
	\chi \in C^2([0,\infty)^2)
	\quad \mbox{and} \quad
	\xi \in C^3([0,\infty)^2),
\ee
that
\be{chi}
	|\chi(v,w)| \le \frac{\Cchi}{v+1}
	\qquad \mbox{for all } (v,w)\in [0,\infty)^2
\ee
and
\be{xi}
	|\xi(v,w)| \le \Cxi
	\quad \mbox{and} \quad
	|\xi_v(v,w)| \le \Cxi
	\qquad \mbox{for all } (v,w)\in [0,\infty)^2,
\ee
and that
\be{psi}
	0 \le s\psi'(s) \le \Cpsi
	\qquad \mbox{for all } s\ge 0
\ee
as well as
\be{phi}
	\phi(u) \ge 0
	\quad \mbox{and} \quad
	u|\phi'(u)| \le \Cphi
	\qquad \mbox{for all } u\ge 0
\ee
with some positive constants $\Cchi,\Cxi,\Cpsi$ and $\Cphi$.\footnote{Observe that the coefficient functions $\chi $ and $\xi$ obtained in \eqref{eq:chi-xi} satisfy these assumptions.}\abs
Here we note that according to \cref{psi},
\bas
	\psi(s) \le \psi(1) + \int_1^s \frac{\Cpsi}{\sigma} d\sigma
	= \psi(1) + \Cpsi \ln s
	\qquad \mbox{for all } s\ge 1,
\eas
whence abbreviating $\tCpsi:=\Cpsi + \max_{s\in [0,1]} \psi(s)$ and $\ln_+(s):=\max\{0, \ln s\}$ for $s>0$,
we obtain the pointwise inequality
\be{psi1}
	\psi(s) \le \tCpsi \cdot (1+\ln_+ s )
	\qquad \mbox{for all } s>0
\ee
for $\psi$ itself, and from \cref{phi} it similarly follows that
\be{phi1}
	\phi(u) \le \tCphi \cdot (1+\ln_+ u )
	\qquad \mbox{for all } u>0
\ee
with $\tCphi:=\Cphi + \max_{s\in [0,1]} \phi(s)$.\abs

\noindent
{\bf Main results.} \quad
\begin{Theorem}\label{theo27}
  Let $\Omega\subset \R^2$ be a bounded convex domain with smooth boundary, assume that $\chi, \xi, \psi$ and $\phi$ 
  satisfy \cref{reg}, \cref{chi}, \cref{xi}, \cref{psi} and \cref{phi}.
  Then given any $(u_0,v_0,w_0,z_0)$ satisfying \cref{init} with some $\vartheta\in (0,1)$, one can find
  a uniquely determined globally defined classical solution 
  $(u,v,w,z) \in (C^{2,1}(\bom\times [0,\infty))^4$ of \cref{0} such that $u,v,w$ and $z$ are nonnegative.
\end{Theorem}
\noindent
{\bf Challenges and ideas.} \quad
A major challenge arising in the anaylsis of \cref{0} is linked to the circumstance that the haptotactic sensitivity coefficient
function $\xi$ therein depends, besides on the haptoattractant $w$ itself, on the quantity $v$ which can be viewed external
with regard to this taxis process. 
In the context of a standard approach toward reduction of \cref{0} to a system free of any explictly appearing 
haptotaxis mechanism (\cite{fontelos_friedman_hu}, \cite{friedman_tello}, \cite{pang_wang_JDE}), this gives rise to 
substantial additional complexity: Namely, the ambition to accordingly join the fluxes $\nabla u - u\xi(v,w)\nabla w$
suggests to introduce the new dependent variable
\be{a}
	a:=u e^{-\Xi(v,w)},
\ee
where in contrast to the case when $\xi=\xi(w)$, the function
\be{Xi}
	\Xi(v,w):=\int_0^w \xi(v,s) ds,
	\qquad (v,w) \in [0,\infty)^2,
\ee
now depends on $v$ as well. As a consequence, the evolution of $a$ is explicitly influenced by the time derivative $v_t$,
which unlike $w_t$ can apparently not be replaced with zero-order quantities, and which in view of the identity
$v_t=\Delta v - \na \cdot (v\na z)$ can hence rather be considered 
as an additional contribution of cross-diffusive and hence particularly delicate character (cf.~also \cref{0a} below).\abs
Thus led to investigating regularity properties of $v_t$, we shall lay a first focus on the chemotaxis-consumption subsystem
of \cref{0} made up by its second and fourth equation. 
Here we intend to make use of an essentially well-known structural property associated with the taxis-absorption interplay therein,
and hence build this part on an analysis of the functional
\be{F}
	{\cal F} := \io v\ln v +  \frac{1}{2} \io \frac{|\na z|^2}{z} .
\ee
Due to the presence of the production term $\phi(u)$ in the fourth equation from \cref{0}, however, the evolution of this functional
along trajectories (cf.~\eqref{44.3a}) is, inter alia, determined by a summand of the form
\be{F1}
	\io \phi'(u) \frac{\na z}{z} \cdot \na u,
\ee
an adequate estimation of which on the basis of \cref{phi} will be preceded by two arguments controlling the singular behavior
of $z$ and $u$ by providing spatio-temporal $L^2$ bounds for $\na \ln z$ and $\na \ln u$ (\cref{lem3} and \cref{lem4}).
This will pave the way for our derivation of an estimate for $v$ in $L\log L(\Omega)$ through an examination of ${\cal F}$
(\cref{lem44}), and in the considered two-dimensional setting this can thereafter be seen to actually entail estimates
for $v_t$ and $D^2 v$ in space-time $L^p$ norms for arbitrarily large $p>1$ (\cref{lem8}).
In \cref{lem11}, this will be found to be sufficient to imply $L^\infty$ bounds for the crucial quantity $a$,
whereupon an adequately arranged bootstrap procedure will provide estimates for all solution components with respect to the norm
in $C^{2+\theta}(\bom)$ for some $\theta\in (0,1)$, and thereby allow for global extensibility of a local-in-time solution
which exists according to a standard argument (\cref{lem25}, \cref{lem26} and  \cref{lem_loc}).
\subsection{An equivalent problem reformulation. Local existence}		
Following the strategy outlined above, we note that the substitution \cref{a}-\cref{Xi})
transforms \cref{0} to the problem
\be{0a}
	\left\{ \begin{array}{ll}
	a_t = e^{-\Xi(v,w)} \nabla \cdot \Big\{ e^{\Xi(v,w)} \nabla a \Big\}
	- e^{-\Xi(v,w)} \nabla \cdot \Big\{ ae^{\Xi(v,w)} \Big(\chi(v,w)-\Xi_v(v,w)\Big) \nabla v \Big\} \\[1mm]
	\hspace*{10mm}
	- a\Xi_v(v,w) v_t
	+ a\xi(v,w) \psi (aw e^{\Xi(v,w)}),
	\qquad & x\in \Omega, \ t>0, \\[1mm]
	v_t = \Delta v - \nabla \cdot (v\nabla z), 
	\qquad & x\in\Omega, \ t>0, \\[1mm]
	w_t = - \psi(awe^{\Xi(v,w)}),
	\qquad & x\in\Omega, \ t>0, \\[1mm]
	z_t = \Delta z - vz + \phi(ae^{\Xi(v,w)}),
	\qquad & x\in\Omega, \ t>0, \\[1mm]
	\frac{\partial a}{\partial\nu} 
	= \frac{\partial v}{\partial\nu} = \frac{\partial z}{\partial\nu} =0,
	\qquad & x\in\pO, \ t>0, \\[1mm]
	a(x,0)=u_0(x)e^{-\Xi(v_0(x),w_0(x))}, \ v(x,0)=v_0(x), \ w(x,0)=w_0(x), \ z(x,0)=z_0(x),
	\ & x\in\Omega,
	\end{array} \right.
\ee
which in the considered framework of classical solvability can indeed be seen to be equivalent to \cref{0}, because
the boundary condition $\frac{\partial w_0}{\partial\nu}=0$ on $\pO$ imposed in \cref{init} ensures that
classical solutions $(u,v,w,z)$ to \cref{0} satisfy $\frac{\partial w}{\partial\nu}=0$ throughout evolution.\abs
Based on the observation that due to the second equation in \cref{0a} the expression $v_t$ appearing in the first equation 
can actually be replaced with the additional cross-diffusion term $\Delta v-\na \cdot (v\na z)$,
we may resort to a standard reasoning in order to assert local existence 
of classical solutions, along with an appropriate extensibility criterion
and some essentially evident basic 
positivity, conservation and dissipation properties, in the following flavor.
\begin{Lemma}\label{lem_loc}
  Assume \cref{reg}, \cref{chi}, \cref{xi}, \cref{psi} and \cref{phi}, and 
  suppose that $u_0, v_0, w_0$ and $z_0$ are such that \cref{init} holds with some $\vartheta\in (0,1)$. 
  Then there exist $\tm\in (0,\infty]$ and a uniquely determined quadruple 
  of nonnegative functions $a,v,w$ and $z$ which belong to $C^{2,1}(\bom\times [0,\tm))$, and which are such that
  $(a,v,w,z)$ solves \cref{0a} in the classical sense in $\Omega\times (0,\tm)$, and that
  \be{ext}
	\mbox{if $\tm=\infty$, \quad then} \quad
	\limsup_{t\nearrow \tm} \Big\{ 
	\|a(\cdot,t)\|_{C^{2+\theta}(\bom)}
	+ \|v(\cdot,t)\|_{C^{2+\theta}(\bom)}
	+ \|w(\cdot,t)\|_{C^{2+\theta}(\bom)}
	+ \|z(\cdot,t)\|_{C^{2+\theta}(\bom)}
	\Big\}
	=\infty
  \ee
  for all $\theta\in (0,1)$.
  Moreover, this solution has the properties that $a>0, v>0$ and $z>0$ in $\bom\times (0,\tm)$, that with
  $u:=e^{\Xi(v,w)} a$ we have
  \be{massu}
	\io u(\cdot,t) = \io u_0
	\qquad \mbox{for all } t\in (0,\tm)
  \ee
  and
  \be{massv}
	\io v(\cdot,t)=\io v_0
	\qquad \mbox{for all } t\in (0,\tm),
  \ee
  and that
  \be{winfty}
	w(x,t) \le \|w_0\|_{L^\infty(\Omega)}
	\qquad \mbox{for all $x\in\Omega$ and } t\in (0,\tm).
  \ee
\end{Lemma}
\proof
  The statements concerning local existence and \cref{ext}
  can be verified by means of a quite well-established argument based on Banach's fixed point theorem, so that we
  may content ourselves with a reference to the reasoning detailed in \cite[Lemma 2.1]{taowin_JDE2014} for a closely
  related situation.\abs
  The claimed positivity properties of $a$, $v$ and $z$ thereafter result from \cref{init} and the strong maximum principle,
  whereas \cref{massu}, \cref{massv} and \cref{winfty} are obvious consequences of \cref{0a} and the nonnegativity
  of $\psi$ entailed by \cref{psi}.
\qed \\[-2ex]

\noindent
Throughout the sequel, without further explicit mentioning we shall suppose that 
\cref{reg}, \cref{chi}, \cref{xi}, \cref{psi}, \cref{phi} and \cref{init} hold,
and let $(a,v,w,z)$, $u$ and $\tm$ be as provided by \cref{lem_loc}, noting that then $(u,v,w,z)$ forms a classical
solution of \cref{0} in $\Omega\times (0,\tm)$.
\subsection{A bound for $v$ in $L\log L(\Omega)$}
\subsubsection{$L^2$ estimates for $\na \ln z$ and $\na \ln u$}
In this section concerned with the time evolution of the functional in \cref{F}, with a view to \cref{F1} we shall first focus
on the derivation of some estimates for $z$ and $u$ which, in some contrast to those explicitly appearing in \cref{ext}, mainly 
concentrate on a control of the behavior near {\em small} values of the respective quantities.\abs
As a preparation, let us first draw a consequence of a basic bound for $\phi(u)$, as implied by \cref{phi1} and \cref{massu},
when combined with nonpositivity of the summand $-vz$ in the forth equation from \cref{0}.
This, namely, already asserts an $L^\infty$ estimate for $z$ through a straightforward argument based on well-known parabolic
regularity theory:
\begin{Lemma} \label{lem2}
  If $\tm<\infty$, then there exists $C>0$ such that
  \be{zinfty}
	z(x,t) \le C
	\qquad \mbox{for all $x\in\Omega$ and } t\in (0,\tm).
  \ee
\end{Lemma} 
\proof
  We first observe that due to the validity of the elementary inequality $\ln s \le \frac{2\sqrt{s}}{e}$ for all $s>0$,
  from \cref{phi1} it follows that
  \bas
	\phi(u) \le \tCphi \cdot \Big(1+\frac{2}{e}\sqrt{u}\Big)
	\qquad \mbox{in } \Omega\times (0,\tm),
  \eas
  and that thus, by Young's inequality and \cref{massu},
  \be{2.2}
	\io \phi^2(u) \le 2\tCphi^2 \io \Big(1+\frac{4}{e^2} u\Big)
	= c_1:=2\tCphi^2 \cdot \bigg\{ |\Omega| + \frac{4}{e^2} \io u_0 \bigg\}
	\qquad \mbox{for all } t\in (0,\tm).
  \ee
  Now by nonnegativity of $v$ and $z$, from \cref{0} we obtain the one-sided inequality
  $z_t \le \Delta z + \phi(u)$ in $\Omega\times (0,\tm)$, which thanks to the comparison principle implies that
  \be{2.3}
	z(\cdot, t)
	\le e^{t\Delta} z_0 + \int_0^t e^{(t-s)\Delta} \phi(u(\cdot,s)) ds
	\quad \mbox{in } \Omega
	\qquad \mbox{for all } t\in (0,\tm),
  \ee
  where $(e^{\sigma\Delta})_{\sigma\ge 0}$ denotes the Neumann heat semigroup over $\Omega$.
  Here we employ a well-known smoothing estimate for $(e^{\sigma\Delta})_{\sigma\ge 0}$ (\cite[Lemma 1.4]{win_JDE}) to find $c_2>0$ such that
  \bas
	\|e^{\sigma \Delta} \varphi\|_{L^\infty(\Omega)} \le c_2 (1+\sigma^{-\frac{1}{2}}) \|\varphi\|_{L^2(\Omega)}
	\qquad \mbox{for all $\sigma>0$ and any } \varphi\in C^0(\bom),
  \eas
  so that again thanks to the comparison principle, \cref{2.3} together with \cref{2.2} entails that
  \bas
	z(\cdot,t)
	&\le& \|e^{t\Delta} z_0\|_{L^\infty(\Omega)}
	+ c_2 \int_0^t \Big(1+(t-s)^{-\frac{1}{2}}\Big) \|\phi(u(\cdot,s))\|_{L^2(\Omega)} ds \\
	&\le& \|z_0\|_{L^\infty(\Omega)}
	+ \sqrt{c_1} c_2 \int_0^t \Big(1+(t-s)^{-\frac{1}{2}}\Big) ds
	\quad \mbox{in } \Omega
	\qquad \mbox{for all } t\in (0,\tm).
  \eas
  Since
  \bas
	\int_0^t \Big(1+(t-s)^{-\frac{1}{2}}\Big) ds
	= t+2\sqrt{t} \le \tm + 2\sqrt{\tm}
	\qquad \mbox{for all } t\in (0,\tm),
  \eas
  in view of the presupposed finiteness of $\tm$ this establishes \cref{zinfty}.
\qed\\[-2ex]

\noindent
As a first consequence, through a fairly simple testing procedure the latter already entails an $L^2$ bound for $\na \ln z$
in the intended flavor:
\begin{Lemma} \label{lem3}
  If $\tm<\infty$, then 
  \be{3.1}
	\int_\tau^{\tm} \frac{|\nabla z|^2}{z^2} <\infty
	\qquad \mbox{for all } \tau\in (0,\tm).
  \ee
\end{Lemma} 
\proof
  Recalling that $z$ is positive in $\bom\times (0,\tm)$, we may multiply the fourth equation in \cref{0} by $-\frac{1}{z}$
  and integrate by parts to see that
  \bas
	-\frac{d}{dt} \io \ln z
	+ \io \frac{|\nabla z|^2}{z^2}
	= \io v - \io \frac{\phi(u)}{z}
	\le \io v_0
	\qquad \mbox{for all } t\in (0,\tm)
  \eas
  according to \cref{massv} and the nonnegativity of $\phi$ required in \cref{phi}. Therefore,
  \bas
	\int_\tau^t \io \frac{|\nabla z|^2}{z^2}
	&\le& \io \ln z(\cdot,t) - \io \ln z(\cdot,\tau) +(t-\tau) \io v_0 \\
	&\le& |\Omega| \cdot \ln \|z(\cdot,t)\|_{L^\infty(\Omega)} 
	- \io \ln z(\cdot,\tau){+T_{max}\io v_0}
	\quad \mbox{for all $\tau\in (0,\tm)$ and } t\in [\tau,\tm),
  \eas
  which implies \cref{3.1} since $\sup_{t\in (0,\tm)} \|z(\cdot,t)\|_{L^\infty(\Omega)}$ is finite by \cref{lem2},
  and since, again by positivity of $z$ in $\bom\times (0,\tm)$, also $-\io \ln z(\cdot,\tau)<\infty$ for all $\tau\in (0,\tm)$.
\qed\\[-2ex]

\noindent
Once again relying on \cref{lem2}, from the latter we immediately obtain an $L^2$ bound for $\na z$ also without any weight.
\begin{Corollary}\label{cor32}
  If $\tm<\infty$, then 
  \bas
	\int_\tau^{\tm} |\nabla z|^2 < \infty
	\qquad \mbox{for all } \tau\in (0,\tm).
  \eas
\end{Corollary}
\proof
  This can immediately be seen by combining \cref{lem3} with \cref{lem2}.
\qed\\[-2ex]

\noindent
This information now provides sufficient regularity of the taxis gradient in the second equation from \cref{0} to ensure
a singular estimate for $\na v$ quite similar to that from \cref{lem3}.
\begin{Lemma} \label{lem33}
  If $\tm<\infty$, then 
  \be{33.1}
	\int_\tau^{\tm} \frac{|\nabla v|^2}{v^2} <\infty
	\qquad \mbox{for all } \tau\in (0,\tm).
  \ee
\end{Lemma} 
\proof
  Thanks to the positivity of $v$ in $\bom\times (0,\tm)$, from the second equation in \cref{0} and Young's inequality we obtain that
  \bas
	-\frac{d}{dt} \io \ln v
	+ \io \frac{|\nabla v|^2}{v^2}
	= \io \frac{\na v}{v} \cdot \na z
	\le \frac{1}{2} \io \frac{|\na v|^2}{v^2} + \frac{1}{2} \io |\na z|^2
	\qquad \mbox{for all } t\in (0,\tm),
  \eas
  and that thus, given $\tau\in (0,\tm)$ we have 
  \bas
	\int_\tau^t \io \frac{|\na v|^2}{v^2} 
	\le 2 \io \ln v(\cdot,t) - 2 \io \ln v(\cdot,\tau)
	+ \int_\tau^t \io |\na z|^2
	\qquad \mbox{for all } t\in [\tau,\tm).
  \eas
  Since here $\sup_{t\in (\tau,\tm)} \int_\tau^t \io |\na z|^2 <\infty$ by  \cref{cor32} and
  $\io \ln v(\cdot,\tau)>-\infty$ by positivity of $v(\cdot,\tau)$ in $\bom$, \cref{33.1} results upon observing that
  \bas
	\io \ln v(\cdot,t) \le \io v(\cdot,t) = \io v_0
	\qquad \mbox{for all } t\in (0,\tm)
  \eas
  due to \cref{massv} and the fact that $\ln s \le s$ for all $s>0$.
\qed\\[-2ex]

\noindent
Now making essential use of our overall assumption \cref{psi}, by analyzing the time evolution of $-\io \ln u + \io |\na w|^2$
we can turn the information from \cref{lem33} into the following second main result of this subsection:
\begin{Lemma} \label{lem4}
  If $\tm<\infty$, then 
  \be{4.1}
	\int_\tau^{\tm} \frac{|\nabla u|^2}{u^2} <\infty
	\qquad \mbox{for all } \tau\in (0,\tm).
  \ee
\end{Lemma} 
\proof
  Since $u$ is positive in $\bom\times (0,\tm)$, on the basis of the first equation in \cref{0}
  we may integrate by parts and use Young's inequality along with \cref{chi}
  and \cref{xi} to see that
  \bea{4.2}
	- \frac{d}{dt} \io \ln u + \io \frac{|\na u|^2}{u^2}
	&=& \io \chi(v,w) \frac{\na u}{u} \cdot \na v
	+ \io \xi(v,w) \frac{\na u}{u} \cdot\na w \nn\\
	&\le& \frac{1}{2} \io \frac{|\na u|^2}{u^2}
	+ \io \chi^2(v,w) |\na v|^2
	+ \io \xi^2(v,w) |\na w|^2 \nn\\
	&\le& \frac{1}{2} \io \frac{|\na u|^2}{u^2}
	+ \Cchi^2 \io \frac{|\na v|^2}{v^2}
	+ \Cxi^2 \io |\na w|^2
	\qquad \mbox{for all } t\in (0,\tm).
  \eea
  To appropriately compensate the rightmost summand herein, we furthermore utilize the third equation from \cref{0}
  to compute
  \bas
	\frac{d}{dt} \io |\na w|^2
	&=& - 2 \io \na w \cdot \na \psi(uw) \\
	&=& - 2 \io u\psi'(uw) |\na w|^2
	- 2 \io w\psi'(uw) \na u\cdot\na w 
	\qquad \mbox{for all } t\in (0,\tm),
  \eas
  so that since \cref{psi} asserts that $u\psi'(uw) \ge 0$ and $0\le uw\psi'(uw) \le \Cpsi$ in $\Omega\times (0,\tm)$, by means of
  Young's inequality we obtain that
  \bas
	\frac{d}{dt} \io |\na w|^2
	&\le& - 2\io uw\psi'(uw) \frac{\na u}{u} \cdot \na w \\
	&\le& \frac{1}{4} \io \frac{|\na u|^2}{u^2}
	+ 4\io \Big\{ uw\psi'(uw)\Big\}^2 |\na w|^2 \\
	&\le& \frac{1}{4} \io \frac{|\na u|^2}{u^2}
	+ 4\Cpsi^2 \io |\na w|^2
	\qquad \mbox{for all } t\in (0,\tm).
  \eas
  In conjunction with \cref{4.2}, this shows that writing $c_1:=\Cxi^2 + 4\Cpsi^2$ we have
  \bas
	\frac{d}{dt} \bigg\{ - \io \ln u + \io |\na w|^2 \bigg\}
	+ \frac{1}{4} \io \frac{|\na u|^2}{u^2}
	\le c_1 \io |\na w|^2
	+ \Cchi^2 \io \frac{|\na v|^2}{v^2}
	\qquad \mbox{for all } t\in (0,\tm),
  \eas
  so that since $\io \ln u \le \io u = c_2:=\io u_0$ for all $t\in (0,\tm)$ by \cref{massu}, it follows that
  $y(t):=-\io \ln u(\cdot,t) + \io |\na w(\cdot,t)|^2$, $g(t):=\frac{1}{4} \io \frac{|\na u(\cdot,t)|^2}{u^2(\cdot,t)}$
  and $h(t):=\Cchi^2 \io \frac{|\na v(\cdot,t)|^2}{v^2(\cdot,t)} + c_1 c_2$, $t\in (0,\tm)$, satisfy
  \be{4.22}
	\io |\na w|^2
	= y(t) + \io \ln u
	\le y(t) + c_2
	\qquad \mbox{for all } t\in (0,\tm)
  \ee
  and hence
  \be{4.3}
	y'(t) + g(t) \le c_1 \cdot (y(t)+c_2) 
	+ \Cchi^2 \io \frac{|\na v|^2}{v^2}
	\le c_1 y(t) + h(t)
	\qquad \mbox{for all } t\in (0,\tm).
  \ee
  For fixed $\tau\in (0,\tm)$, upon simply estimating $g\ge 0$ and noting that
  \bas
	c_3:=y(\tau) 
	\quad \mbox{and} \quad
	c_4:=\int_\tau^{\tm} h(s) ds
  \eas
  are both finite by positivity of $u(\cdot,\tau)$ and $\bom$ and \cref{lem33}, from this we firstly infer through an ODE
  comparison argument that
  \bas
	y(t)
	&\le& y(\tau) e^{c_1 \cdot (t-\tau)} + \int_\tau^t e^{c_1\cdot (t-s)} h(s) ds \\
	&\le& c_5:=c_3 e^{c_1 \tm} + c_4 e^{c_1 \tm}
	\qquad \mbox{for all } t\in [\tau,\tm).
  \eas
  Thereupon, once more going back to \cref{4.3} we infer that, again by \cref{4.22},
  \bas
	\int_\tau^t g(s) ds
	&\le& y(\tau) -y(t) + c_1 \int_\tau^t y(s) ds + \int_\tau^t h(s)ds \\
	&\le& c_5 + \bigg\{ c_2-\io |\na w|^2 \bigg\} + c_1 c_5 \cdot (t-\tau) + c_4 \\
	&\le& c_5 + c_2 + c_1 c_5 \tm + c_4
	\qquad \mbox{for all } t\in [\tau,\tm),
  \eas
  which by definition of $g$ implies \cref{4.1}.
\qed
\subsubsection{The chemotaxis-consumption interaction in \cref{0}. A bound for $v$ in $L\log L$}\label{subsec:bound-v}
%
%
%
%
%
%
%
We are now prepared for our analysis of the chemotaxis-consumption subsystem of \cref{0} at the fundamental level of the 
functional ${\cal F}$ in \cref{F}, 
 well-known 
to play the role of a genuine energy in the accordingly unperturbed case
when $\phi\equiv 0$ in convex domains (\cite{taowin_consumption}).
In particular, \cref{lem3} and \cref{lem4} will enable us to 
appropriately control the respective additional contribution due to the coupling to the crucial quantity $u$, as
foreshadowed in \cref{F1}, and to thereby derive the following key estimate for $v$.
\begin{Lemma} \label{lem44}
  If $\tm<\infty$, then there exists $C>0$ such that
  \be{44.1}
	\io v(\cdot,t) \big| \ln v(\cdot,t)\big| \le C
	\qquad \mbox{for all } t\in (0,\tm).
  \ee
\end{Lemma} 
\proof
  Again relying on the positivity of $v$ and $z$ in $\bom\times (0,\tm)$, by means of several integrations by parts in the second and
  fourth equation from \cref{0} we compute
  \be{44.2}
	\frac{d}{dt} \io v\ln v
	+ \io \frac{|\na v|^2}{v} = \io \na v\cdot\na z
	\qquad \mbox{for all } t\in (0,\tm)
  \ee
  and
  \bea{44.3a}
	\frac{1}{2} \frac{d}{dt} \io \frac{|\na z|^2}{z}
	&=& \io \frac{\na z}{z} \cdot\na \Big\{ \Delta z - vz + \phi(u) \Big\} 
	- \frac{1}{2} \io \frac{|\na z|^2}{z^2} \cdot \Big\{ \Delta z - vz + \phi(u)\Big\} \nn\\
	&=& - \io z|D^2 \ln z|^2
	+ \frac{1}{2} \int_{\pO} \frac{1}{z} \cdot \frac{\partial |\na z|^2}{\partial\nu} \nn\\
	& & - \io \na v\cdot\na z
	- \frac{1}{2} \io \frac{v}{z} |\na z|^2 \nn\\
	& & + \io \phi'(u) \frac{\na z}{z} \cdot\na u
	- \frac{1}{2} \io \frac{|\na z|^2}{z^2} \phi(u)
	\qquad \mbox{for all } t\in (0,\tm),
  \eea
  where we have used that $\na z\cdot\na\Delta z = \frac{1}{2} \Delta |\na z|^2 - |D^2 z|^2$ in $\Omega\times (0,\tm)$,
  and that
  \bas
	\frac{1}{2} \io \frac{1}{z^2} \na z\cdot\na |\na z|^2
	- \frac{1}{z} |D^2 z|^2
	- \frac{1}{2} \io \frac{1}{z^2} |\na z|^2 \Delta z
	= - \io z|D^2 ln z|^2
	\qquad \mbox{for all } t\in (0,\tm)
  \eas
  (cf.~\cite[Lemma 3.2]{win_CPDE2012} for a detailed derivation of the latter identity).
  Here since $\frac{\partial |\na z|^2}{\partial\nu} \le 0$ on $\pO\times (0,\tm)$ by convexity of $\Omega$ (\cite{lions_ARMA}),
  in view of the nonnegativity of $\phi$ postulated in \cref{phi} it follows that the second, fourth and sixth summands on the
  right of \eqref{44.3a} are all nonpositive, whereas in the crucial second last integral therein we make use of the upper bound
  for $\phi'$ in \cref{phi} to see that thanks to Young's inequality,
  \bas
	\io \phi'(u) \frac{\na z}{z} \cdot\na u
	&\le& \io \frac{|\na z|^2}{z^2} 
	+ \frac{1}{4} \io (\phi'(u))^2 |\na u|^2 \\
	&\le& \io \frac{|\na z|^2}{z^2} 
	+ \frac{\Cphi^2}{4} \io \frac{|\na u|^2}{u^2}
	\qquad \mbox{for all } t\in (0,\tm).
  \eas
  When adding \eqref{44.3a} to \cref{44.2} and taking benefit from a favorable cancellation thereby induced on the right-hand side 
  of \cref{44.2}, we thus obtain that
  \bas
	& & \hspace*{-20mm}
	\frac{d}{dt} \bigg\{ \io v\ln v + \frac{1}{2} \io \frac{|\na z|^2}{z} \bigg\}
	+ \io  \frac{|\na v|^2}{v} + \io z|D^2 \ln z|^2 \\
	&\le& \io \frac{|\na z|^2}{z^2} + \frac{\Cphi^2}{4} \io \frac{|\na u|^2}{u^2}
	\qquad \mbox{for all } t\in (0,\tm),
  \eas
  which upon an integration in time implies that, writing $\tau:=\frac{1}{2}\tm$, for all $t\in [\tau,\tm)$, we have
  \bas
	\io v(\cdot,t) \ln v(\cdot,t)
	+ \frac{1}{2} \io \frac{|\na z(\cdot,t)|^2}{z(\cdot,t)}
	&\le& c_1:=\io v(\cdot,\tau)\big|\ln v(\cdot,\tau)\big|
	+ \frac{1}{2} \io \frac{|\na z(\cdot,\tau)|^2}{z(\cdot,\tau)} \\
	& & + \int_\tau^{\tm} \io \frac{|\na z|^2}{z^2}
	+ \frac{\Cphi^2}{4} \int_\tau^{\tm} \io \frac{|\na u|^2}{u^2},
  \eas
  with $c_1$ being finite and positive due to \cref{lem3} and \cref{lem4}. 
  Since $\io v|\ln v| = \io v\ln v - 2\int_{\{v<1\}} v\ln v \le \io v\ln v + \frac{2|\Omega|}{e}$ for all $t\in (0,\tm)$
  due to the fact that $s\ln s \ge -\frac{1}{e}$ for all $s>0$, this entails that
  \bas
	\io v(\cdot,t) \big| \ln v(\cdot,t)\big| \le c_1 + \frac{2|\Omega|}{e}
	\qquad \mbox{for all } t\in [\tau,\tm),
  \eas
  and that thus \cref{44.1} holds if we let 
  $C:=\max \Big\{ \sup_{t\in (0,\tau)} \io v(\cdot,t) |\ln v(\cdot,t)| \, , \, c_1+\frac{2|\Omega|}{e}\Big\}$.	
\qed
\subsection{Maximal Sobolev regularity results for $z$ and $v$. $L^p$ bounds for $v_t$}
In this section we shall use the outcome of \cref{lem44}, and again that of \cref{lem2}, as a starting point
for a bootstrap procedure applied to the second and fourth equations in \cref{0}. 
Thanks to the availability of suitable embedding properties,
in the considered two-dimensional framework this initial regularity information will actually be seen to imply
estimates for both $(z,\na z,D^2 z,z_t)$ and $(v,\na v,D^2 v, v_t)$ in arbitrary space-time $L^p$ spaces.\abs
Our first step in this direction establishes the following by analyzing the evolution of $\io v^2 + \frac{1}{2} \io |\na z|^4$.
\begin{Lemma} \label{lem5}
  If $\tm<\infty$, then there exists $C>0$ such that
  \be{5.1}
	\io |\na z(\cdot,t)|^4 \le C
	\qquad \mbox{for all } t\in (0,\tm).
  \ee
\end{Lemma} 
\proof
  Following the strategy from \cite[Lemma 3.3]{BBTW}, we use the second and fourth equations in \cref{0} to see that again since
  $\na z\cdot\na \Delta z=\frac{1}{2}\Delta |\na z|^2 - |D^2 z|^2$ in $\Omega\times (0,\tm)$ and 
  $\frac{\partial |\na z|^2}{\partial\nu}\le 0$ on $\pO\times (0,\tm)$, according to Young's inequality we have
  \be{5.2}
	\frac{1}{2} \frac{d}{dt} \io v^2 + \io |\na v|^2
	= \io v\na v\cdot\na z 
	\le \frac{1}{2} \io |\na v|^2
	+ \frac{1}{2} \io v^2 |\na z|^2
	\qquad \mbox{for all } t\in (0,\tm)
  \ee
  and
  \bea{5.3-eq}
	\frac{1}{4} \frac{d}{dt} \io |\na z|^4
	&=& \io |\na z|^2 \na z\cdot \na \Big\{ \Delta z - vz + \phi(u) \Big\} \nn\\
	&=& - \frac{1}{2} \io \Big|\na |\na z|^2 \Big|^2
	+ \frac{1}{2} \int_{\pO} |\na z|^2 \frac{\partial |\na z|^2}{\partial\nu}
	- \io |\na z|^2 |D^2 z|^2 \nn\\
	& & + \io \phi'(u) |\na z|^2 \na u\cdot\na z +\io vz\na \cdot (|\na z|^2\na z)\nn \\
	&\le& - \io |\na z|^2 |D^2 z|^2
	+ (2+\sqrt{2}) c_1 \io v|\na z|^2 |D^2 z| \nn\\
	& & + \io \phi'(u) |\na z|^2 \na u\cdot\na z \nn\\
	&\le& - \frac{1}{2} \io |\na z|^2 |D^2 z|^2
	+ \frac{(2+\sqrt{2})^2 c_1^2}{2} \io v^2 |\na z|^2 \nn\\
	& & + \io \phi'(u) |\na z|^2 \na u\cdot\na z
	\qquad \mbox{for all } t\in (0,\tm),
  \eea
  where $c_1:=\|z\|_{L^\infty(\Omega\times (0,\tm))}$ is finite by \cref{zinfty}, 
  and where we have used that for each $\gamma\ge 1$,
  \bea{5.4}
	\Big| \na \cdot \Big( |\na\varphi|^{2\gamma} \na\varphi\Big) \Big|
	&=& \Big| 2\gamma |\na\varphi|^{2\gamma-2} \na\varphi\cdot (D^2\varphi\cdot\na\varphi) 
		+ |\na\varphi|^{2\gamma} \Delta\varphi \Big| \nn\\
	&\le& (2\gamma+\sqrt{2}) |\na\varphi|^{2\gamma} |D^2\varphi|	
	\quad \mbox{in } \Omega
	\qquad \mbox{for all } \varphi\in C^2(\bom).
  \eea
  A combination of \eqref{5.2} with \eqref{5.3-eq} thus shows that
  \bea{5.5}
	& & \hspace*{-20mm}
	\frac{d}{dt} \bigg\{ \io v^2 + \frac{1}{2} \io |\na z|^4 \bigg\}
	+ \io |\na v|^2 + \io |\na z|^2 |D^2 z|^2 \nn\\
	&\le& c_2 \io v^2 |\na z|^2
	+ 2\io \phi'(u) |\na z|^2 \na u\cdot\na z
	\qquad \mbox{for all } t\in (0,\tm)
  \eea
  with $c_2:=1+(2+\sqrt{2})^2 c_1^2$, and to proceed from this we note that again thanks to \eqref{5.4}, for arbitrary 
  $\varphi\in C^2(\bom)$ satisfying $\frac{\partial\varphi}{\partial\nu}=0$ on $\pO$ we have
  \bas
	\io |\na\varphi|^6
	&=& \io |\na\varphi|^4 \na\varphi\cdot\na\varphi \\
	&=& - \io \varphi \na\cdot (|\na\varphi|^4\na\varphi) \\
	&\le& (4+\sqrt{2}) \io |\varphi| \cdot |\na\varphi|^4 \cdot |D^2\varphi| \\
	&\le& (4+\sqrt{2}) \|\varphi\|_{L^\infty(\Omega)} \cdot \bigg\{ \io |\na\varphi|^6 \bigg\}^\frac{1}{2} \cdot
		\bigg\{ \io |\na \varphi|^2 |D^2\varphi|^2 \bigg\}^\frac{1}{2}
  \eas
  as a consequence of the Cauchy-Schwarz inequality, and that thus
  \bas
	\io |\na\varphi|^6
	\le c_3 \|\varphi\|_{L^\infty(\Omega)}^2 \io |\na\varphi|^2 |D^2\varphi|^2
  \eas
  holds for any such $\varphi$ if we let $c_3:=(4+\sqrt{2})^2$.
  Therefore, namely, by definition of $c_1$ we obtain that with $c_4:=c_1^2 c_3$,
  \be{5.6}
	\io |\na z|^6  \le c_4 \io |\na z|^2 |D^2 z|^2
	\qquad \mbox{for all } t\in (0,\tm),
  \ee
  and in line with this we apply Young's inequality on the right-hand side of \eqref{5.5} in such a way that
  \bas
	c_2 \io v^2 |\na z|^2
	\le \frac{1}{2c_4} \io |\na z|^6
	+ \frac{c_2^2 c_4}{2} \io v^3
	\qquad \mbox{for all } t\in (0,\tm)
  \eas
  and
  \bas
	2 \io \phi'(u) |\na z|^2 \na u \cdot\na z
	\le \frac{1}{2c_4} \io |\na z|^6 
	+ 2c_4 \io (\phi'(u))^2 |\na u|^2
	\qquad \mbox{for all } t\in (0,\tm),
  \eas
  so that recalling \eqref{phi} we infer that due to \eqref{5.6},
  \bea{5.7}
	\hspace*{-10mm}
	c_2 \io v^2 |\na z|^2
	+ 2\io \phi'(u) |\na z|^2 \na u \cdot\na z
	&\le& \frac{1}{c_4} \io |\na z|^6
	+ \frac{c_2^2 c_4}{2} \io v^3
	+ 2c_4 \Cphi^2 \io \frac{|\na u|^2}{u^2} \nn\\
	&\le& \io |\na z|^2 |D^2 z|^2
	+ \frac{c_2^2 c_4}{2} \io v^3
	+ 2c_4 \Cphi^2 \io \frac{|\na u|^2}{u^2} 
  \eea
  for all $t\in (0,\tm)$.
  Here the second last summand can be controlled by relying on the $L\log L$ estimate provided by \cref{lem44}, which
  in conjunction with a well-known variant of the Gagliardo-Nirenberg inequality (\cite{BHN}), namely, entails the existence of $c_5>0$
  such that 
  \bas
	\frac{c_2^2 c_4}{2} \io v^3
	\le \io |\na v|^2 + c_5
	\qquad \mbox{for all } t\in (0,\tm).
  \eas
  Inserted into \eqref{5.7}, this shows that \eqref{5.5} implies the inequality
  \bas
	\frac{d}{dt} \bigg\{ \io v^2 + \frac{1}{2} \io |\na z|^4 \bigg\}
	\le c_5 + 2c_4 \Cphi^2 \io \frac{|\na u|^2}{u^2}
	\qquad \mbox{for all } t\in (0,\tm),
  \eas
  from wich by integration over $[\tau,t]$ with $\tau:=\frac{1}{2}\tm$ and $t\in [\tau,\tm)$ it particularly follows that
  for any such $t$,
  \bas
	\frac{1}{2} \io |\na z(\cdot,t)|^4
	\le c_6:= \io v^2(\cdot,\tau)
	+ \frac{1}{2} \io |\na z(\cdot,\tau)|^2
	+ c_5 \cdot (\tm-\tau)
	+ 2c_4 \Cphi^2 \int_\tau^{\tm} \io \frac{|\na u|^2}{u^2},
  \eas
  with finiteness of $c_6$ asserted by \cref{lem4}.
  Since also $c_7:=\sup_{t\in (0,\tau)} \io |\na z(\cdot,t)|^4$ is finite by \cref{lem_loc}, this entails \eqref{5.1}
  with $C:=\max\{ c_7, 2c_6\}$.
\qed\\[-2ex]

\noindent
Using that the integrability exponent in \eqref{5.1} exceeds the considered spatial dimension, 
by means of an essentially well-established argument based on semigroup estimates 
we can turn this into an $L^\infty$ bound for $v$.
\begin{Lemma} \label{lem6}
  If $\tm<\infty$, then there exists $C>0$ such that
  \be{vinfty}
	v(x,t) \le C
	\qquad \mbox{for all $x\in\Omega$ and } t\in (0,\tm).
  \ee
\end{Lemma} 
\proof
  We abbreviate $c_1:=\io v_0$ and take $c_2>0$ such that in accordance with \cref{lem5} we have
  $\|\na z(\cdot,t)\|_{L^4(\Omega)} \le c_2$ for all $t\in (0,\tm)$, and to estimate the numbers
  \bas
	M(T):=\sup_{t\in (0,\tm)} \|v(\cdot,t)\|_{L^\infty(\Omega)},
	\qquad T\in (0,\tm),
  \eas
  we recall a well-known smoothing property of the Neumann heat semigroup (\cite{FIWY}) to fix $c_3>0$ such that 
  for all $\sigma>0$,
  \bas
	\|e^{\sigma\Delta} \na \cdot\varphi\|_{L^\infty(\Omega)}
	\le c_3 \cdot (1-\sigma^{-\frac{5}{6}}) \|\varphi\|_{L^3(\Omega)}
	\qquad \mbox{for all $\varphi\in C^1(\Omega;\R^2)$ such that $\varphi\cdot\nu=0$ on $\pO$.}
  \eas
  On the basis of a Duhamel representation, we then see that for each $T\in (0,\tm)$, due to the maximum principle we have
  \bas
	\|v(\cdot,t)\|_{L^\infty(\Omega)}
	&=& \bigg\| e^{t\Delta} v_0 
	- \int_0^t e^{(t-s)\Delta} \na \cdot \Big( v(\cdot,s)\na z(\cdot,s)\Big) ds \bigg\|_{L^\infty(\Omega)} \\
	&\le& \|v_0\|_{L^\infty(\Omega)}
	+ c_3 \int_0^t \Big(1+(t-s)^{-\frac{5}{6}}\Big) \|v(\cdot,s)\na z(\cdot,s)\|_{L^3(\Omega)} ds \\
	&\le& \|v_0\|_{L^\infty(\Omega)}
	+ c_3 \int_0^t \Big(1+(t-s)^{-\frac{5}{6}}\Big) \|v(\cdot,s)\|_{L^\infty(\Omega)}^\frac{11}{12}
		\|v(\cdot,s)\|_{L^1(\Omega)}^\frac{1}{12} \|\na z(\cdot,s)\|_{L^4(\Omega)} ds \\
	&\le& c_4 + c_4 M^\frac{11}{12}(T)
	\qquad \mbox{for all } t\in (0,T)
  \eas
  with 
  $c_4:=\max \Big\{ \|v_0\|_{L^\infty(\Omega)} \, , \, c_1^\frac{1}{12} c_2 c_3 \int_0^{\tm} (1+\sigma^{-\frac{5}{6}})d\sigma \Big\}$.
  Consequently, $M(T) \le c_4 + c_4 M^\frac{11}{12}(T)$ and hence $M(T) \le \max \{1 \, , \, (2c_4)^{12} \}$ for all
  $T\in (0,\tm)$, as intended.
\qed\\[-2ex]

\noindent
The latter now provides sufficient regularity information on the expression $z_t-\Delta z$ to warrant
accessibility to maximal Sobolev regularity theory:
\begin{Lemma} \label{lem7}
  If $\tm<\infty$, then for all $p>1$,
  \be{7.1}
	\int_0^{\tm} \Big\{ \|z(\cdot,t)\|_{W^{2,p}(\Omega)}^p + \|z_t(\cdot,t)\|_{L^p(\Omega)}^p \Big\} dt <\infty.
  \ee
\end{Lemma} 
\proof
In view of the logarithmic bound in \eqref{phi1} and of \eqref{massu} it readily follows from \eqref{zinfty} and \eqref{lem6}
  that $h:=-vz+\phi(u)$ belongs to $L^p(\Omega\times (0,\tm))$ for each finite $p>1$. The claim is an immediate consequence
  of a standard result on maximal Sobolev regularity in the Neumann problem for the inhomogeneous linear heat equation
  $z_t=\Delta z+h(x,t)$ with initial data satisfying the regularity and compatibility conditions stated in \eqref{init} (\cite{giga_sohr}).
\qed\\[-2ex]

\noindent
Along with \eqref{vinfty}, this in turn facilitates an application of the same token to the more delicate second equation
in \eqref{0}:
\begin{Lemma} \label{lem8}
  If $\tm<\infty$, then for all $p>1$,
  \be{8.1}
	\int_0^{\tm} \Big\{ \|v(\cdot,t)\|_{W^{2,p}(\Omega)}^p + \|v_t(\cdot,t)\|_{L^p(\Omega)}^p \Big\} dt <\infty.
  \ee
\end{Lemma} 
\proof
  We fix $p>1$ and then infer from \cref{lem7} and \cref{lem6} that there exist $c_1(p)>0,\ c_2(p)>0$ and $c_3>0$ such that
  \be{8.3}
	\int_0^T \io |\na z|^{2p} \le c_1(p)
	\qquad \mbox{for all } T\in (0,\tm)
  \ee
  and
  \be{8.4}
	\int_0^T \io |\Delta z|^p \le c_2(p)
	\qquad \mbox{for all } T\in (0,\tm)
  \ee
  as well as
  \be{8.5}
	v(x,t) \le c_3
	\qquad \mbox{for all $x\in\Omega$ and } t\in (0,\tm),
  \ee
  while maximal Sobolev regularity theory (\cite{giga_sohr}) provides $c_4(p)>0$ with the property
  that whenever $h\in C^0(\bom\times [0,\tm))$ and $\varphi\in C^{2,1}(\bom\times [0,\tm))$ are such that
  \bas
	\left\{ \begin{array}{ll}
	\varphi_t = \Delta\varphi + h(x,t),
	\qquad & x\in\Omega, \ t\in (0,\tm), \\[1mm]
	\frac{\partial\varphi}{\partial\nu}=0,
	& x\in\pO, \ t\in (0,\tm), \\[1mm]
	\varphi(x,0)=v_0(x),
	& x\in\Omega,
	\end{array} \right.
  \eas
  we have
  \be{8.6}
	\int_0^T \Big\{ \|\varphi(\cdot,t)\|_{W^{2,p}(\Omega)}^p + \|\varphi_t(\cdot,t)\|_{L^p(\Omega)}^p \Big\} dt
	\le c_4(p) + c_4(p) \int_0^T \io |h|^p
	\qquad \mbox{for all } T\in (0,\tm).
  \ee
  Furthermore, using the Gagliardo-Nirenberg inequality to choose $c_5(p)>0$ such that
  \be{8.7}
	\io |\nabla\varphi|^{2p} 
	\le c_5(p) \|\varphi\|_{W^{2,p}(\Omega)}^p \|\varphi\|_{L^\infty(\Omega)}^p
	\qquad \mbox{for all } \varphi\in W^{2,p}(\Omega),
  \ee
  we infer from \eqref{8.6} that, using the second equation in \eqref{0}, along with \eqref{8.4}, \eqref{8.5}, the Cauchy-Schwarz inequality, \eqref{8.3} and \eqref{8.7},
  \bas
	& & \hs{-20mm}
	\int_0^T \Big\| \|v(\cdot,t)\|_{W^{2,p}(\Omega)}^p + \|v_t(\cdot,t)\|_{L^p(\Omega)}^p \Big\} dt \\
	&\le& c_4(p) + c_4(p) \int_0^T \io \Big| - \na \cdot (v\na z) \Big|^p \\
	&\le& c_4(p) + 2^{p-1} c_4(p) \int_0^T \io |v\Delta z|^p
	+ 2^{p-1} c_4(p) \int_0^T \io |\na v\cdot\na z|^p \\
	&\le& c_4(p) + 2^{p-1} c_4(p) \|v\|_{L^\infty(\Omega\times (0,T))} \int_0^T \io |\Delta z|^p \\
	& & + 2^{p-1} c_4(p) \cdot \bigg\{ \int_0^T \io |\na v|^{2p} \bigg\}^\frac{1}{2}
		\cdot \bigg\{ \io |\na z|^{2p} \bigg\}^\frac{1}{2} \\
	&\le& c_4(p) + 2^{p-1} c_4(p) \|v\|_{L^\infty(\Omega\times (0,T))} \int_0^T \io |\Delta z|^p \\
	& & + 2^{p-1} c_4(p) c_5^\frac{1}{2}(p) \|v\|_{L^\infty(\Omega\times (0,T))}^\frac{p}{2}
	\cdot \bigg\{ \int_0^T \|v(\cdot,t)\|_{W^{2,p}(\Omega)}^p dt \bigg\}^\frac{1}{2}
	\cdot \bigg\{ \int_0^T \io |\na z|^{2p} \bigg\}^\frac{1}{2} \\
	&\le& c_6(p) + c_6(p) \cdot \bigg\{ \int_0^T \|v(\cdot,t)\|_{W^{2,p}(\Omega)}^p dt \bigg\}^\frac{1}{2}
	\qquad \mbox{for all } T\in (0,\tm)
  \eas
  with 
  $c_6(p):=\max \Big\{ c_4(p) + 2^{p-1} c_2(p) c_3 c_4(p) \, , \, 2^{p-1} c_1^\frac{1}{2}(p) c_3^\frac{p}{2} c_4(p) c_5^\frac{1}{2}(p)
  \Big\}$. 
  Since  
  \bas
	c_6(p) \cdot \bigg\{ \int_0^T \|v(\cdot,t)\|_{W^{2,p}(\Omega)}^p dt \bigg\}^\frac{1}{2}
	\le \frac{1}{2} \int_0^T \|v(\cdot,t)\|_{W^{2,p}(\Omega)}^p dt + \frac{c_6^2(p)}{2}
	\qquad \mbox{for all } T\in (0,\tm)
  \eas
  by Young's inequality, this implies \eqref{8.1}.
\qed\\[-2ex]

\noindent
For later reference, let us also explicitly state here the following immediate consequence concerning H\"older regularity of $v$ and
$\na v$.
\begin{Lemma} \label{lem9}
  If $\tm<\infty$, then there exists $\theta\in (0,1)$ such that
  \be{9.1}
	v\in C^{1+\theta,\theta}(\bom\times [0,\tm]).
  \ee
\end{Lemma} 
\proof  
  Due to a well-known embedding property (\cite{amann}), this directly results from an application of \cref{lem8} to
  some suitably large $p>1$.
\qed
\subsection{Pointwise boundedness of $a$}
Next, approaching the core of our analysis, we shall use the regularity information gained so far, and especially the 
$L^p$ bounds for $v_t$ included in \cref{lem8}, to derive an $L^\infty$ estimate for the apparently most crucial
solution component $a$. This will be achieved by means of an appropriately designed $L^p$ iteration of
Moser type on the basis of the first equation in \eqref{0a}, noting that, inter alia, due to the explicit presence of $v_t$ therein,
classical results in this regard (\cite{alikakos}, \cite{taowin_subcritical}) apparently do not apply directly to the current
situation.
\begin{Lemma} \label{lem11}
  If $\tm<\infty$, then there exists $C>0$ such that
  \be{11.1}
	a(x,t) \le C
	\qquad \mbox{for all $x\in\Omega$ and } t\in (0,\tm).
  \ee
\end{Lemma} 
\proof
  For integers $k\ge 0$, we let $p_k:=2^k$ and
  \be{11.2}
	M_k(T):=\max \bigg\{ 1 \, , \, \sup_{t\in (0,T)} \io a^{p_k}(\cdot,t) \bigg\},
	\qquad T\in (0,\tm).
  \ee  
  To appropriately control the latter quantities we note that according to  \eqref{winfty}, \eqref{xi}, \eqref{Xi},
  \eqref{chi} and \eqref{psi1} we can find positive constants $c_i, i\in\{1,...,6\}$, such that
  \be{11.3}
	0 \le w \le c_1,
	\quad
	|\xi(v,w)| \le c_2,
	\quad 
	|\Xi(v,w)| \le c_3,
	\quad
	|\Xi_v(v,w)| \le c_4
	\quad \mbox{and} \quad
	|\chi(v,w)| \le c_5
	\quad \mbox{in } \Omega\times (0,\tm)
  \ee
  as well as
  \be{11.33}
	\psi(s) \le c_6 s^\frac{1}{3} + c_6
	\qquad \mbox{for all } s\ge 0.
  \ee
  Moreover, \cref{lem9} and \cref{lem8} allow us to fix $c_7>0$ and $c_8>0$ fulfilling
  \be{11.4}
	|\na v| \le c_7 
	\qquad \mbox{in } \Omega\times (0,\tm)
  \ee
  and
  \be{11.5}
	\int_0^T \io v_t^4 \le c_8
	\qquad \mbox{for all } T\in (0,\tm).
  \ee
  Therefore, namely, we particularly know that due to \eqref{a} and \eqref{massu},
  \bas
	\io a = \io u e^{-\Xi(v,w)} \le e^{c_3} \io u = e^{c_3} \io u_0
	\qquad \mbox{for all } t\in (0,\tm),
  \eas
  and that thus
  \be{11.6}
	M_0(T) \le \max \bigg\{ 1 \, , \, e^{c_3} \io u_0 \bigg\}
	\qquad \mbox{for all } T\in (0,\tm).
  \ee
  Moreover, \eqref{11.3}-\eqref{11.5} enable us to estimate $M_k(T)$ for $k\ge 1$ and $T\in (0,\tm)$ by means of a refined testing
  procedure on the basis of \eqref{0a}: Indeed, for any such $k$ we see that abbreviating $p:=p_k$ we have
  \bea{11.7}
	\frac{d}{dt} \io e^{\Xi(v,w)} a^p
	&=& p \io e^{\Xi(v,w)} a^{p-1} \cdot \bigg\{ e^{-\Xi(v,w)} \na \cdot \Big\{ e^{\Xi(v,w)} \na a \Big\} \nn\\
	& & \hs{20mm}
	- e^{-\Xi(v,w)} \na \cdot \Big\{ a e^{\Xi(v,w)} \Big(\chi(v,w)-\Xi_v(v,w)\Big) \na v \Big\} \nn\\
	& & \hs{20mm}
	- a \Xi_v(v,w) v_t + a\xi(v,w) \psi(aw e^{\Xi(v,w)}) \bigg\} \nn\\
	& & + \io e^{\Xi(v,w)} a^p \cdot \Big\{ \Xi_v(v,w) v_t - \xi(v,w) \psi(awe^{\Xi(v,w)}) \Big\} \nn\\[1mm]
	&=& p \io a^{p-1} \na \cdot \Big\{ e^{\Xi(v,w)}\na a\Big\}
	- p\io a^{p-1} \na \cdot \Big\{ a e^{\Xi(v,w)} \Big(\chi(v,w)-\Xi_v(v,w)\Big) \na v \Big\} \nn\\
	& & - p\io e^{\Xi(v,w)} \Xi_v(v,w) a^p v_t
	+ p\io e^{\Xi(v,w)} \xi(v,w) a^p \psi(aw e^{\Xi(v,w)}) \nn\\
	& & + \io e^{\Xi(v,w)} \Xi_v(v,w) a^p v_t
	- \io e^{\Xi(v,w)} \xi(v,w) a^p \psi(awe^{\Xi(v,w)}) \nn\\[1mm]
	&=& -p(p-1) \io e^{\Xi(v,w)} a^{p-2} |\na a|^2 \nn\\
	& & + p(p-1) \io e^{\Xi(v,w)} \Big(\chi(v,w)-\Xi_v(v,w)\Big) a^{p-1} \na a \cdot\na v \nn\\
	& & - (p-1) \io e^{\Xi(v,w)} \Xi_v(v,w) a^p v_t \nn\\
	& & + (p-1) \io e^{\Xi(v,w)} \xi (v,w) a^p \psi (awe^{\Xi(v,w)})
	\qquad \mbox{for all } t\in (0,\tm).
  \eea
By Young's inequality, \eqref{11.3} and \eqref{11.4},
  \bea{11.8}
	& & \hs{-20mm}
	p(p-1) \io e^{\Xi(v,w)} \Big(\chi(v,w)-\Xi_v(v,w)\Big) a^{p-1} \na a \cdot \na v \nn\\
	&\le& \frac{p(p-1)}{2} \io e^{\Xi(v,w)} a^{p-2} |\na a|^2 \nn\\
	& & + \frac{p(p-1)}{2} \io e^{\Xi(v,w)} \Big(\chi(v,w)-\Xi_v(v,w)\Big)^2 a^p |\na v|^2 \nn\\
	&\le& \frac{p(p-1)}{2} \io e^{\Xi(v,w)} a^{p-2} |\na a|^2 \nn\\
	& & + \frac{p(p-1)}{2} e^{c_3} (c_5+c_4)^2 c_7^2 \io a^p
	\qquad \mbox{for all } t\in (0,\tm),
  \eea
  and due to \eqref{11.3} and \eqref{11.33}, 
  \bea{11.9}
	& & \hs{-20mm}
	(p-1) \io e^{\Xi(v,w)} \xi(v,w) a^p \psi( awe^{\Xi(v,w)}) \nn\\
	&\le& (p-1) e^{c_3} c_2 \io a^p \cdot \Big\{ c_6 \cdot (aw e^{\Xi(v,w)})^\frac{1}{3} + c_6 \Big\} \nn\\
	&\le& (p-1) e^\frac{4c_3}{3} c_2 c_1^\frac{1}{3} c_6 \io a^{p+\frac{1}{3}}
	+ (p-1) e^{c_3} c_2 c_6 \io a^p
	\qquad \mbox{for all } t\in (0,\tm).
  \eea
  Since furthermore an application of the H\"older inequality shows that thanks to \eqref{11.3} we have
  \bas
	-(p-1) \io e^{\Xi(v,w)} \Xi_v(v,w) a^p v_t	
	&\le& (p-1) \|v_t\|_{L^4(\Omega)} \cdot \bigg\{ \io e^{\frac{4}{3} \Xi(v,w)} |\Xi_v(v,w)|^\frac{4}{3} a^\frac{4p}{3} 
		\bigg\}^\frac{3}{4} \\
	&\le& (p-1) e^{c_3} c_4 \|v_t\|_{L^4(\Omega)} \cdot \bigg\{ \io a^\frac{4p}{3} \bigg\}^\frac{3}{4}
	\qquad \mbox{for all } t\in (0,\tm),
  \eas
  and since
  \bas
	\io a^p \le \io a^{p+\frac{1}{3}} + |\Omega|
	\qquad \mbox{for all } t\in (0,\tm)
  \eas
  by Young's inequality, from \eqref{11.7}-\eqref{11.9} we infer the existence of $c_9\in (0,2)$ and $c_{10}>0$ such that for
  any $k\ge 1$,
  \bea{11.10}
	& & \hs{-20mm}
	\frac{d}{dt} \io e^{\Xi(v,w)} a^{p_k}
	+ c_9 \io |\nabla a^\frac{p_k}{2}|^2 \nn\\
	&\le& c_{10} p_k^2 \io a^{p_k+\frac{1}{3}}
	+ c_{10} p_k \|v_t\|_{L^4(\Omega)} \cdot \bigg\{ \io a^\frac{4p_k}{3} \bigg\}^\frac{3}{4}
	+ c_{10} p_k^2
	\qquad \mbox{for all } t\in (0,\tm).
  \eea
  We now employ the Gagliardo-Nirenberg inequality to fix $c_{11} \ge 1$ such that
  \be{11.11}
	\|\varphi\|_{L^\frac{8}{3}(\Omega)}^2
	\le c_{11} \|\na \varphi\|_{L^2(\Omega)}^\frac{5}{4} \|\varphi\|_{L^1(\Omega)}^\frac{3}{4} 
	+ c_{11} \|\varphi\|_{L^1(\Omega)}^2
	\qquad \mbox{for all } \varphi\in W^{1,2}(\Omega),
  \ee
  and note that due to the H\"older inequality and Young's inequality, this implies that for each $p\ge 1$,
  \bas
	\|\varphi\|_{L^\frac{2(p+\frac{1}{3})}{p}(\Omega)}^\frac{2(p+\frac{1}{3})}{p}
	&\le& \|\varphi\|_{L^\frac{8}{3}(\Omega)}^\frac{8(p+\frac{2}{3})}{5p} \|\varphi\|_{L^1(\Omega)}^\frac{2(p-1)}{5p} \\
	&\le& c_{11}^\frac{4(p+\frac{2}{3})}{5p} \cdot 
	\bigg\{ \|\na \varphi\|_{L^2(\Omega)}^\frac{5}{4} \|\varphi\|_{L^1(\Omega)}^\frac{3}{4}
	+ \|\varphi\|_{L^1(\Omega)}^2 \bigg\}^\frac{4(p+\frac{2}{3})}{5p} 
	\|\varphi\|_{L^1(\Omega)}^\frac{2(p-1)}{5p} \\
	&\le& 2^\frac{4(p+\frac{2}{3})}{5p} c_{11}^\frac{4(p+\frac{2}{3})}{5p} \cdot
	\bigg\{ \|\na \varphi\|_{L^2(\Omega)}^\frac{p+\frac{2}{3}}{p} \|\varphi\|_{L^1(\Omega)}^\frac{3(p+\frac{2}{3})}{5p}
		\|\varphi\|_{L^1(\Omega)}^\frac{2(p-1)}{5p}
	+ \|\varphi\|_{L^1(\Omega)}^\frac{8(p+\frac{2}{3})}{5p} \|\varphi\|_{L^1(\Omega)}^\frac{2(p-1)}{5p} \bigg\} \\
	&=& 2^\frac{4(p+\frac{2}{3})}{5p} c_{11}^\frac{4(p+\frac{2}{3})}{5p} \cdot
	\bigg\{ \|\na \varphi\|_{L^2(\Omega)}^\frac{p+\frac{2}{3}}{p} \|\varphi\|_{L^1(\Omega)}
	+ \|\varphi\|_{L^1(\Omega)}^\frac{2(p+\frac{1}{3})}{p} \bigg\}
	\qquad \mbox{for all } \varphi\in W^{1,2}(\Omega),
  \eas
  because for any such $p$ we have $1 \le \frac{2(p+\frac{1}{3})}{p} \le \frac{8}{3}$.
  Since, apart from that, $\frac{4(p+\frac{2}{3})}{5p} \le \frac{4}{3}$ for all $p\ge 1$, this entails that writing
  $c_{12}:=(2c_{11})^\frac{4}{3}$, for any choice of $k\ge 1$ we can estimate
  \be{11.12}
	\|\varphi\|_{L^\frac{2(p_k+\frac{1}{3})}{p_k}(\Omega)}^\frac{2(p_k+\frac{1}{3})}{p_k}
	\le c_{12} \|\na \varphi\|_{L^2(\Omega)}^\frac{p_k+\frac{2}{3}}{p_k} \|\varphi\|_{L^1(\Omega)}
	+ c_{12} \|\varphi\|_{L^1(\Omega)}^\frac{2(p_k+\frac{1}{3})}{p_k}
	\qquad \mbox{for all } \varphi\in W^{1,2}(\Omega).
  \ee
  With these interpolation properties at hand, we can proceed to estimate the first two summands on the right-hand side of \eqref{11.10}
  by making use of the observation that whenever $k\ge 1$ and $T>0$,
  \be{11.13}
	\io a^\frac{p_k}{2} = \io a^{p_{k-1}} \le M_{k-1}(T)
	\qquad \mbox{for all } t\in (0,T).
  \ee
  According to \eqref{11.12}, namely, this firstly implies that due to an application of 
  the Gagliardo-Nirenberg and Young inequalities,
  relying on the fact that $\frac{p_k+\frac{2}{3}}{p_k} \le \frac{5}{3} \le 2$ for all $k\ge 1$,
  \bas
	c_{10} p_k^2 \io a^{p_k+\frac{1}{3}}
	&=& c_{10} p_k^2 \|a^\frac{p_k}{2}\|_{L^\frac{2(p_k+\frac{1}{3})}{p_k}(\Omega)}^\frac{2(p_k+\frac{1}{3})}{p_k} \\
	&\le& c_{10} c_{12} p_k^2 \|\na a^\frac{p_k}{2}\|_{L^2(\Omega)}^\frac{p_k+\frac{2}{3}}{p_k} \|a^\frac{p_k}{2}\|_{L^1(\Omega)}
	+ c_{10} c_{12} p_k^2 \|a^\frac{p_k}{2}\|_{L^1(\Omega)}^\frac{2(p_k+\frac{1}{3})}{p_k} \\
	&\le& c_{10} c_{12} p_k^2 M_{k-1}(T) \|\na a^\frac{p_k}{2}\|_{L^2(\Omega)}^\frac{p_k+\frac{2}{3}}{p_k}
	+ c_{10} c_{12} p_k^2 M_{k-1}^\frac{2(p_k+\frac{1}{3})}{p_k} (T) \\
	&=& \bigg\{ \frac{c_9}{2} \io |\na a^\frac{p_k}{2}|^2 \bigg\}^\frac{p_k+\frac{2}{3}}{2p_k}
	\cdot \Big(\frac{2}{c_9}\Big)^\frac{p_k+\frac{2}{3}}{2p_k} c_{10} c_{12} p_k^2 M_{k-1}(T)
	+ c_{10} c_{12} p_k^2 M_{k-1}^\frac{2(p_k+\frac{1}{3})}{p_k} (T) \\
	&\le& \frac{c_9}{2} \io |\na a^\frac{p_k}{2}|^2 
	+\Big(\frac{2}{c_9}\Big)^\frac{p_k+\frac{2}{3}}{p_k-\frac{2}{3}} (c_{10} c_{12})^\frac{2p_k}{p_k-\frac{2}{3}} 
	p_k^\frac{4p_k}{p_k-\frac{2}{3}}  M_{k-1}^\frac{2p_k}{p_k-\frac{2}{3}}(T) \\
	& & + c_{10} c_{12} p_k^2 M_{k-1}^\frac{2(p_k+\frac{1}{3})}{p_k} (T)
	\qquad \mbox{for all $t\in (0,T)$, any $T\in (0,\tm)$ and each } k\ge 1.
  \eas
  Since $\frac{p_k+\frac{2}{3}}{p_k-\frac{2}{3}} \le 5$ and $1\le \frac{p_k}{p_k-\frac{2}{3}} \le 3$ as well as
  $\frac{2(p_k+\frac{1}{3})}{p_k} \le \frac{2p_k}{p_k-\frac{2}{3}}$ for all $k\ge 1$, and since $M_{k-1}(T)\ge 1$ for all
  $k\ge 1$ and $T\in (0,\tm)$, using that $\frac{2}{c_9} \ge 1$ we thus obtain that if we let 
  \bas
	c_{13}:=\Big(\frac{2}{c_9}\Big)^5 \cdot (c_{10} c_{12} )^3 + c_{10} c_{12},
  \eas
  then
  \be{11.14}
	c_{10} p_k^2 \io a^{p_k+\frac{1}{3}}
	\le \frac{c_9}{2} \io |\na a^\frac{p_k}{2}|^2 
	+ c_{13} p_k^{12} M_{k-1}^\frac{2p_k}{p_k-\frac{2}{3}}(T)
	\qquad \mbox{for all $t\in (0,T)$, $T\in (0,\tm)$ and } k\ge 1.
  \ee
  We next use \eqref{11.13} together with \eqref{11.11} to similarly see that again thanks to Young's inequality,
  whenever $T\in (0,\tm)$ and $k\ge 1$,
  \bas
	& & \hs{-20mm}
	c_{10} p_k \|v_t\|_{L^4(\Omega)} \cdot \bigg\{ \io a^\frac{4p_k}{3} \bigg\}^\frac{3}{4} \nn\\
	&=& c_{10} p_k \|v_t\|_{L^4(\Omega)} \|a^\frac{p_k}{2}\|_{L^\frac{8}{3}(\Omega)}^2 \\
	&\le& c_{10} c_{11} p_k M_{k-1}^\frac{3}{4}(T) \|v_t\|_{L^4(\Omega)} \|\na a^\frac{p_k}{2}\|_{L^2(\Omega)}^\frac{5}{4}
	+ c_{10} c_{11} p_k M_{k-1}^2(T) \|v_t\|_{L^4(\Omega)} \\
	&=& \bigg\{ \frac{c_9}{2} \io |\na a^\frac{p_k}{2} |^2 \bigg\}^\frac{5}{8}
	\cdot \Big(\frac{2}{c_9}\Big)^\frac{5}{8} c_{10} c_{11} p_k M_{k-1}^\frac{3}{4}(T) \|v_t\|_{L^4(\Omega)} 
	+ c_{10} c_{11} p_k M_{k-1}^2(T) \|v_t\|_{L^4(\Omega)} \\
	&\le& \frac{c_9}{2} \io |\na a^\frac{p_k}{2} |^2
	+ \Big(\frac{2}{c_9}\Big)^\frac{5}{3} (c_{10} c_{11})^\frac{8}{3} p_k^\frac{8}{3} M_{k-1}^2(T) 
		\|v_t\|_{L^4(\Omega)}^\frac{8}{3} \\
	& & + c_{10} c_{11} p_k M_{k-1}^2(T) \|v_t\|_{L^4(\Omega)} 
	\qquad \mbox{for all } t\in (0,T),
  \eas
  so that since $p_k\ge 1$ for all $k\ge 1$ and $\max \{s^\frac{8}{3},s\} \le s^4 +1$ for all $s\ge 0$ by Young's inequality, 
  it follows that
  \be{11.15}
	c_{10} p_k \|v_t\|_{L^4(\Omega)} \cdot \bigg\{ \io a^\frac{4p_k}{3} \bigg\}^\frac{3}{4} 
	\le \frac{c_9}{2} \io |\na a^\frac{p_k}{2} |^2
	+ c_{14} p_k^\frac{8}{3} M_{k-1}^2(T) \cdot \bigg\{ \io v_t^4 + 1 \bigg\}
  \ee
  for any such $t, T$ and $k$, with $c_{14}:=(\frac{2}{c_9})^\frac{5}{3} (c_{10} c_{11})^\frac{8}{3} + c_{10} c_{11}$.\\
  We now combine this with \eqref{11.14} and \eqref{11.10} to infer that as
  \bas
	p_k^\frac{8}{3} M_{k-1}^2(T) \le p_k^{12} M_{k-1}^\frac{2p_k}{p_k-\frac{2}{3}}(T)
	\quad \mbox{and} \quad
	p_k^2 \le p_k^{12}
	\qquad \mbox{for all $T\in (0,\tm)$ and } k\ge 1,
  \eas
  writing $c_{15}:=c_{10}+c_{13}+c_{14}$ and $\theta_k:=\frac{2p_k}{p_k-\frac{2}{3}}$ for $k\ge 1$, we have
  \bas
	\frac{d}{dt} \io e^{\Xi(v,w)} a^{p_k}
	\le c_{15} p_k^{12} M_{k-1}^{\theta_k}(T) \cdot \bigg\{ \io v_t^4 +1\bigg\}
	\qquad \mbox{for all $t\in (0,T)$, $T\in (0,\tm)$ and } k\ge 1.
  \eas
Upon an integration this entails that according to \eqref{11.5},
  \bas
	\io e^{\Xi(v,w)} a^{p_k}
	\le \io e^{\Xi(v_0,w_0)} a^{p_k}(\cdot,0)
	+ c_{16} p_k^{12} M_{k-1}^{\theta_k}(T)
	\qquad \mbox{for all $t\in (0,T)$, $T\in (0,\tm)$ and } k\ge 1,
  \eas
  with $c_{16}:= c_{15} \cdot (c_8+\tm)$. 
  Since by \eqref{11.3}
  \bas
	e^{-c_3} \io a^{p_k} \le \io e^{\Xi(v,w)} a^{p_k}
	\le e^{c_3} \|a\|_{L^\infty(\Omega)}^{p_k}
	\qquad \mbox{for all $t\in [0,\tm)$ and } k\ge 1,
  \eas
  once more in view of \eqref{11.2} this reveals that
  \bea{11.16}
	M_k(T)
	&\le& 1 + e^{c_3} \cdot \bigg\{ \io e^{\Xi(v_0,w_0)} a^{p_k}(\cdot,0)
	+ c_{16} p_k^{12} M_{k-1}^{\theta_k}(T) \bigg\} \nn\\
	&\le& 1 + e^{2c_3} \|a(\cdot,0)\|_{L^\infty(\Omega)}^{p_k}
	+ e^{c_3} c_{16} p_k^{12} M_{k-1}^{\theta_k}(T) \nn\\
	&\le& e^{2c_3} \|a(\cdot,0)\|_{L^\infty(\Omega)}^{p_k} 
	+ c_{17} p_k^{12} M_{k-1}^{\theta_k}(T)
	\qquad \mbox{for all $T\in (0,\tm)$ and } k\ge 1,
  \eea
  where $c_{17}:=1+e^{c_3} c_{16}$.
  In light of \eqref{11.6}, by induction this firstly implies that
  \bas
	\ov{M}_k:=\sup_{T\in (0,\tm)} M_k(T)
  \eas
  is finite for all $k\ge 1$, and that, secondly, these numbers satisfy
  \be{11.17}
	\ov{M}_k \le e^{2c_3} \|a(\cdot,0)\|_{L^\infty(\Omega)}^{p_k}
	+ c_{17} p_k^{12} \ov{M}_{k-1}^{\theta_k}
	\qquad \mbox{for all } k\ge 1.
  \ee
  Now if $\ov{M}_k \le 2 e^{2c_3} \|a(\cdot,0)\|_{L^\infty(\Omega)}^{p_k}$ for infinitely many $k\in\N$, then, since
  $p_k\to\infty$ as $k\to\infty$, it holds that
  \bas
	\|a(\cdot,t)\|_{L^\infty(\Omega)}
	= \liminf_{k\to\infty} \|a(\cdot,t)\|_{L^{p_k}(\Omega)}
	\le \liminf_{k\to\infty} \ov{M}_k^\frac{1}{p_k}
	\le \|a(\cdot,0)\|_{L^\infty(\Omega)}
	\qquad \mbox{for all } T\in (0,\tm).
  \eas
  In the opposite case, there exists $k_0\in\N$ such that $e^{2c_3} \|a(\cdot,0)\|_{L^\infty(\Omega)}^{p_k} \le \frac{1}{2} \ov{M}_k$
  and hence, by \eqref{11.17},
  \bas
	\ov{M}_k \le 2c_{17} p_k^{12} \ov{M}_{k-1}^{\theta_k}
  \eas
  for all $k\ge k_0$.
  Recalling that $p_k=2^k$ for $k\ge 0$, we can then fix $b>1$ large enough such that
  \bas
	\ov{M}_k \le b^k \ov{M}_{k-1}^{\theta_k}
	\qquad \mbox{for all } k\ge 1,
  \eas
  where by definition of $(\theta_k)_{k\ge 1}$,
  \bas
	\frac{\theta_k}{2}-1
	= \frac{p_k}{p_k-\frac{2}{3}} -1
	= \frac{2}{3p_k-2}
	= \frac{2}{3\cdot 2^k -2}
	\le \frac{2}{3\cdot 2^k - 2^k}
	= \frac{1}{2^k}
	\qquad \mbox{for all } k\ge 1.
  \eas
  Therefore, according to an elementary estimate derived in \cite[Lemma 4.3]{win_NON_exp},
  \bas
	\ov{M}_k \le b^{k+ e\cdot 2^{k+1}} \ov{M}_0^{e\cdot 2^k}
	\qquad \mbox{for all } k\ge 1,
  \eas
  so that
  \bas	
	\limsup_{k\to\infty} \ov{M}_k^\frac{1}{p_k}
	\le \limsup_{k\to\infty} \bigg\{ b^{k\cdot 2^{-k} +2e} \ov{M}_0^e \bigg\},
  \eas
  and hence we may conclude that \eqref{11.1} also holds in this case.
\qed
\subsection{Bounds in $C^{2+\theta}(\bom)$. Proof of \cref{theo27}}
Having at hand the above information, and especially $L^\infty$ bounds for all solution components,
we can build our subsequent regularity arguments on the following elementary observation.
\begin{Lemma} \label{lem111}
  We have
  \bea{111.1}
	a_t 
	&=& \Delta a + A_1(x,t) \na v\cdot\na a + A_2(x,t) \na w\cdot\na a
	+ A_3(x,t) |\na v|^2
	+ A_4(x,t) \na v\cdot\na w \nn\\
	& & + A_5(x,t) \Delta v + A_6(x,t) v_t 
	+ A_7(x,t)
	\qquad \mbox{for all $x\in\Omega$ and } t\in (0,\tm)
  \eea
  and
  \be{211.1}
	\na w_t = B_1(x,t) \na a + B_2(x,t) \na w + B_3(x,t)\na v
	\qquad \mbox{for all $x\in\Omega$ and } t\in (0,\tm),
  \ee
  where
  \begin{align}\label{A}
	A_1(x,t)&:= 2\Xi_v(v,w) - \chi(v,w), \nn\\
	A_2(x,t)&:=\xi(v,w), \nn\\
	A_3(x,t)&:= a\big( -\chi(v,w) \Xi_v(v,w) + \Xi_v^2(v,w) -\chi_v(v,w) + \Xi_{vv}(v,w) \big),
	\nn\\
	A_4(x,t)&:= a \big( -\chi(v,w) \xi(v,w) +\Xi_v(v,w) \xi(v,w) - \chi_w(v,w) + \xi_v(v,w) \big), 
	\\
	A_5(x,t)&:= a \big( -\chi(v,w) + \Xi_v(v,w) \big), 
	\nn\\
	A_6(x,t)&:= a \big( -\Xi_v(v,w) \big)
	\qquad \mbox{and} \nn\\
	A_7(x,t)&:= a\big( \xi(v,w) \psi(awe^{\Xi(v,w)}) \big),
	\nn
  \end{align}

  as well as
  \be{B}
	\begin{array}{l}
	B_1(x,t):=-\psi'(awe^{\Xi(v,w)}) \cdot we^{\Xi(v,w)}, \\
	B_2(x,t):=-\psi'(awe^{\Xi(v,w)}) \cdot ae^{\Xi(v,w)} \cdot (1+w\xi(v,w))
	\qquad \mbox{and} \\
	B_3(x,t):=-\psi'(awe^{\Xi(v,w)}) \cdot awe^{\Xi(v,w)} \cdot \Xi_v(v,w)
	\end{array}
  \ee
  for $(x,t)\in\Omega\times (0,\tm)$.
\end{Lemma} 
\proof
  This can be seen by straightforward differentiation in \eqref{0a}.\hfill 
\qed\\[-2ex]

\noindent
In fact, due to the boundedness features gathered in the previous sections, the latter lemma implies the following.
\begin{Corollary}\label{cor112}
  Assume that $\tm<\infty$. Then there exists $C>0$ such that
  \be{112.1}
	|a_t-\Delta a|
	\le C \cdot \Big\{ |\na w\cdot\na a| + |\na a| + |\na w| + |\Delta v| + |v_t| + 1 \Big\}
	\qquad \mbox{in } \Omega\times (0,\tm)
  \ee
  and that
  \be{112.2}
	|\na w_t| \le C \cdot \Big\{ |\na a| + |\na w| + 1 \Big\}
	\qquad \mbox{in } \Omega\times (0,\tm).
  \ee
\end{Corollary}
\proof
  Since \eqref{chi}, \eqref{xi}, \eqref{psi} and \eqref{psi1} together with \eqref{winfty}, \eqref{vinfty} and \eqref{lem11}
  warrant boundedness of all the functions $A_i$ and $B_j$, $i\in\{1,...,7\}$, $j\in \{1,2,3\}$ in \eqref{A} and \eqref{B}, 
  and since moreover also $\na v$ is bounded
  in $\Omega\times (0,\tm)$ by \cref{lem9}, both inequalities are direct consequences of \cref{lem111}.
\qed\\[-2ex]

\noindent
Once again relying on the planarity of the spatial setting, through an analysis of $\io |\na a|^2 + c\io |\na w|^4$ 
along trajectories, with suitably chosen $c>0$, we shall next obtain, inter alia, some information on $a$ at the level
of second-order derivatives.
\begin{Lemma} \label{lem12}
  If $\tm<\infty$, then there exists $C>0$ such that
  \be{12.1}
	\io |\na w(\cdot,t)|^4 \le C
	\qquad \mbox{for all } t\in (0,\tm),
  \ee
  and moreover we have
  \be{12.11}
	\int_0^{\tm} \io |\Delta a|^2 <\infty.
  \ee
\end{Lemma} 
\proof
  We employ Young's inequality and \cref{cor112} to see that with some $c_1>0$,
  \bea{12.2}
	\frac{1}{2} \frac{d}{dt} \io |\na a|^2
	+ \io |\Delta a|^2
	&=& - \io \Delta a \cdot (a_t-\Delta a) \nn\\
	&\le& \frac{1}{4} \io |\Delta a|^2 + \io |a_t-\Delta a|^2 \nn\\
	&\le& \frac{1}{4} \io |\Delta a|^2
	+ c_1 \io |\na w|^2 |\na a|^2
	+ c_1 \io |\na a|^2
	+ c_1 \io |\na w|^2 \nn\\
	& & + c_1 \io |\Delta v|^2
	+ c_1 \io v_t^2
	+ c_1
	\qquad \mbox{for all } t\in (0,\tm).
  \eea
  Writing $c_2:=\|a\|_{L^\infty(\Omega\times (0,\tm))}<\infty$ and taking $c_3>0$ such that in accordance with the
  Gagliardo-Nirenberg inequality and elliptic regularity theory (\cite{GT}) we have
  \be{12.3}
	\io |\na \varphi|^4 \le c_3\|\Delta\varphi\|_{L^2(\Omega)}^2 \|\varphi\|_{L^\infty(\Omega)}^2
	\qquad \mbox{for all $\varphi\in W^{2,2}(\Omega)$ such that $\frac{\partial\varphi}{\partial\nu}=0$ on $\pO$,}
  \ee
  again by Young's inequality we estimate
  \bas
	c_1 \io |\na w|^2 |\na a|^2 + c_1 \io |\na a|^2
	&\le& \frac{1}{4c_2^2 c_3} \io |\na a|^4 
	+ 2c_1^2 c_2^2 c_3 \io |\na w|^4 
	+ 2c_1^2 c_2^2 c_3 |\Omega| \\
	&\le& \frac{1}{4c_2^2} \|\Delta a\|_{L^2(\Omega)}^2 \|a\|_{L^\infty(\Omega)}^2
	+ 2c_1^2 c_2^2 c_3 \io |\na w|^4 
	+ 2c_1^2 c_2^2 c_3 |\Omega| \\
	&\le& \frac{1}{4} \io |\Delta a|^2
	+ 2c_1^2 c_2^2 c_3 \io |\na w|^4 
	+ 2c_1^2 c_2^2 c_3 |\Omega|
	\qquad \mbox{for all } t\in (0,\tm).
  \eas
  Since Young's inequality moreover entails that
  \bas
	c_1 \io |\na w|^2
	\le \frac{1}{2} \io |\na w|^4
	+ \frac{c_1^2 |\Omega|}{2}
	\qquad \mbox{for all } t\in (0,\tm),
  \eas
  from \eqref{12.2} we thus obtain the inequality
  \be{12.4}
	\frac{d}{dt} \io |\na a|^2
	+ \io |\Delta a|^2
	\le c_3 \io |\na w|^4
	+ 2c_1 \io |\Delta v|^2
	+ 2c_1 \io v_t^2 + c_4
	\qquad \mbox{for all } t\in (0,\tm)
  \ee
  with $c_3:=4c_1^2 c_2^2 c_3 +1$ and $c_4:=4c_1^2 c_2^2 c_3 |\Omega| + c_1^2 |\Omega| +2c_1$.
  In order to appropriately control the growth induced by the first summand on the right-hand side herein, we combine
  \cref{cor112} with Young's inequality to see that with some $c_5>0$ we have
  \bea{12.5}
	\frac{d}{dt} \io |\na w|^4
	&=& 4 \io |\na w|^2 \na w \cdot\na w_t \nn\\
	&\le& c_5 \io |\na w|^3 \cdot \Big\{ |\na a| + |\na w| + 1 \Big\} \nn\\
	&\le& 3c_5 \io |\na w|^4 + c_5 \io |\na a|^4 + c_5|\Omega|
	\qquad \mbox{for all } t\in (0,\tm),
  \eea
  for which once more by \eqref{12.3} we estimate
  \bas
	c_5 \io |\na a|^4 \le c_6 \io |\Delta a|^2
	\qquad \mbox{for all } t\in (0,\tm)
  \eas
  with $c_6:=c_2^2 c_3 c_5$.
  Therefore, \eqref{12.4} in conjunction with \eqref{12.5} show that
  \bas
	& & \hs{-20mm}
	\frac{d}{dt} \bigg\{ \io |\na a|^2 + \frac{1}{2c_6} \io |\na w|^4 \bigg\}
	+ \frac{1}{2} \io |\Delta a|^2 \\
	&\le& \Big(c_3 + \frac{3c_5}{2c_6}\Big) \io |\na w|^4
	+ 2c_1 \io |\Delta v|^2 + 2c_1 \io v_t^2 + c_4 + \frac{c_5|\Omega|}{2c_6}
	\qquad \mbox{for all } t\in (0,\tm),
  \eas
  and that hence the functions given by $y(t):=\io |\na a(\cdot,t)|^2 + \frac{1}{2c_6} \io |\na w(\cdot,t)|^4$,
  $g(t):=\frac{1}{2} \io |\Delta a(\cdot,t)|^2$ and 
  $h(t) :=2c_1 \io |\Delta v(\cdot,t)|^2 + 2c_1 \io v_t^2(\cdot,t) + c_4 + \frac{c_5|\Omega|}{2c_6}$, $t\in [0,\tm)$, satisfy
  $y'(t) + g(t)  \le 2c_6 (c_3+\frac{3c_5}{2c_6}) y(t) + h(t)$ for all $t\in (0,\tm)$.
  Since $\int_0^{\tm} h(t) dt<\infty$ by \cref{lem8}, this can readily be seen to firstly imply that 
  $\sup_{t\in (0,\tm)} y(t)<\infty$, and to secondly entail that also $\int_0^{\tm} g(t) dt$ is finite, whereby
  both \eqref{12.1} and \eqref{12.11} are established.
\qed\\[-2ex]

\noindent
In order to prepare an extension of the latter to a result in the flavor of that from \cref{lem9},
let us draw the following further conclusion of \cref{cor112} when combined with maximal Sobolev regularity
theory.
\begin{Lemma} \label{lem13}
  Suppose that $\tm<\infty$. Then for all $p>1$ one can find $C(p)>0$ such that
  \be{13.1}
	\int_0^{\tm} \Big\{ \|a(\cdot,t)\|_{W^{2,p}(\Omega)}^p + \|a_t(\cdot,t)\|_{L^p(\Omega)}^p \Big\} dt
	\le C(p) \cdot \bigg\{ 1 + \int_0^{\tm} \io |\na w|^{2p} \bigg\}.
  \ee
\end{Lemma} 
\proof
  We once again recall a standard result from maximal Sobolev regularity theory (\cite{giga_sohr}) to fix $c_1(p)>0$ such that
  if $\varphi\in C^{2,1}(\bom\times [0,\tm))$ and $h\in C^0(\bom\times [0,\tm))$ are such that
  \be{13.2}
	\left\{ \begin{array}{ll}
	\varphi_t = \Delta\varphi + h(x,t),
	\qquad & x\in \Omega, \ t\in (0,\tm), \\[1mm]
	\frac{\partial\varphi}{\partial\nu}=0,
	& x\in\pO, \ t\in (0,\tm), \\[1mm]
	\varphi(x,0)=a(x,0),
	\qquad & x\in\Omega,
	\end{array} \right.
  \ee
  then
  \be{13.3}
	\int_0^T \Big\{ \|\varphi(\cdot,t)\|_{W^{2,p}(\Omega)}^p + \|\varphi_t(\cdot,t)\|_{L^p(\Omega)}^p \Big\} dt
	\le c_1(p) \cdot \bigg\{ 1 + \int_0^T \io |h|^p \bigg\}
	\qquad \mbox{for all } T\in (0,\tm).
  \ee
  To adequately make use of this in the present situation, we once more invoke \cref{cor112} to choose $c_2>0$ in such a way
  that
  \bas
	|a_t-\Delta a|
	\le c_2 \cdot \Big\{ |\na w\cdot\na a| + |\na a| + |\na w| + |\Delta v| + |v_t| +1 \Big\}	
	\qquad \mbox{in } \Omega\times (0,\tm),
  \eas
  and employ the Gagliardo-Nirenberg inequality along with \cref{lem11} to see that with some $c_3(p)>0$ and $c_4(p)>0$ we have
  \be{13.4}
	\|\na a\|_{L^{2p}(\Omega)}^p
	\le c_3(p) \|a\|_{W^{2,p}(\Omega)}^\frac{p}{2} \|a\|_{L^\infty(\Omega)}^\frac{p}{2}
	\le c_4(p) \|a\|_{W^{2,p}(\Omega)}^\frac{p}{2}
	\qquad \mbox{for all } t\in (0,\tm).
  \ee
  Therefore, an application of \eqref{13.3} to $\varphi:=a$ and $h:=a_t-\Delta a$ shows that due to the Cauchy-Schwarz inequality,
  \bea{13.5}
	& & \hs{-20mm}
	\int_0^T \Big\{ \|a(\cdot,t)\|_{W^{2,p}(\Omega)}^p + \|a_t(\cdot,t)\|_{L^p(\Omega)}^p \Big\} dt \nn\\
	&\le& c_1(p) + c_1(p) c_2^p \int_0^T \io \Big\{|\na w\cdot\na a| + |\na a| + |\na w| + |\Delta v| + |v_t| + 1 \Big\}^p \nn\\
	&\le& c_1(p) + 6^p c_1(p) c_2^p 
	\int_0^T \Big\{ \|\na w(\cdot,t)\cdot\na a(\cdot,t)\|_{L^p(\Omega)}^p
	+ \|\na a(\cdot,t)\|_{L^p(\Omega)}^p
	+ \|\na w(\cdot,t)\|_{L^p(\Omega)}^p \nn\\
	& & \hs{40mm}
	+ \|\Delta v(\cdot,t)\|_{L^p(\Omega)}^p
	+ \|v_t(\cdot,t)\|_{L^p(\Omega)}^p
	+ 1 \Big\} dt \nn\\
	&\le& c_1(p) + 6^p c_1(p) c_2^p 
	\int_0^T \Big\{ \|\na w(\cdot,t)\|_{L^{2p}(\Omega)}^p \|a(\cdot,t)\|_{L^{2p}(\Omega)}^p \nn\\
	& & \hs{40mm}
	+ |\Omega|^\frac{1}{2} \|\na a(\cdot,t)\|_{L^{2p}(\Omega)}^p
	+ |\Omega|^\frac{1}{2} \|\na w(\cdot,t)\|_{L^{2p}(\Omega)}^p \nn\\
	& & \hs{40mm}
	+ \|\Delta v(\cdot,t)\|_{L^p(\Omega)}^p
	+ \|v_t(\cdot,t)\|_{L^p(\Omega)}^p
	+ 1 \Big\} dt 
  \eea
  for all $T\in (0,\tm)$.
  Here, \eqref{13.5} together with Young's inequality ensures that if we let 
  $c_5(p):=\frac{1}{2} \cdot 6^{2p} c_1^2(p) c_2^{2p} c_4^2(p)$, then
  \bas
	& & \hs{-20mm}
	6^p c_1(p) c_2^p
	\int_0^T \Big\{ \|\na w(\cdot,t)\|_{L^{2p}(\Omega)}^p \|\na a(\cdot,t)\|_{L^{2p}(\Omega)}^p
	+ |\Omega|^\frac{1}{2} \|\na a(\cdot,t)\|_{L^{2p}(\Omega)}^p \Big\} dt \\
	&\le& 6^p c_1(p) c_2^p c_4(p) \int_0^T \Big\{ \|\na w(\cdot,t)\|_{L^{2p}(\Omega)}^p + |\Omega|^\frac{1}{2} \Big\}
		\cdot \|a(\cdot,t)\|_{W^{2,p}(\Omega)}^\frac{p}{2} dt \\
	&\le& \frac{1}{2} \int_0^T \|a(\cdot,t)\|_{W^{2,p}(\Omega)}^p dt
	+ c_5(p) \int_0^T \Big\{ \|\na w(\cdot,t)\|_{L^{2p}(\Omega)}^p + |\Omega|^\frac{1}{2} \Big\}^2 dt \\
	&\le& \frac{1}{2} \int_0^T \|a(\cdot,t)\|_{W^{2,p}(\Omega)}^p dt
	+ 2 c_5(p) \int_0^T \Big\{ \|\na w(\cdot,t)\|_{L^{2p}(\Omega)}^{2p} + |\Omega| \Big\} dt
  \eas
  for all $T\in (0,\tm)$.
  Since Young's inequality furthermore guarantees that
  \bas
	\int_0^T \|\na w(\cdot,t)\|_{L^{2p}(\Omega)}^p dt
	\le \int_0^T \io |\na w|^{2p} + \tm
	\qquad \mbox{for all } T\in (0,\tm),
  \eas
  from \eqref{13.5} we accordingly infer that
  \bas
	& & \hs{-20mm}
	\int_0^T \Big\{ \|a(\cdot,t)\|_{W^{2,p}(\Omega)}^p + \|a_t(\cdot,t)\|_{L^p(\Omega)}^p \Big\} dt \nn\\
	&\le& \Big\{ 2c_5(p) + 6^p c_1(p) c_2^p \Big\} \int_0^T \io |\na w|^{2p} \nn\\
	& & + 6^p c_1(p) c_2^p \int_0^T \io \Big\{ |\Delta v|^p + |v_t|^p \Big\} \nn\\
	& & + c_1(p) + 2c_5(p) |\Omega|\tm 
	+ 6^p c_1(p) c_2^p \tm
	+ 6^p c_1(p) c_2^p |\Omega| \tm
	\qquad \mbox{for all } T\in (0,\tm),
  \eas
  and that thus \eqref{13.1} holds due to the fact that $\int_0^T \io \big\{ |\Delta v|^p + |v_t|^p \big\}<\infty$ thanks to
  \cref{lem8}.
\qed\\[-2ex]

\noindent
Here a criterion for boundedness of the expression appearing on the right of \eqref{13.1} can be obtained on the basis
of \eqref{211.1}:
\begin{Lemma} \label{lem22}
  Assume that $\tm<\infty$, and that $q \ge 1$ is such that			
  \be{22.1}
	\int_0^{\tm} \|\na a(\cdot,t)\|_{L^q(\Omega)} dt < \infty.
  \ee
  Then there exists $C>0$ such that
  \be{22.2}
	\|\na w(\cdot,t)\|_{L^q(\Omega)} \le C
	\qquad \mbox{for all } t\in (0,\tm).
  \ee
\end{Lemma} 
\proof
  We observe that due to \eqref{winfty} and \cref{lem11} when combined with our overall regularity 
  and boundedness
  assumptions
  on $\psi$ and $\xi$, the functions $B_1$, $B_2$ and $B_3$ in \eqref{B} are all bounded in $\Omega\times (0,\tm)$.
  Since furthermore also $\na v$ is bounded in $\Omega\times (0,\tm)$ by \cref{lem9}, from \eqref{211.1} and our hypothesis
  \eqref{22.1} we thus infer that with some $c_1>0$ and $c_2>0$ we have
  \bas
	\|\na w(\cdot,t)\|_{L^q(\Omega)}
	&=& \bigg\| \na w_0 + \int_0^t \na w_t (\cdot,s) ds \bigg\|_{L^q(\Omega)} \\
	&\le& \|\na w_0\|_{L^q(\Omega)}
	+ \int_0^t \|B_1(\cdot,s) \na a(\cdot,s) \|_{L^q(\Omega)} ds \\
	& & + \int_0^t \|B_2(\cdot,s) \na w(\cdot,s) \|_{L^q(\Omega)} ds
	+ \int_0^t \|B_3(\cdot,s) \na v(\cdot,s) \|_{L^q(\Omega)} ds \\
	&\le& \|\na w_0\|_{L^q(\Omega)}
	+ c_1 \int_0^t \|\na a(\cdot,s)\|_{L^q(\Omega)} ds  \\
	& & + c_1 \int_0^t \|\na w(\cdot,s)\|_{L^q(\Omega)} ds 
	+ c_1 \int_0^t \|\na v(\cdot,s)\|_{L^q(\Omega)} ds \\
	&\le& c_1 \int_0^t \|\na w(\cdot,s)\|_{L^q(\Omega)} ds + c_2
	\qquad \mbox{for all } t\in (0,\tm).
  \eas
  An application of Gronwall's lemma hence asserts \eqref{22.2}.
\qed\\[-2ex]

\noindent
Thanks to the fact that the expression on the left of \eqref{12.11} essentially
dominates the quantity in \eqref{22.1} for arbitrary finite $q>1$ due to two-dimensional embeddings, 
from \cref{lem13} we can immediately draw the following conclusion.
\begin{Lemma} \label{lem23}
  Assume that $\tm<\infty$. Then there exists $\theta\in (0,1)$ such that
  \be{23.1}
	a \in C^{1+\theta,\theta}(\bom\times [0,\tm]).
  \ee
\end{Lemma} 
\proof
  Since \cref{lem12} implies that $\int_0^{\tm} \|a(\cdot,t)\|_{W^{2,2}(\Omega)}^2 dt < \infty$ and that thus clearly also
  the integral $\int_0^{\tm} \|\na a(\cdot,t)\|_{L^q(\Omega)} dt$ is finite for each $q\in [1,\infty)$ due to the two-dimensional
  Sobolev embedding theorem, from \cref{lem22} it follows that $\sup_{t\in (0,\tm)} \|\na w(\cdot,t)\|_{L^q(\Omega)}<\infty$
  for any such $q$.
  An application of \cref{lem13} therefore shows that
  \bas
	\int_0^{\tm} \Big\{ \|a(\cdot,t)\|_{W^{2,p}(\Omega)}^p + \|a_t(\cdot,t)\|_{L^p(\Omega)}^p \Big\} dt<\infty
	\qquad \mbox{for all } p\in (1,\infty),
  \eas
  so that \eqref{23.1} becomes a consequence of the embedding result from \cite{amann}.
\qed\\[-2ex]

\noindent
We next go back to \eqref{211.1} to see that our present information on H\"older regularity of $a,\na a, v$ and $\na v$
implies that if $\tm$ was finite, then actually also $w$ and $\na w$ must be H\"older continuous in $\bom\times [0,\tm]$.
\begin{Lemma} \label{lem24}
  If $\tm<\infty$, then there exists $\theta\in (0,1)$ such that
  \be{24.1}
	w\in C^{\theta}(\bom\times [0,\tm])
	\qquad \mbox{and} \qquad
	\na w \in C^\theta(\bom\times [0,\tm];\R^2).
  \ee
\end{Lemma} 
\proof
  We first recall that $\sup_{t\in (0,\tm)} \|\na w(\cdot,t)\|_{L^4(\Omega)}$ is finite by \cref{lem12},
  which together with (\cref{winfty}) and a Sobolev embedding theorem ensures that
  \be{24.11}
	\sup_{t\in (0,\tm)} \|w(\cdot,t)\|_{C^\frac{1}{2}(\bom)} <\infty.
  \ee
  Since also $\sup_{t\in (0,\tm)} \|a(\cdot,t)\|_{C^1(\bom)}$ is finite by \cref{lem23}, relying on $\psi'\in C^1([0,\infty))$ and the fact that $\xi,\Xi$ and $\Xi_v$ belong to $C^1([0,\infty)^2)$, from this we readily infer that
  there exists $\theta_1\in (0,1)$ such that in \eqref{B} we have
  \bas
	\sup_{t\in (0,\tm)} \|B_i(\cdot,t)\|_{C^{\theta_1}(\bom)} <\infty
	\qquad \mbox{for } i\in\{1,2,3\}.
  \eas
  Apart from that, \cref{lem23} and \cref{lem8} provide $\theta_2\in (0,\theta_1)$ fulfilling
  \bas
	\sup_{t\in (0,\tm)} \|\na a(\cdot,t)\|_{C^{\theta_2}(\bom)} <\infty
	\qquad \mbox{and} \qquad
	\sup_{t\in (0,\tm)} \|\na v(\cdot,t)\|_{C^{\theta_2}(\bom)} <\infty,
  \eas
  so that since 
  $\|B\varphi\|_{C^{\theta_2}(\bom)} \le \|B\|_{L^\infty(\Omega)} \|\varphi\|_{C^{\theta_2}(\bom)}
  + \|B\|_{C^{\theta_2}(\bom)} \|\varphi\|_{L^\infty(\Omega)}$ for all $B\in C^{\theta_2}(\bom)$ and $\varphi\in C^{\theta_2}(\bom)$,
  from \eqref{211.1} we infer that with some $c_1>0$ we have
  \bea{above-2.106}
	\|\na w_t\|_{C^{\theta_2}(\bom)}
	&\le& 
	\|B_1\na a\|_{C^{\theta_2}(\bom)} 
	+ \|B_2\na w\|_{C^{\theta_2}(\bom)} 
	+ \|B_3\na v\|_{C^{\theta_2}(\bom)} \nn\\
	&\le& c_1 + c_1 \|\na w\|_{C^{\theta_2}(\bom)}
	\qquad \mbox{for all } t\in (0,\tm).
  \eea
  Therefore,
  \bas
	\|\na w(\cdot,t)\|_{C^{\theta_2}(\bom)}
	\le \|\na w_0\|_{C^{\theta_2}(\bom)}
	+ c_1\tm + c_1 \int_0^t \|\na w(\cdot,s)\|_{C^{\theta_2}(\bom)} ds
	\qquad \mbox{for all } t\in (0,\tm),
  \eas
  so that Gronwall's lemma ensures the existence of $c_2>0$ fulfilling
  \be{24.2}
	\|\na w(\cdot,t)\|_{C^{\theta_2}(\bom)} \le c_2
	\qquad \mbox{for all } t\in (0,\tm).
  \ee
  As a consequence of this, \eqref{above-2.106} now guarantees that also $c_3:=\|\na w_t\|_{L^\infty(\Omega\times (0,\tm))}$ is finite,
  which together with \eqref{24.2} shows that
  \bas
	|\na w(x_1,t_1)-\na w(x_2,t_2)|
	&\le& |\na w(x_1,t_1)-\na w(x_1,t_2)| + |\na w(x_1,t_2)-\na w(x_2,t_2)| \\
	&\le& c_3 |t_1-t_2| + c_2|x_1-x_2|^{\theta_2} \\
	& & \hs{14mm} \mbox{for all $(x_1,t_1)\in\Omega\times (0,\tm)$ and } (x_2,t_2)\in\Omega\times (0,\tm).
  \eas
  Along with a similar argument directly applied to $w$ on the basis of \eqref{24.11} and the evident boundedness of $w_t$ in
  $\Omega\times (0,\tm)$, this asserts \eqref{24.1}.
\qed\\[-2ex]

\noindent
This now warrants accessibility of the first, second and fourth equations from \eqref{0a} to standard parabolic Schauder theory:
\begin{Lemma} \label{lem25}
  If $\tm<\infty$, then there exists $\theta\in (0,1)$ such that
  \be{25.1}
	a\in C^{2+\theta,1+\frac{\theta}{2}}(\bom\times [0,\tm])
  \ee
  and
  \be{25.2}
	v\in C^{2+\theta,1+\frac{\theta}{2}}(\bom\times [0,\tm])
  \ee
  as well as
  \be{25.3}
	z\in C^{2+\theta,1+\frac{\theta}{2}}(\bom\times [0,\tm]).
  \ee
\end{Lemma} 
\proof
  From \cref{lem9} together with \cref{lem23} and \cref{lem24} we know that in the identity 
  $z_t=\Delta z+g_1(x,t) z + g_2(x,t)$, the functions $g_1:=-v$ and $g_2:=\phi(u) \equiv \phi(ae^{\Xi(v,w)})$ belong to
  $C^{\theta_1,\frac{\theta_1}{2}}(\bom\times [0,\tm)$ with some $\theta_1\in (0,1)$.
  
  \noindent
  Therefore, parabolic Schauder estimates (\cite{LSU}) provide $\theta_2\in (0,1)$ fulfilling
  $z\in C^{2+\theta_2,1+\frac{\theta_2}{2}}(\bom\times [0,\tm])$, so that, in particular, the coefficients $h_1:=-\na z$ and
  $h_2:=-\Delta z$ in $v_t=\Delta v + h_1(x,t)\cdot\na v + h_2(x,t)v$ are H\"older continuous in $\bom\times [0,\tm]$,
  and that hence, again by Schauder theory, $v$ must be an element of $C^{2+\theta_3,1+\frac{\theta_3}{2}}(\bom\times [0,\tm])$
  for some $\theta_3\in (0,1)$.
  
    \noindent
  Together with the outcomes of \cref{lem23} and \cref{lem24}, through \eqref{111.1} and \eqref{A} this in turn warrants
  H\"older continuity of $a_t-\Delta a$ in $\bom\times [0,\tm]$, whence the proof can be completed by a third application
  of parabolic Schauder theory.
\qed\\[-2ex]

\noindent
In consequence, due to \eqref{211.1} also $w$ must remain bounded in $C^{2+\theta}(\bom)$ with some $\theta\in (0,1)$
if $\tm<\infty$:
\begin{Lemma} \label{lem26}
  If $\tm<\infty$, then there exist $\theta\in (0,1)$ and $C>0$ such that
  \be{26.1}
	\|w(\cdot,t)\|_{C^{2+\theta}(\bom)} \le C
	\qquad \mbox{for all } t\in (0,\tm).
  \ee
\end{Lemma} 
\proof
  This can be seen by a straightforward adaptation of the argument from \cref{lem24}, based on the result from \cref{lem25}
  and a differentiation of the identity in \cref{211.1}.
\qed
\medskip

\noindent
In light of \cref{lem_loc}, our main result thereby becomes obvious:\abs
\proofc of \cref{theo27}. \quad
  As a consequence of \cref{lem25} and \cref{lem26}, assuming $\tm$ to be finite would 
  yield a contradiction to the extensibility criterion \cref{ext}, so that the claim results from  \cref{lem_loc}.
\qed
\section{Numerical simulations}\label{sec:numerics}

We perform numerical simulations of system \eqref{q-macro-u}, \eqref{q-macro-v}, \eqref{-1p2}, together with the initial and no-flux boundary conditions from \eqref{0}. For the discretization we use a finite difference approach. We consider a domain $\Omega =[0,100]\times [0,100]$ (in $\mu m$) and a regular lattice. The standard scheme is improved using for the equations concerning cell motility the method in \cite{weickert1998anisotropic}, in order to ensure preservation of non-negativity for the solution. The advection terms are discretized with a first order upwind scheme. The time discretization uses an implicit-explicit IMEX method (see e.g., \cite{Ascher1995}) which handles the diffusion implicitly, while the drift and reaction terms are treated in an explicit manner. The implementation was performed in MATLAB\footnote{Version 9.11.0.1769968 (R2021b), The MathWorks Inc.}.\\[-2ex]

\noindent
For the functions $\phi, \psi$ in \eqref{-1p2} we make the following choices, which satisfy \eqref{psi}, \eqref{phi}:
\begin{align*}
\phi (u):=\frac{u}{1+u},\quad \psi(s):=\frac{\beta s}{1+s},
\end{align*}
where $K_u>0$ denotes the tumor carrying capacity and $\beta>0$ is a constant.\\[-2ex]

\begin{figure}[h!]
	\includegraphics[width=0.24\linewidth]{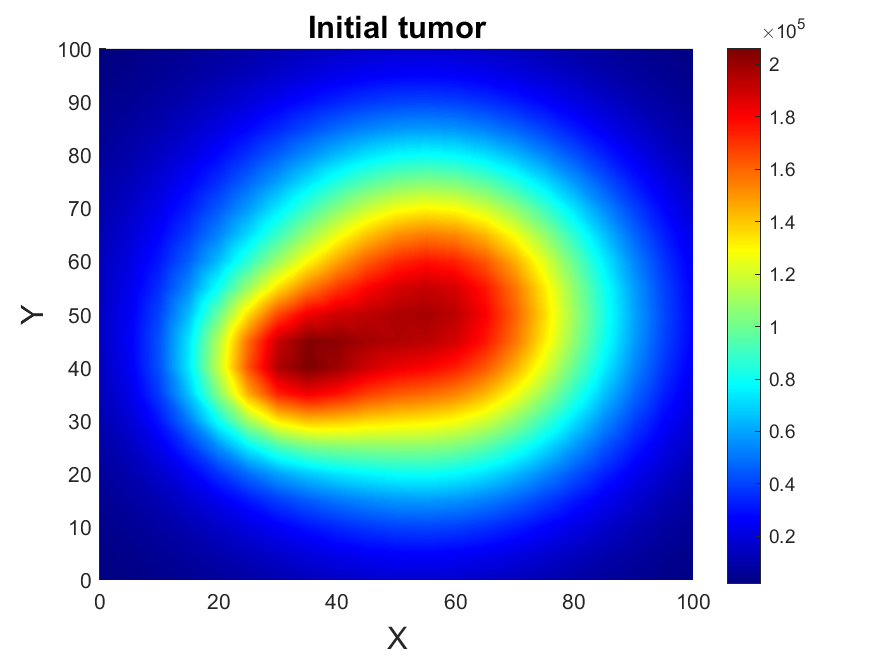}\hspace*{-0.3cm} \includegraphics[width=0.24\linewidth]{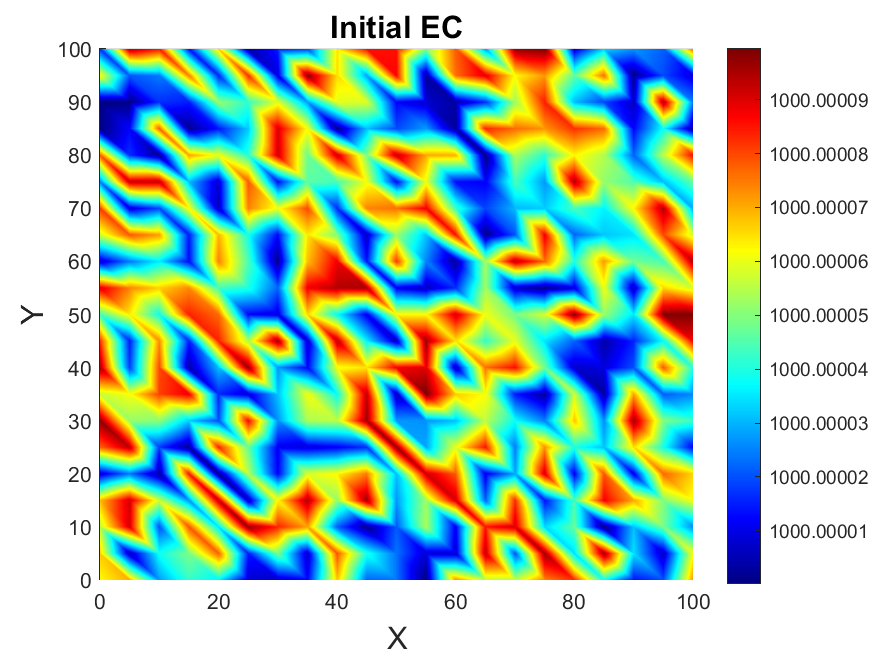}\quad \includegraphics[width=0.24\linewidth]{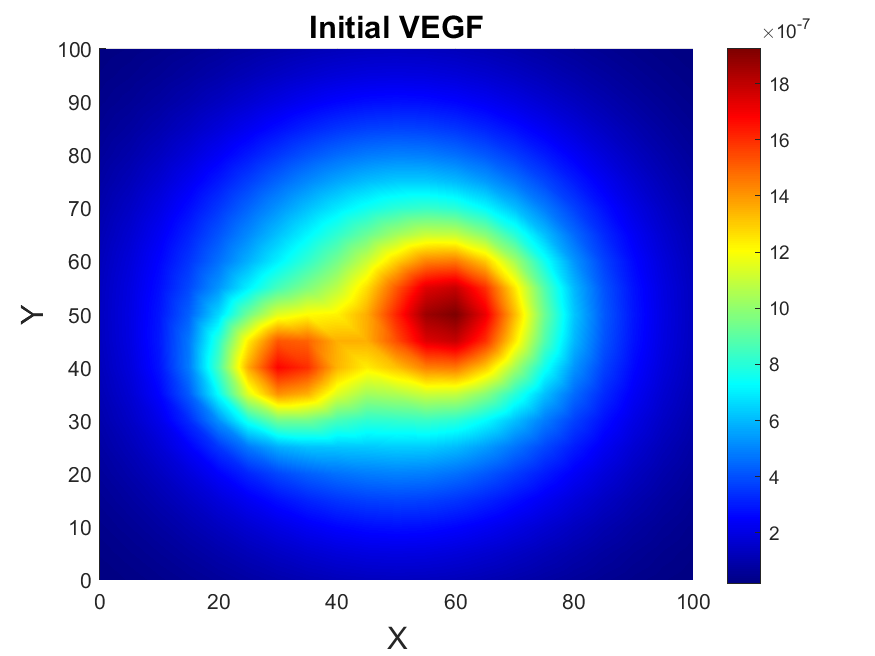}\hspace*{-0.3cm} \includegraphics[width=0.24\linewidth]{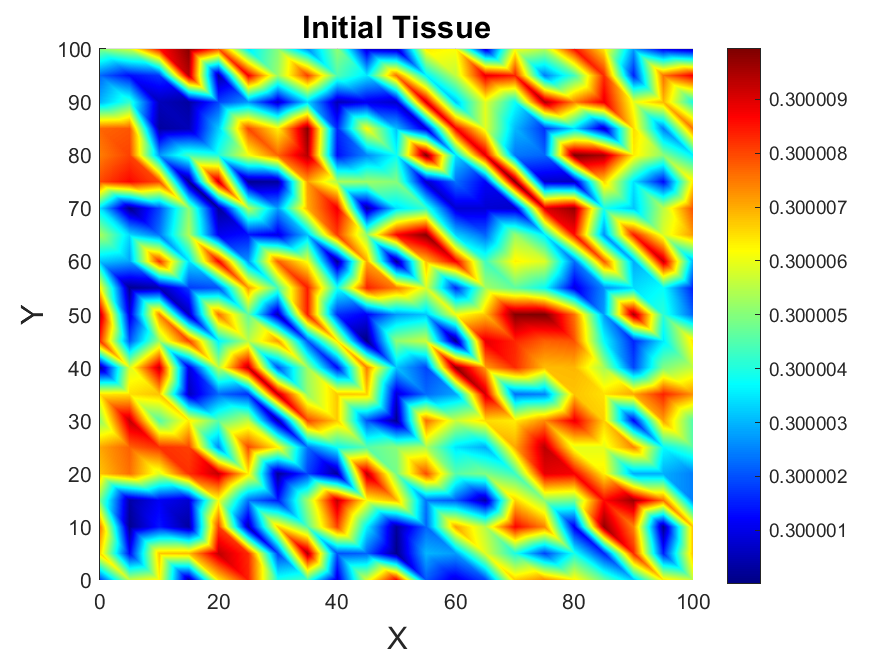}\\
	\includegraphics[width=0.24\linewidth]{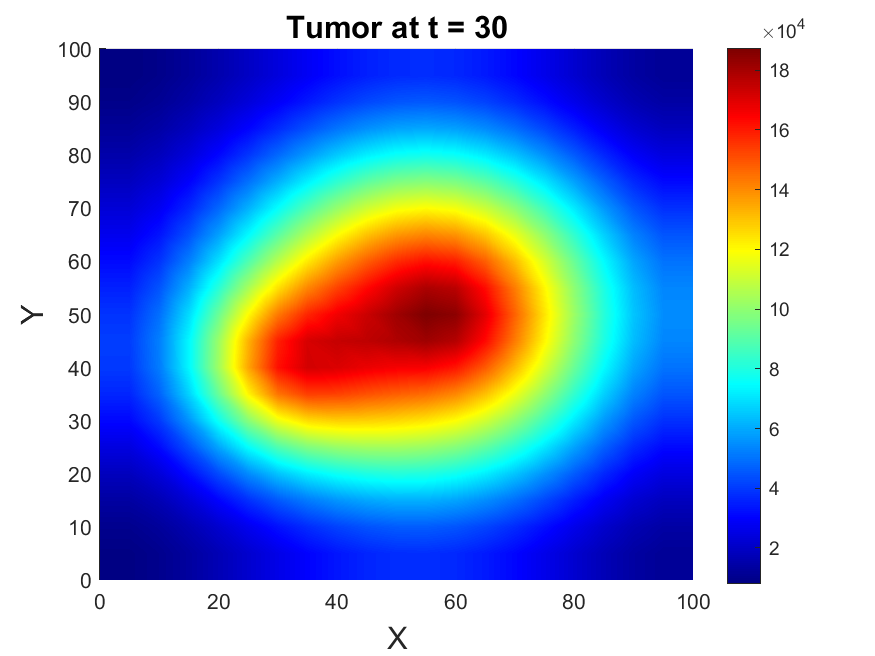}\hspace*{-0.3cm} \includegraphics[width=0.24\linewidth]{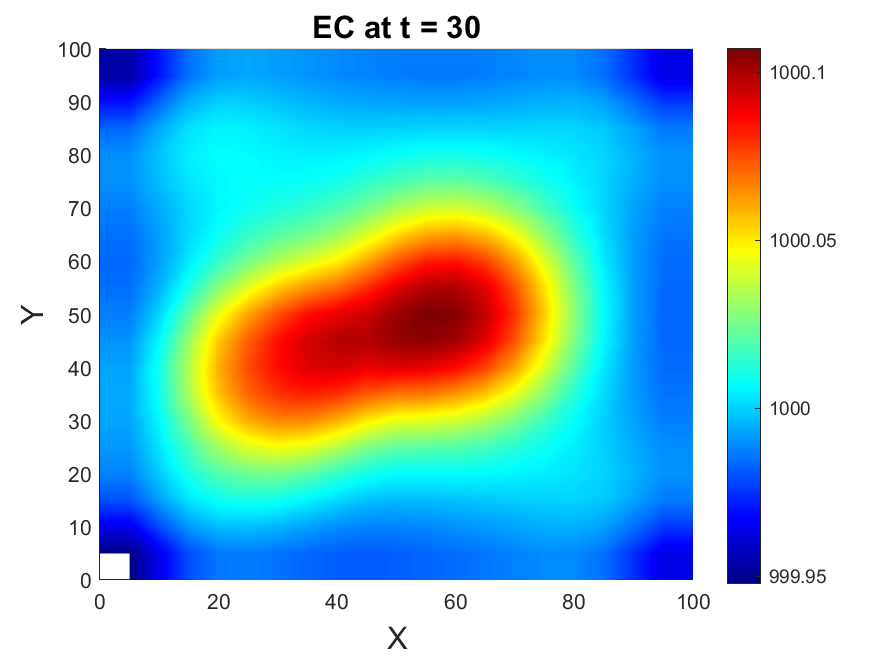}\quad \includegraphics[width=0.24\linewidth]{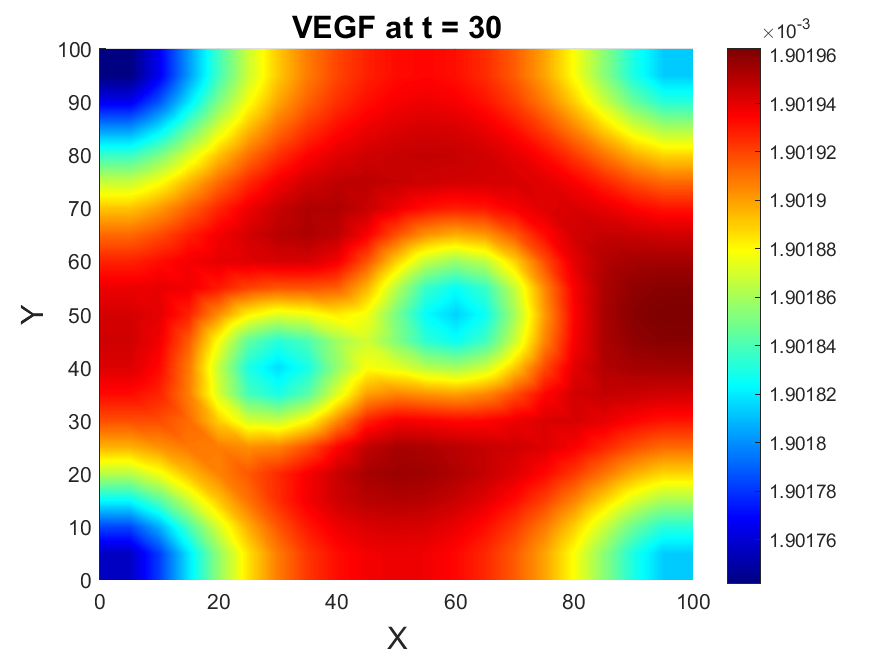}\hspace*{-0.3cm} \includegraphics[width=0.24\linewidth]{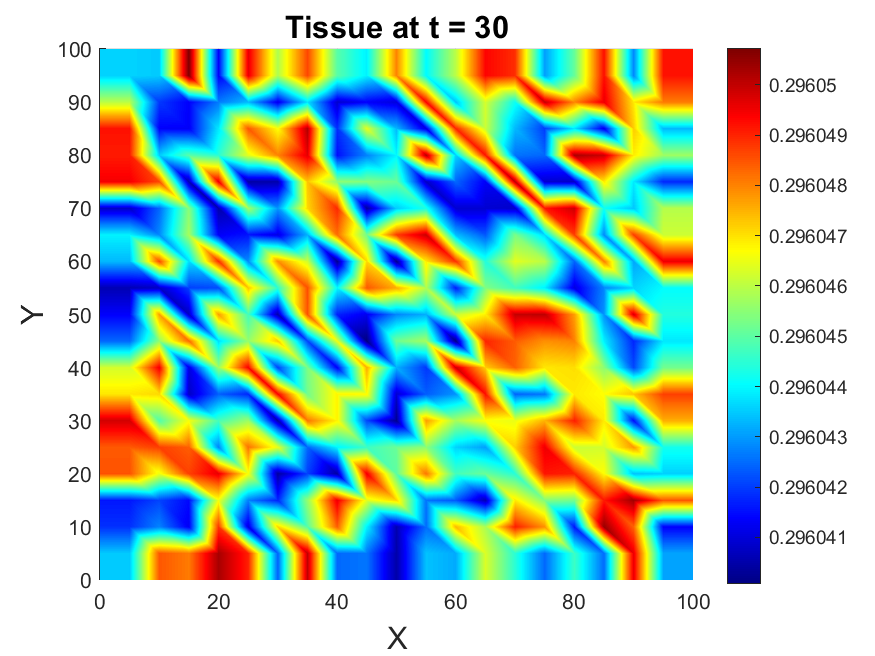}\\
	\includegraphics[width=0.24\linewidth]{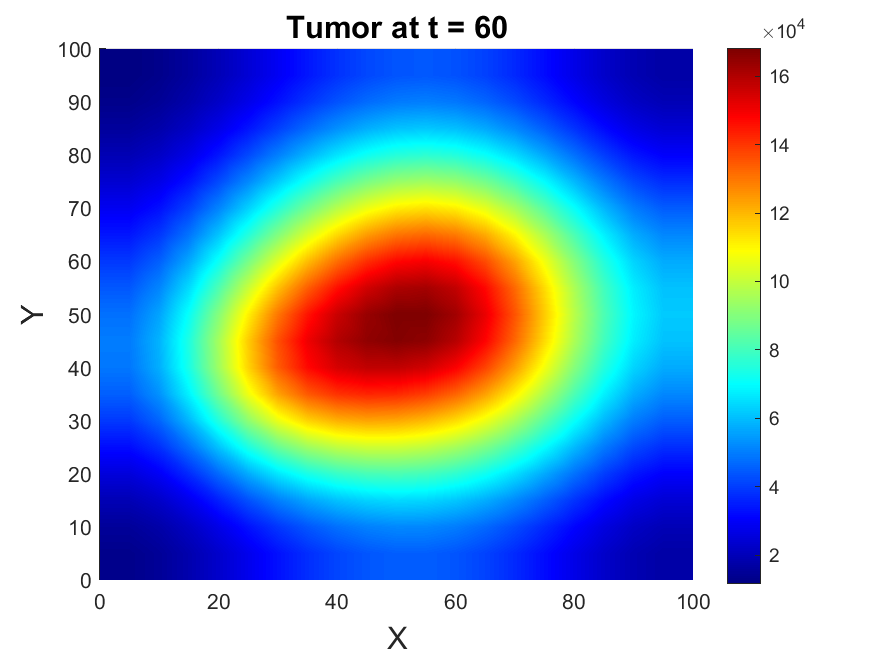}\hspace*{-0.3cm} \includegraphics[width=0.24\linewidth]{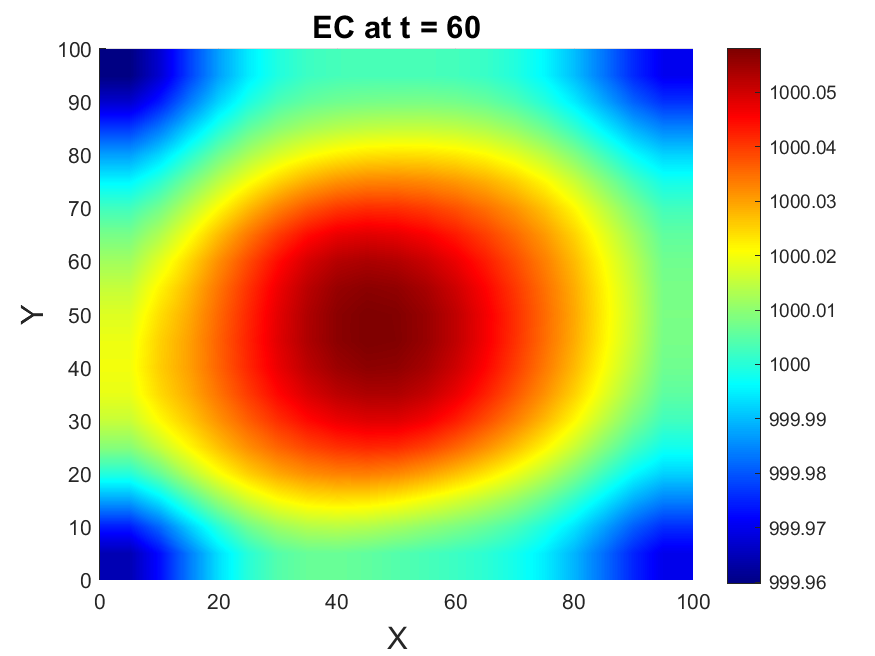}\quad \includegraphics[width=0.24\linewidth]{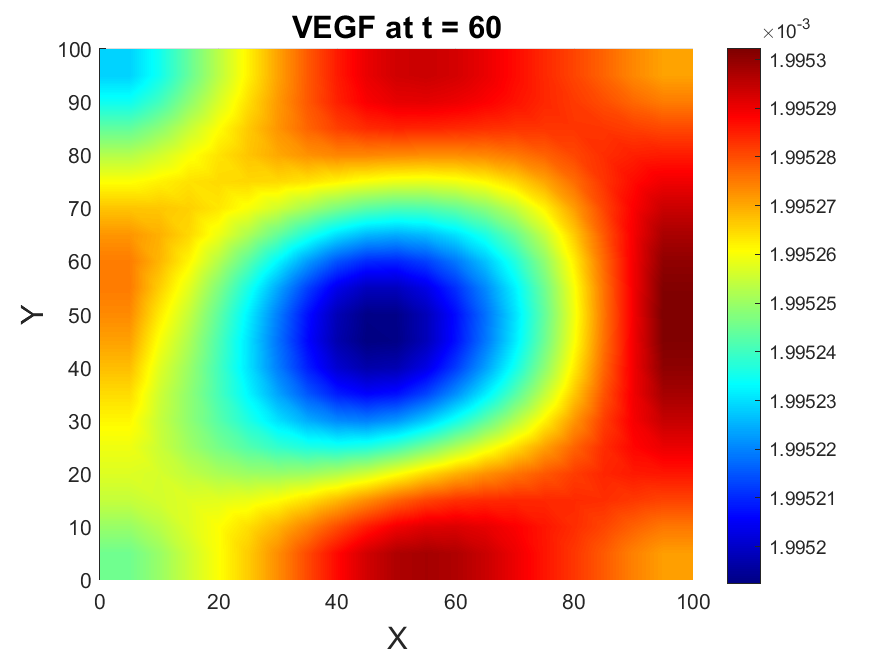}\hspace*{-0.3cm} \includegraphics[width=0.24\linewidth]{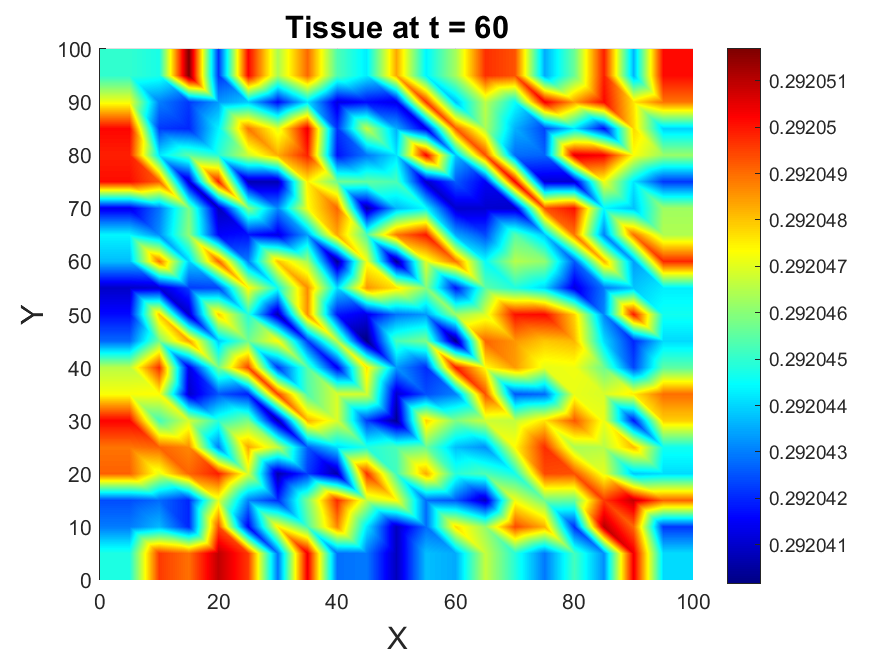}\\
	\includegraphics[width=0.24\linewidth]{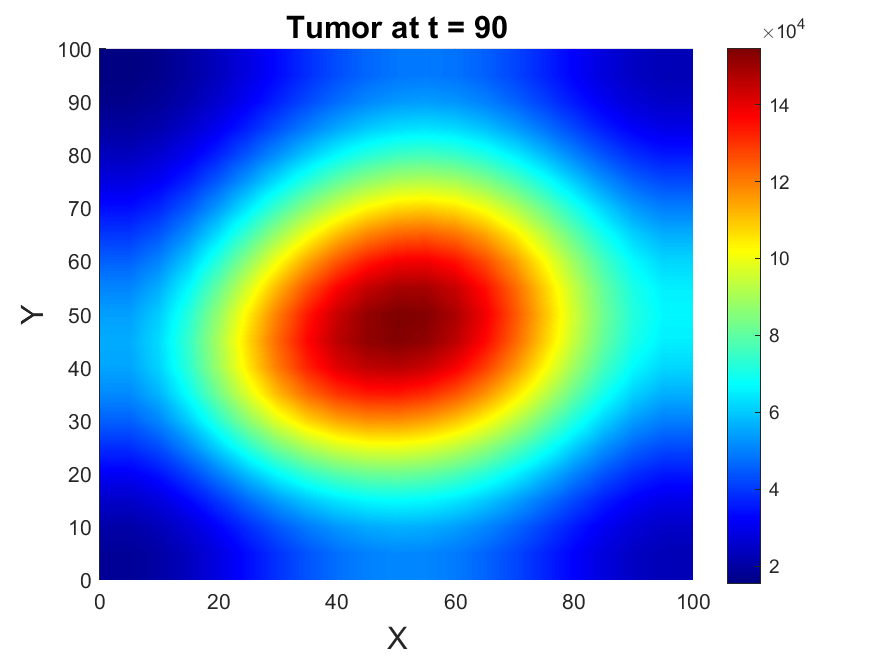}\hspace*{-0.3cm} \includegraphics[width=0.24\linewidth]{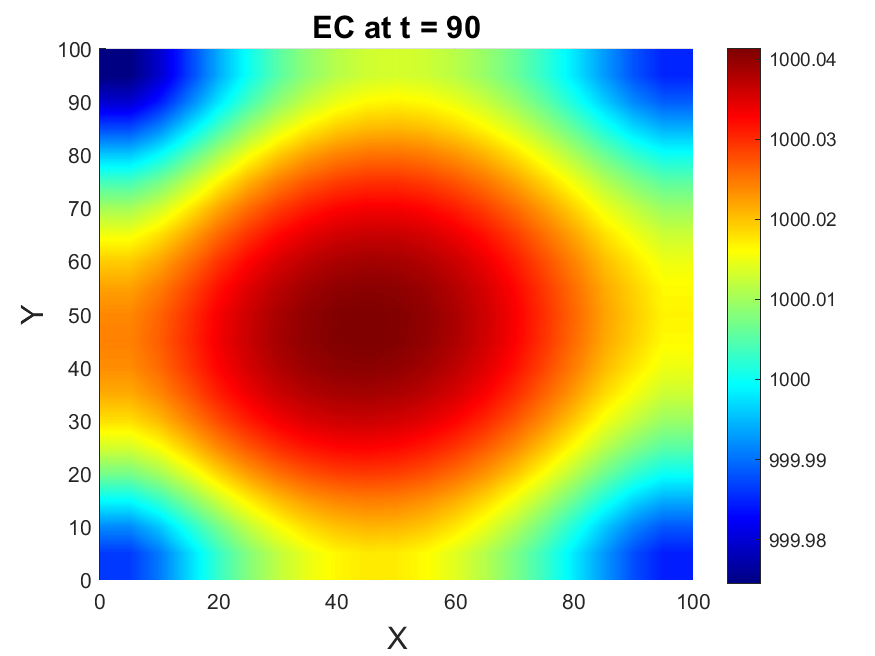}\quad \includegraphics[width=0.24\linewidth]{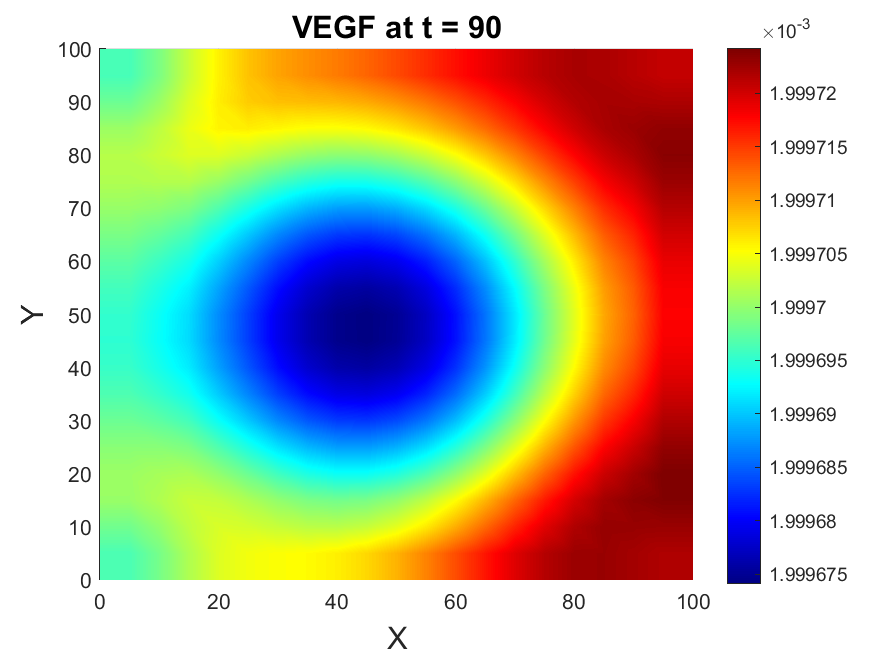}\hspace*{-0.3cm} \includegraphics[width=0.24\linewidth]{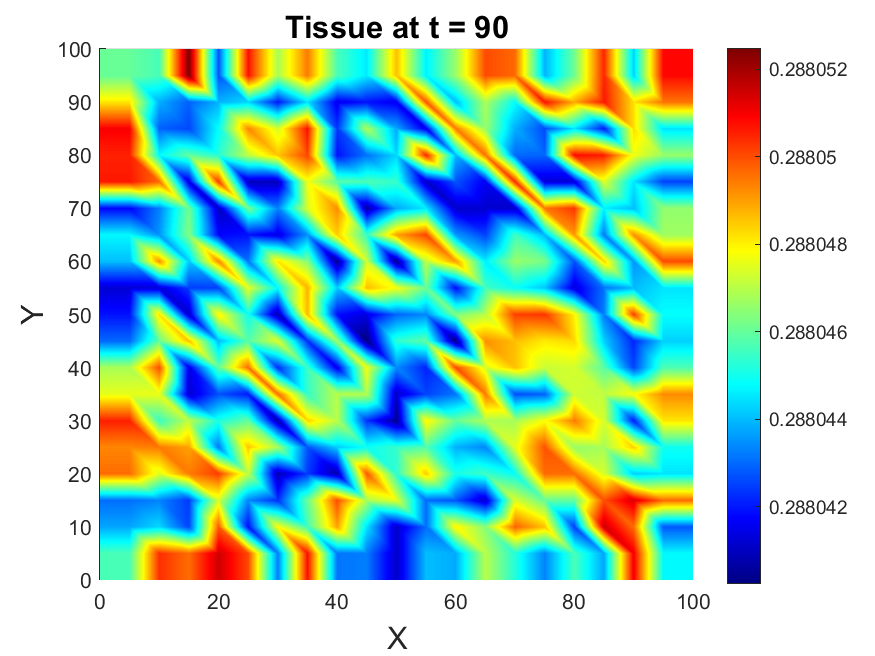}\\
	\includegraphics[width=0.24\linewidth]{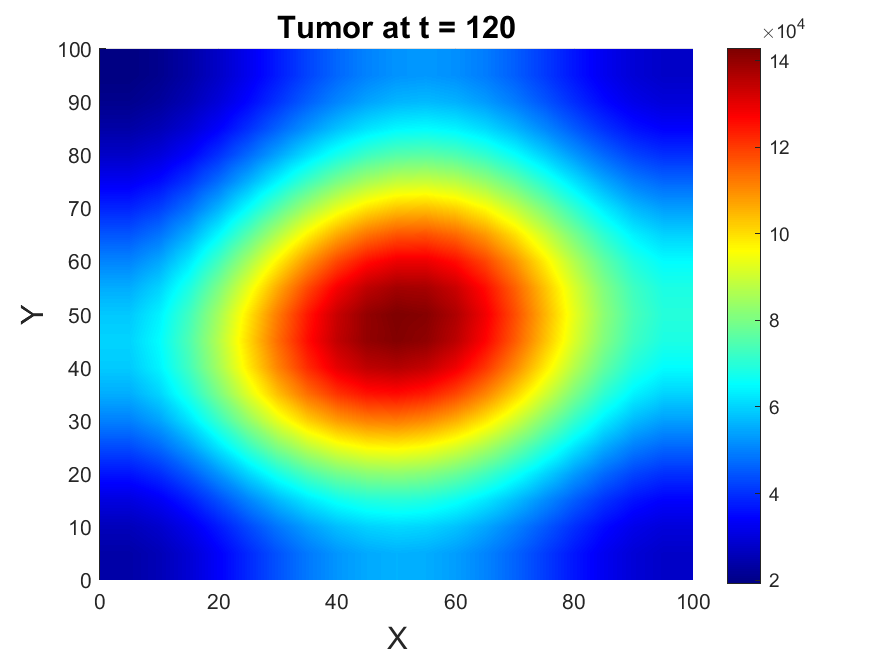}\hspace*{-0.3cm} \includegraphics[width=0.24\linewidth]{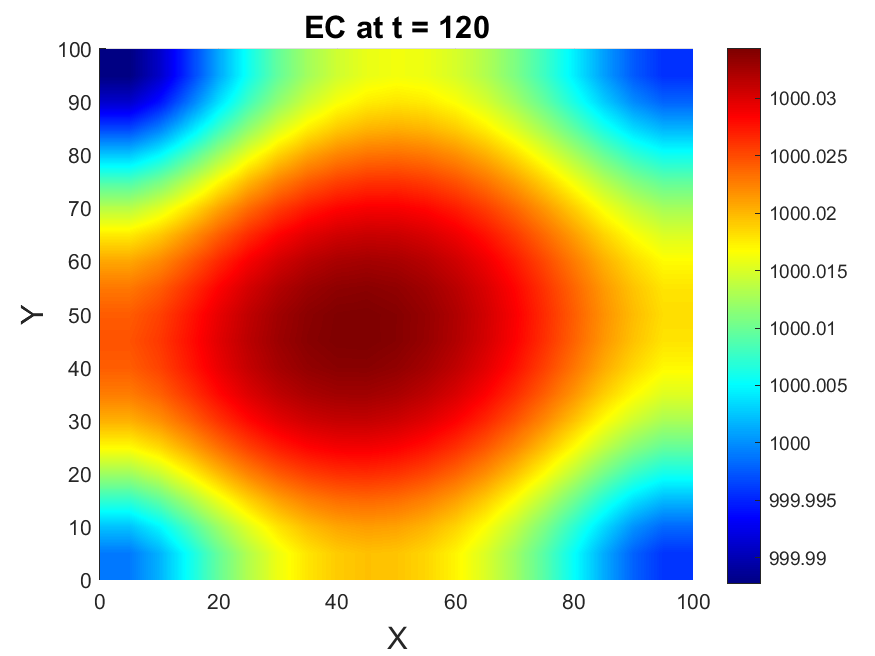}\quad \includegraphics[width=0.24\linewidth]{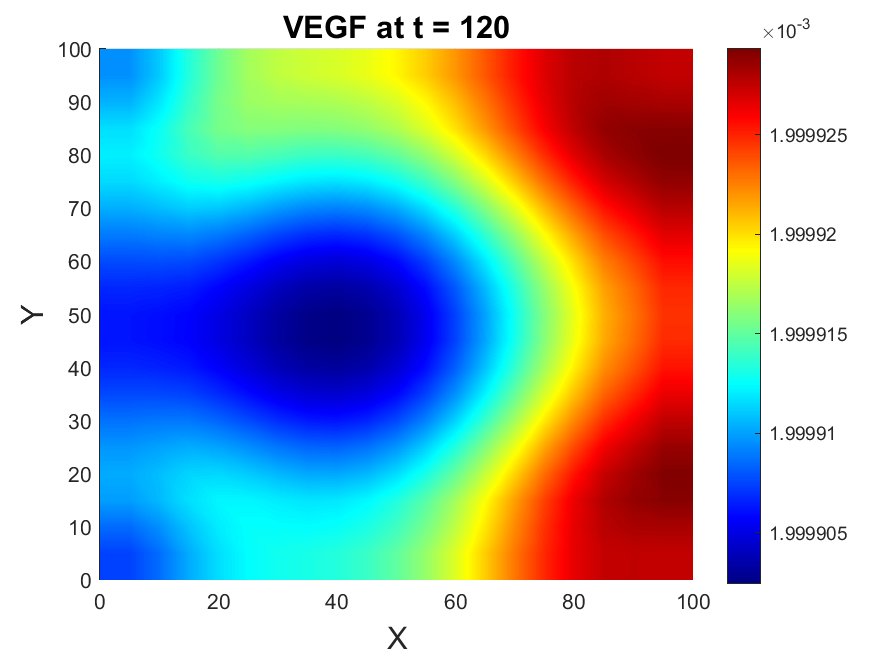}\hspace*{-0.3cm} \includegraphics[width=0.24\linewidth]{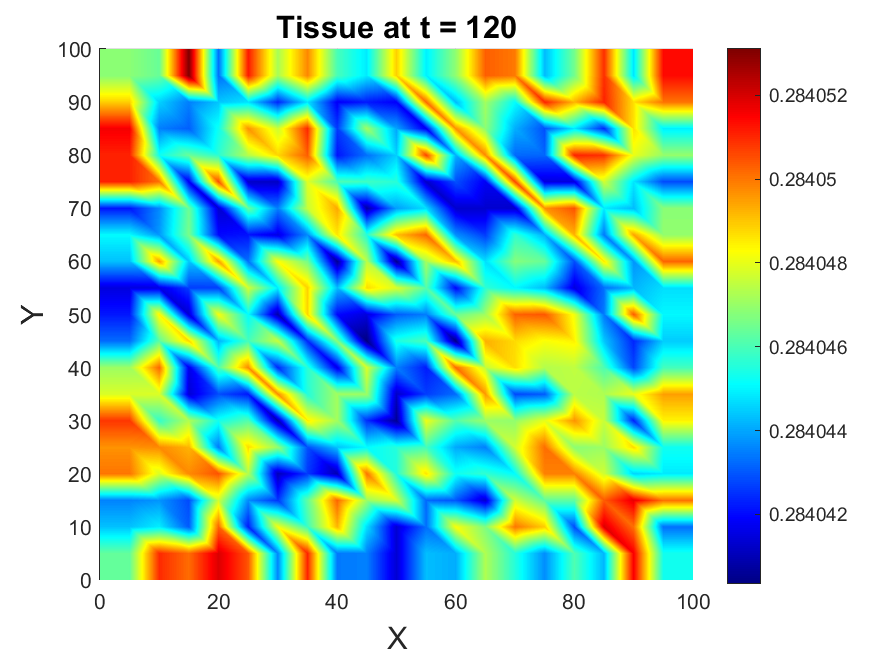}\\
	\includegraphics[width=0.24\linewidth]{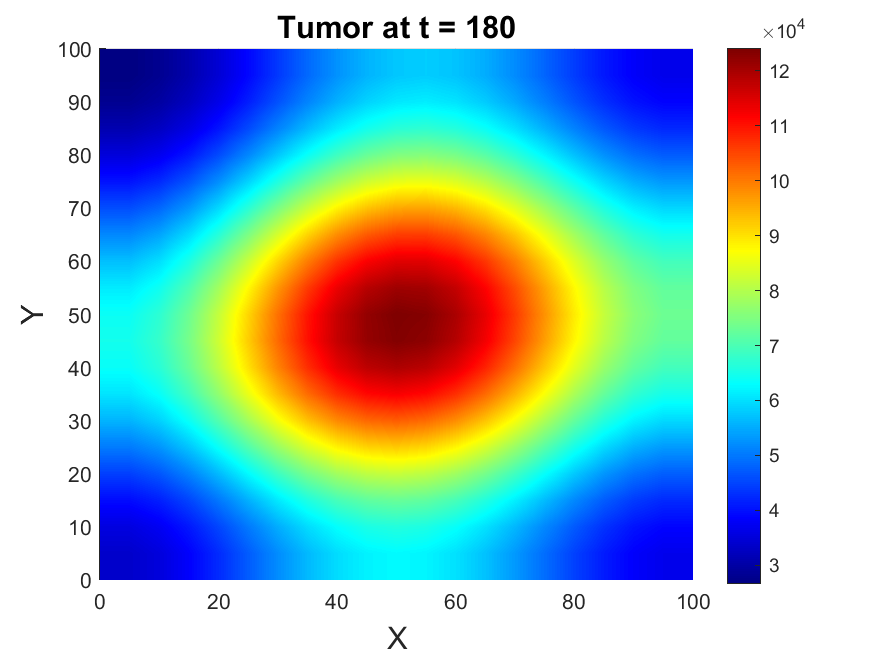}\hspace*{-0.3cm} \includegraphics[width=0.24\linewidth]{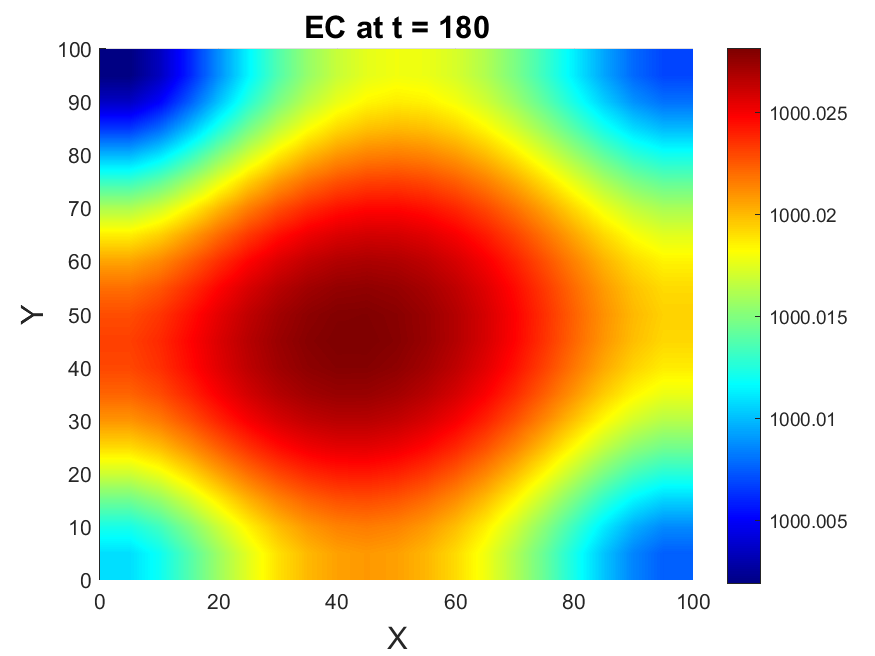}\quad \includegraphics[width=0.24\linewidth]{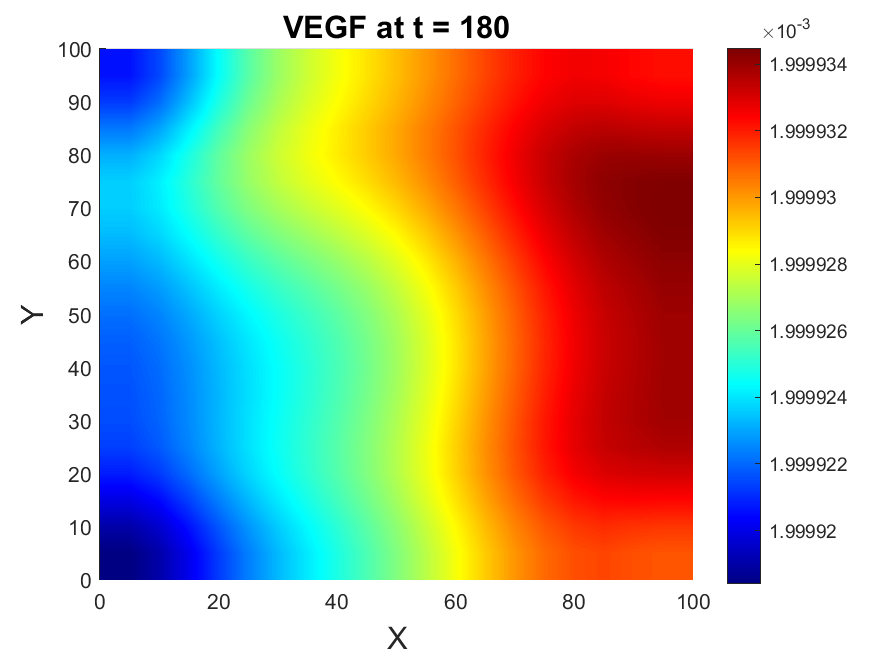}\hspace*{-0.3cm} \includegraphics[width=0.24\linewidth]{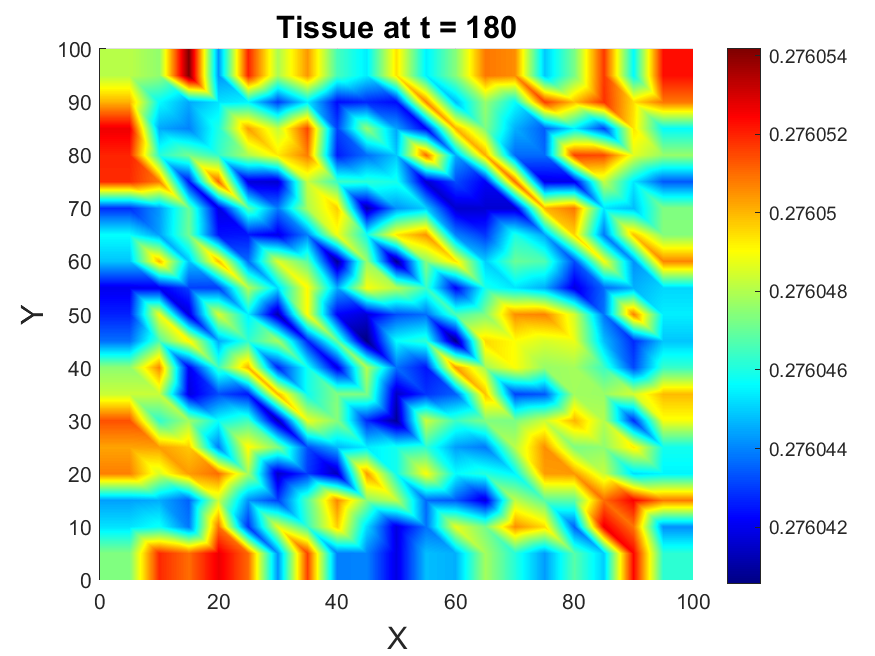}\\
	\includegraphics[width=0.24\linewidth]{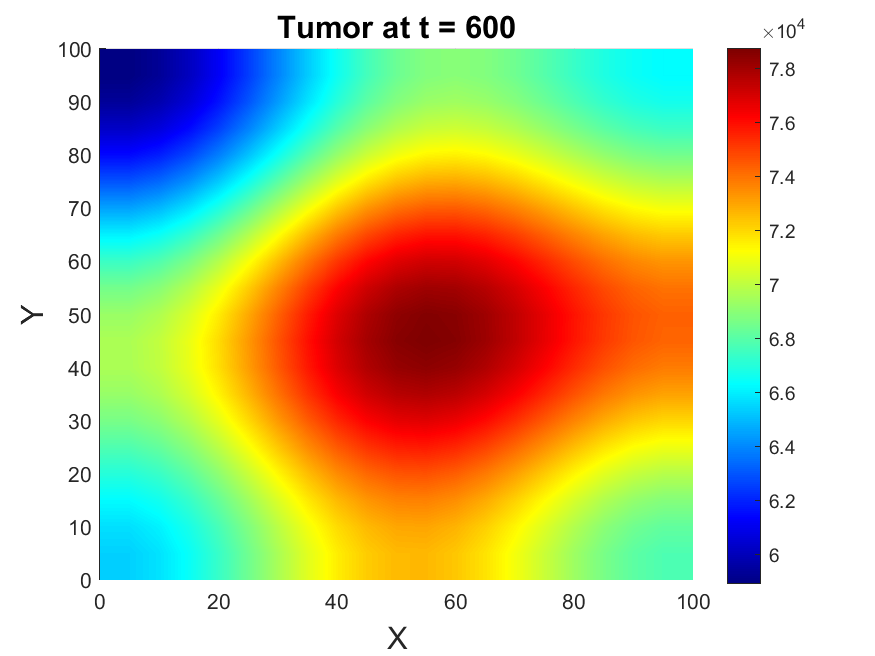}\hspace*{-0.3cm} \includegraphics[width=0.24\linewidth]{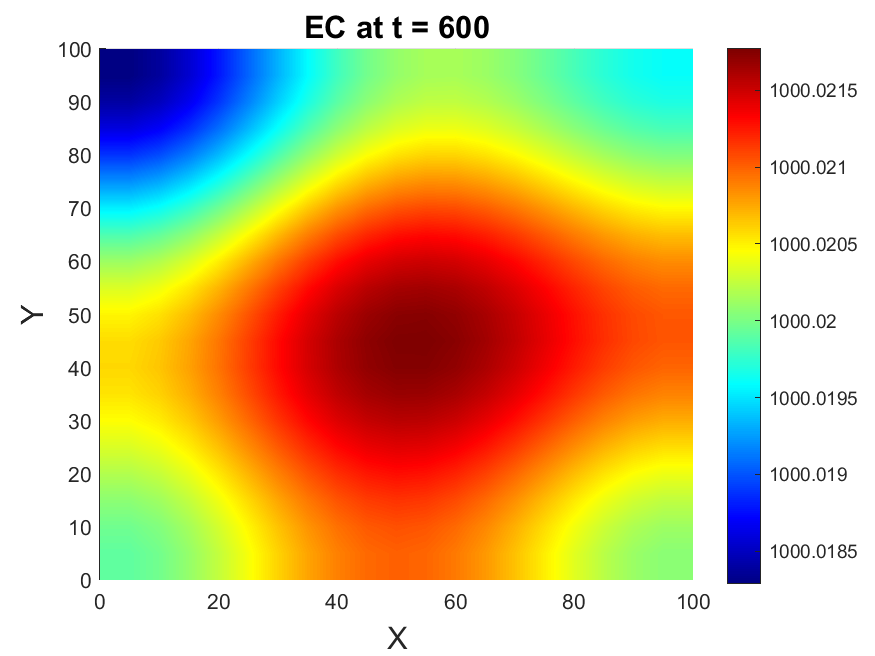}\quad \includegraphics[width=0.24\linewidth]{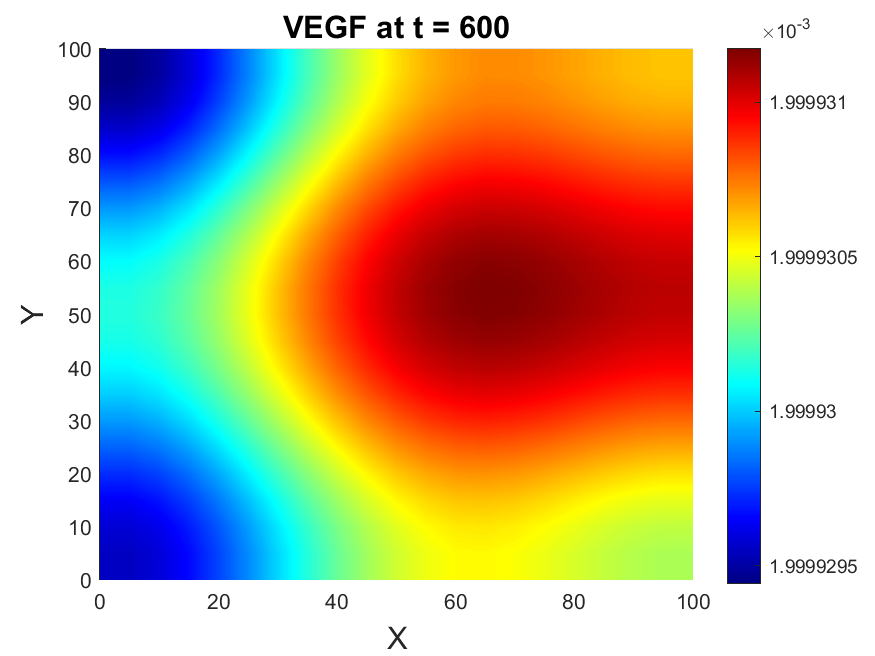}\hspace*{-0.3cm} \includegraphics[width=0.24\linewidth]{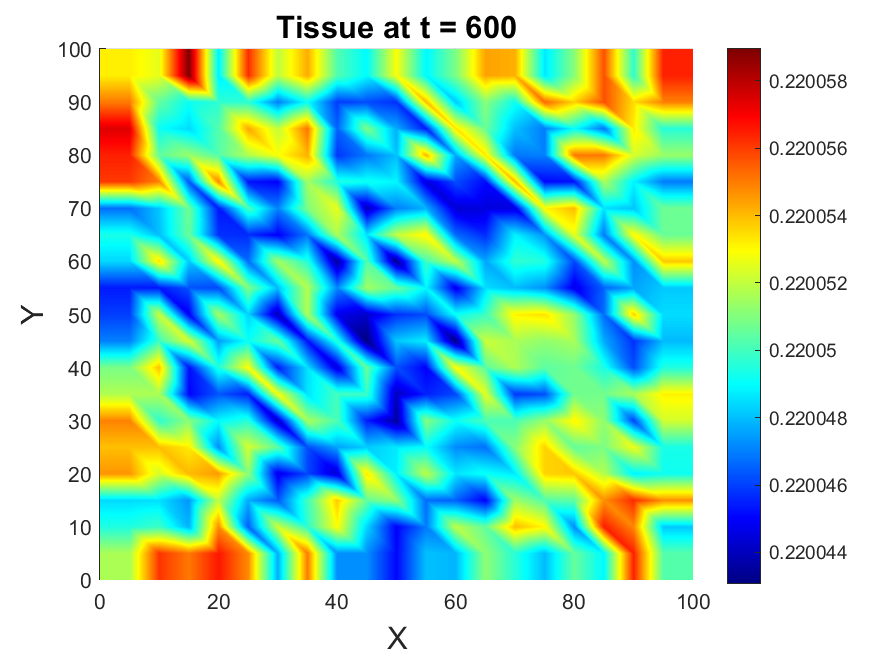}\\
	\includegraphics[width=0.24\linewidth]{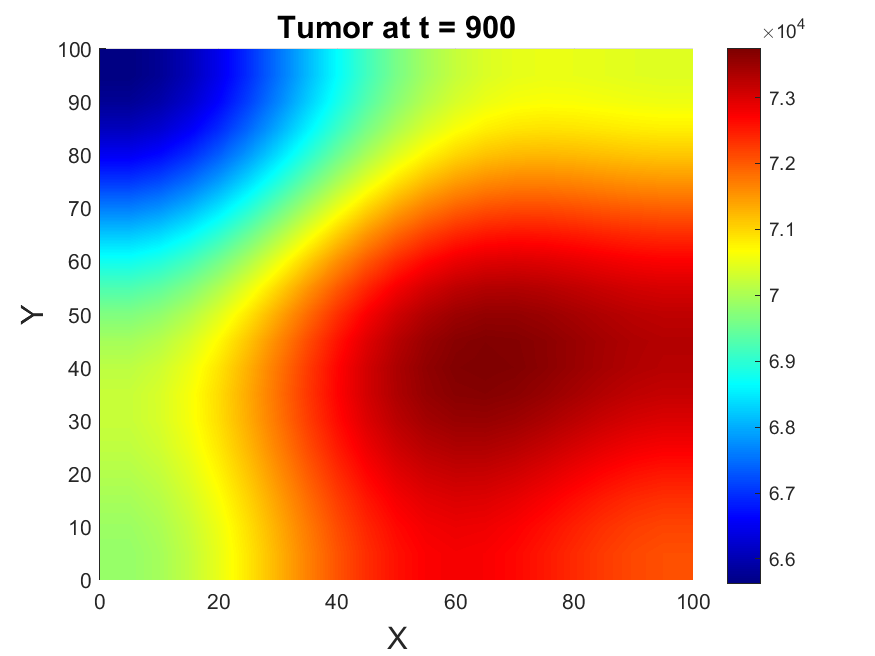}\hspace*{-0.3cm} \includegraphics[width=0.24\linewidth]{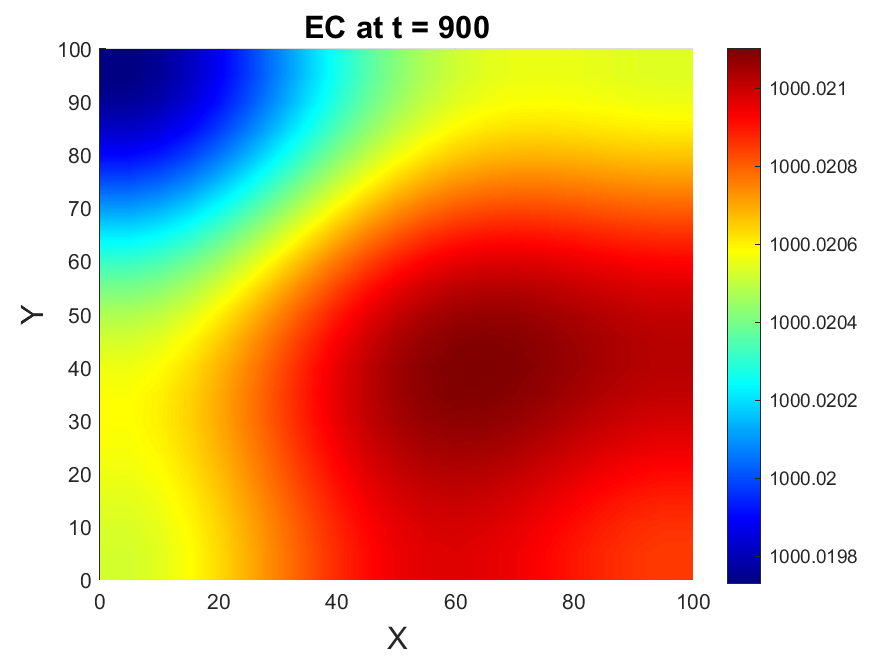}\quad \includegraphics[width=0.24\linewidth]{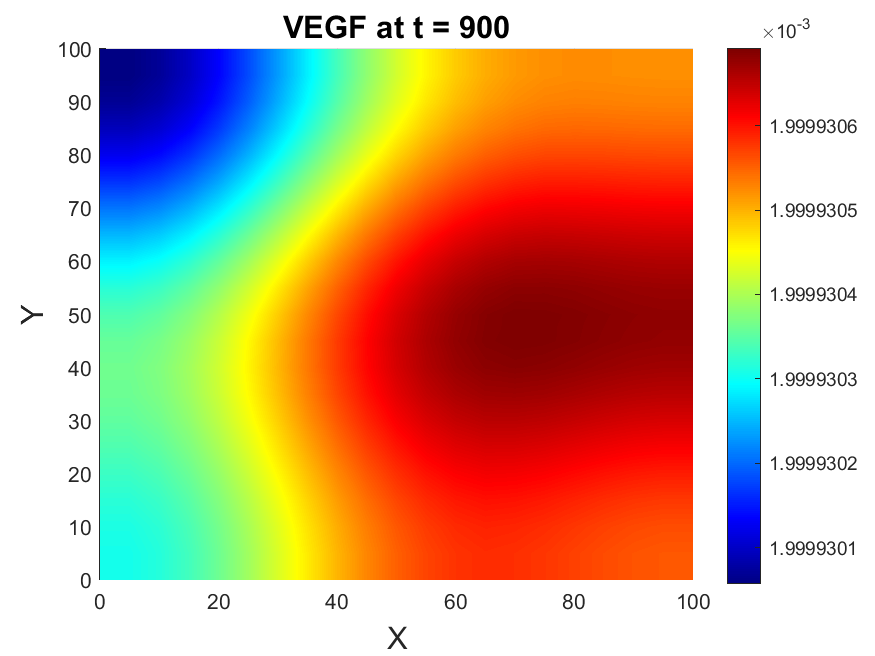}\hspace*{-0.3cm} \includegraphics[width=0.24\linewidth]{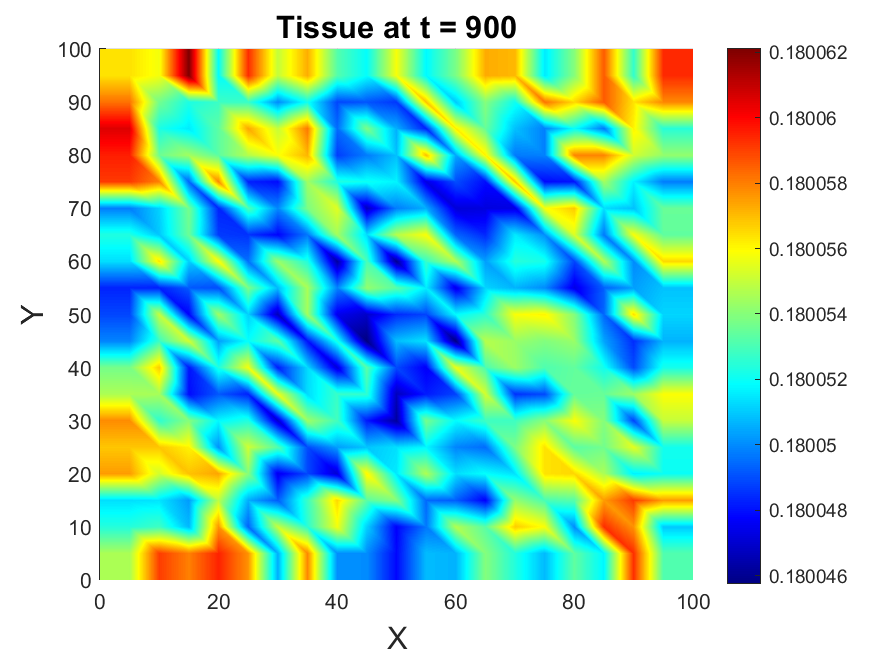}\\
			\caption{Simulations of model \eqref{q-macro-u}-\eqref{-1p2} without source terms for tumor and endothelial cells.}\label{fig:1}
\end{figure}

\noindent
First we run the simulations for the model \eqref{q-macro-u}-\eqref{-1p2} without source terms for tumor and endothelial cells: this is the situation for which the global existence and uniqueness proof was performed. The initial conditions are as illustrated in the first row of Figure \ref{fig:1}. Thus, ECs and tissue are assumed to be uniformly distributed within $\Omega $, while the tumor and, correspondigly, the VEGF distribution are more 'localized'. The rest of Figure \ref{fig:1} shows the densities $u$, $v$, $z$, $w$ of tumor cells (1st column), ECs (2nd column), VEGF (3rd column), and tissue (last column), respectively. The predicted behavior of the solution components is as expected from the corresponding biological processes: in the absence of source terms, the tumor cells spread over the domain, following EC and tissue gradients. They release VEGF, which is relatively fast diffusing. The ECs are chemotactically attracted by VEGF; at the beginning of the simulation they are quickly grouping at areas with abundant VEGF (the latter being produced by the tumor cells), then spread due to diffusion and taxis. The two middle columns show how the EC population is following (with some delay) the VEGF bulk and also the uptake of VEGF at the sites with high EC density. The ECM is degraded by the tumor cells, more substantially at the sites where the tumor bulk is more concentrated, and weaker when the cancer cells have spread and filled the space at densities which are not high enough to noticeably degrade the tissue any further. Tumor, ECs, and VEGF eventually accumulate in the lower right corner, the tissue in the upper left corner is correspondingly spared, and the dynamics is not changing anymore. In fact, all solutions components seem to remain within numerical ranges which are close to certain values. This suggests an asymptotic stabilization of the solution around those values. \\[-2ex]

\begin{figure}[h!]
	\includegraphics[width=0.24\linewidth]{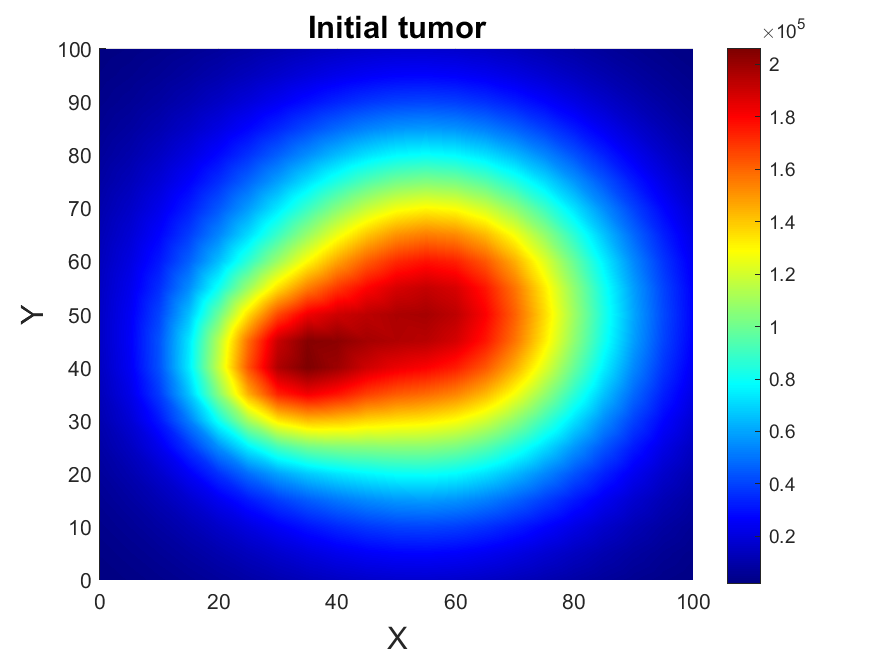}\hspace*{-0.3cm} \includegraphics[width=0.24\linewidth]{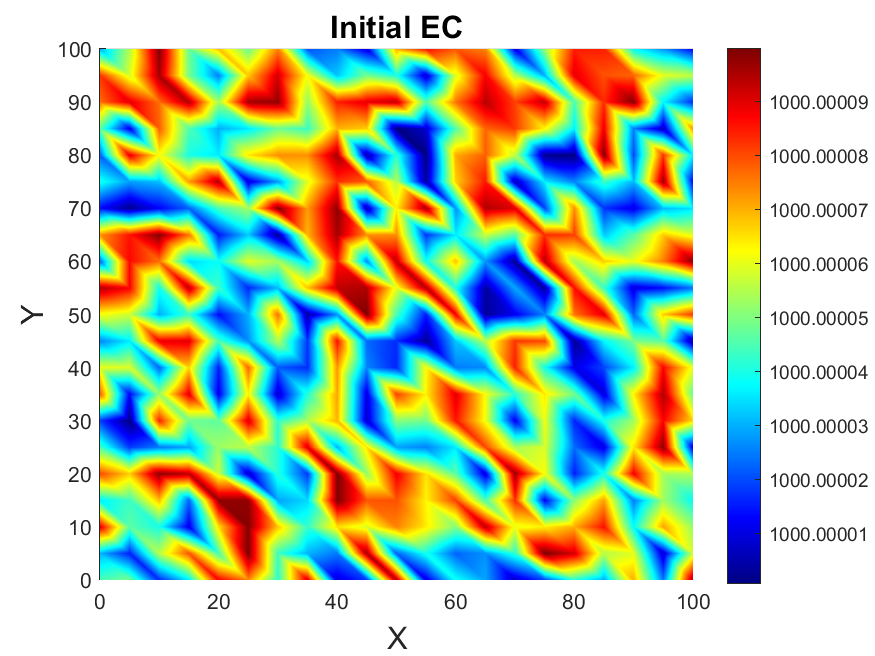}\quad \includegraphics[width=0.24\linewidth]{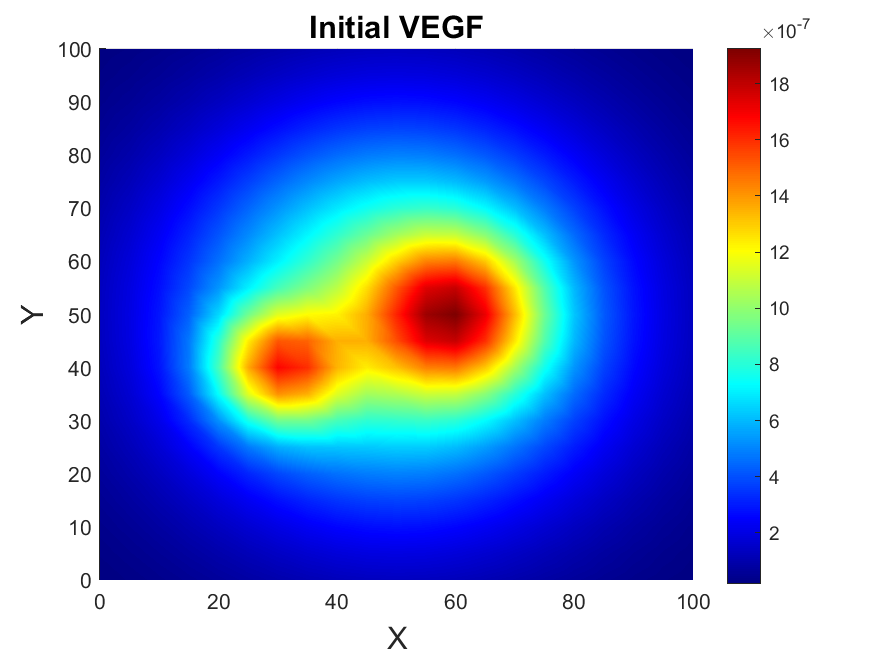}\hspace*{-0.3cm} \includegraphics[width=0.24\linewidth]{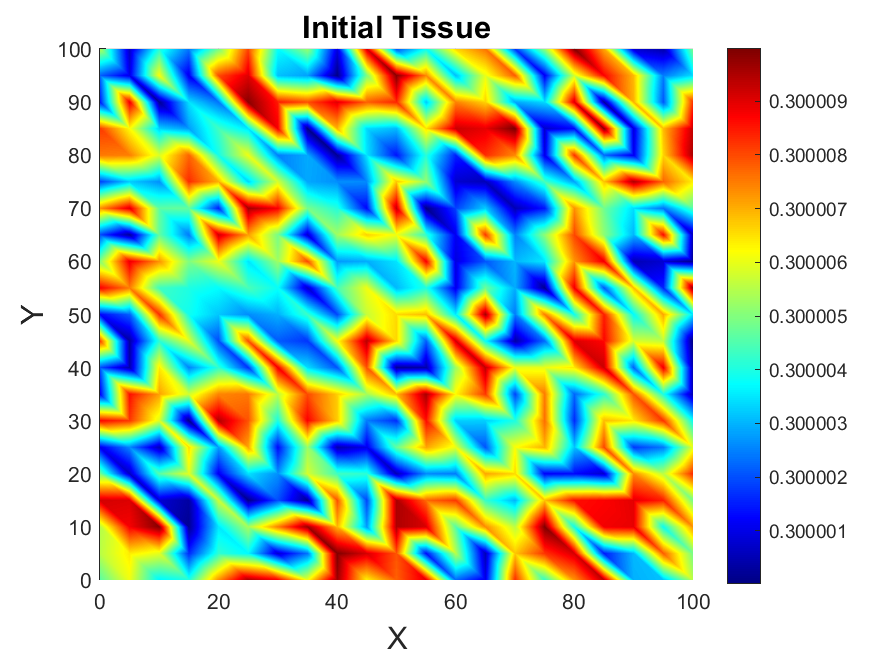}\\
	\includegraphics[width=0.24\linewidth]{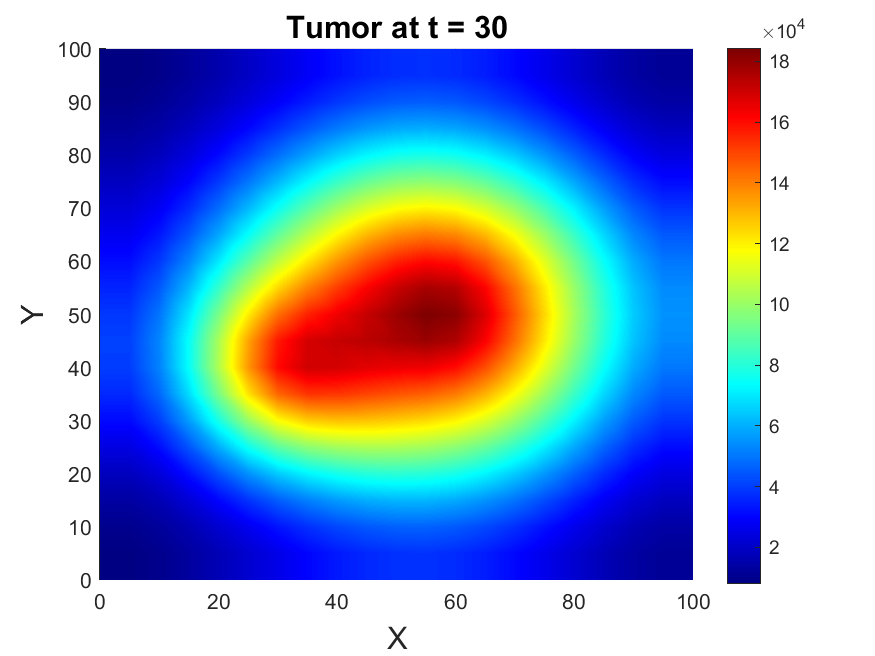}\hspace*{-0.3cm} \includegraphics[width=0.24\linewidth]{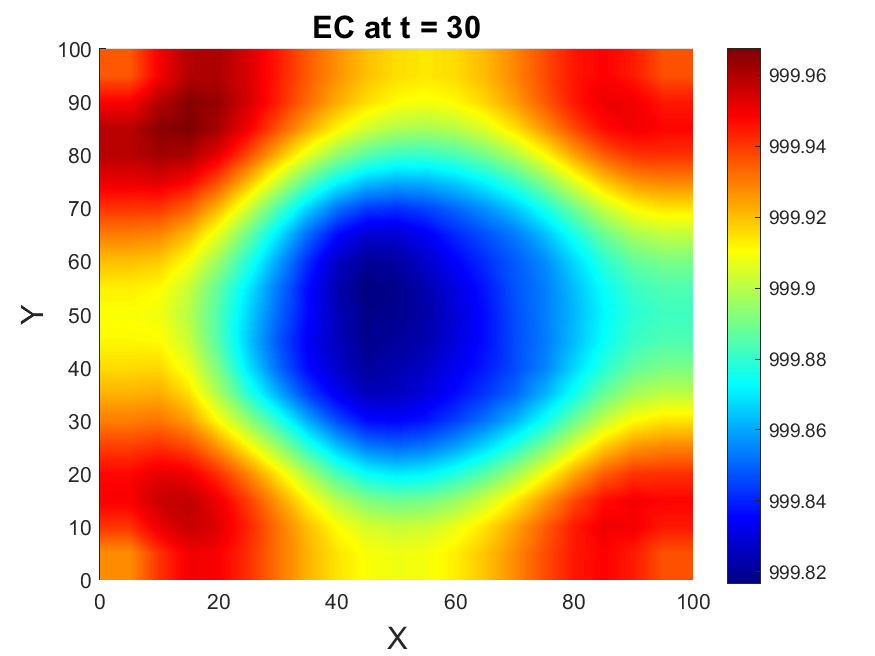}\quad \includegraphics[width=0.24\linewidth]{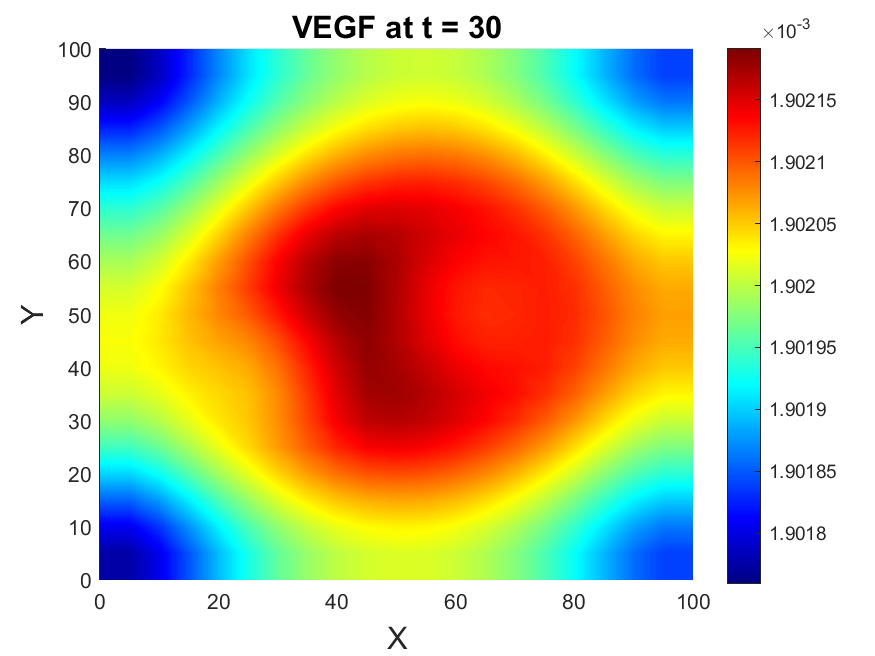}\hspace*{-0.3cm} \includegraphics[width=0.24\linewidth]{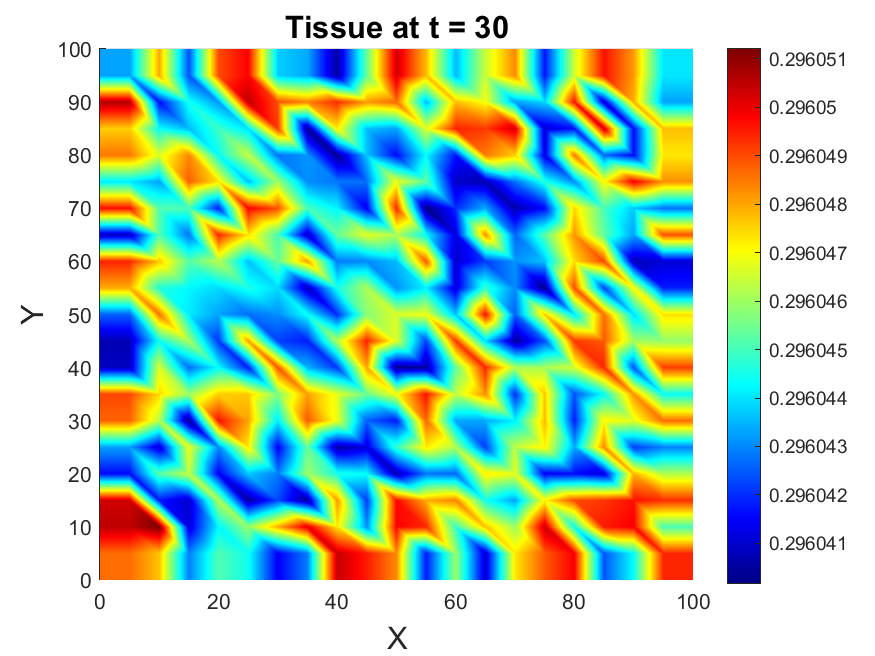}\\
	\includegraphics[width=0.24\linewidth]{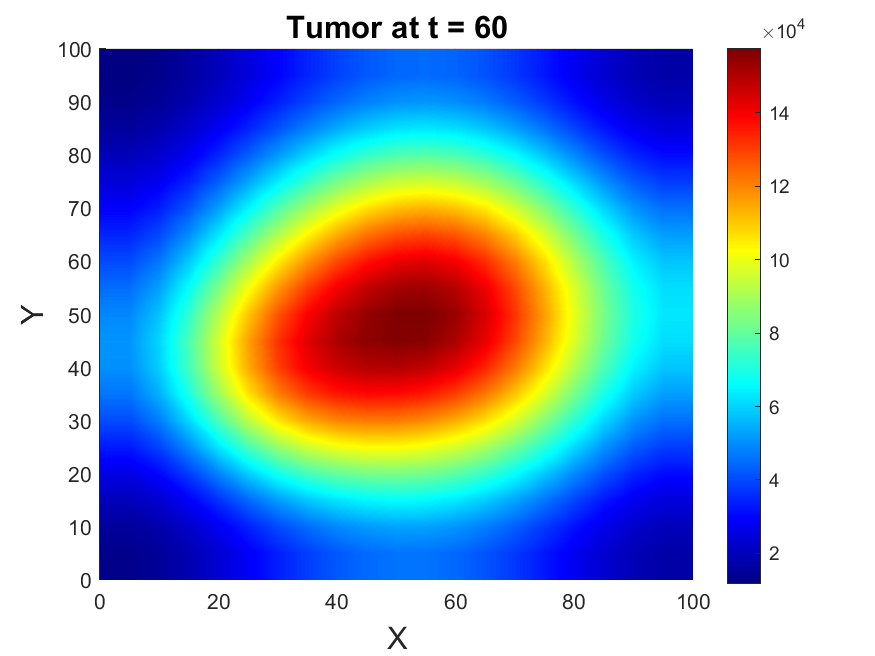}\hspace*{-0.3cm} \includegraphics[width=0.24\linewidth]{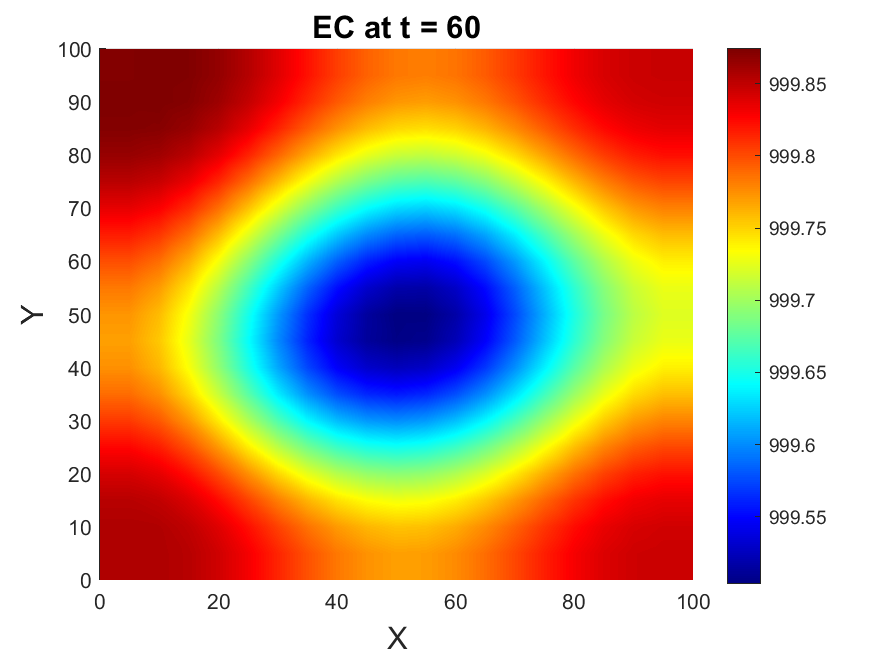}\quad \includegraphics[width=0.24\linewidth]{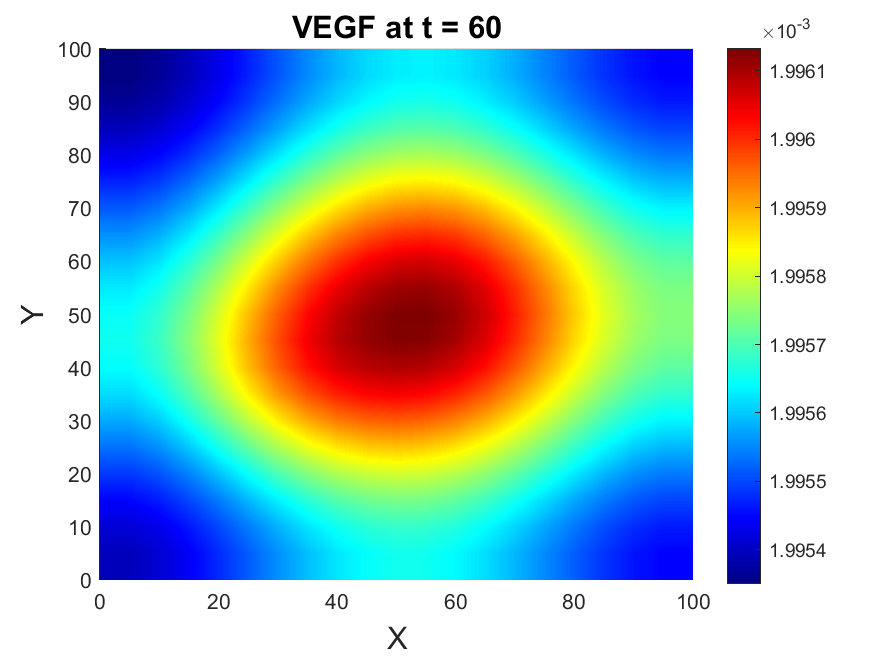}\hspace*{-0.3cm} \includegraphics[width=0.24\linewidth]{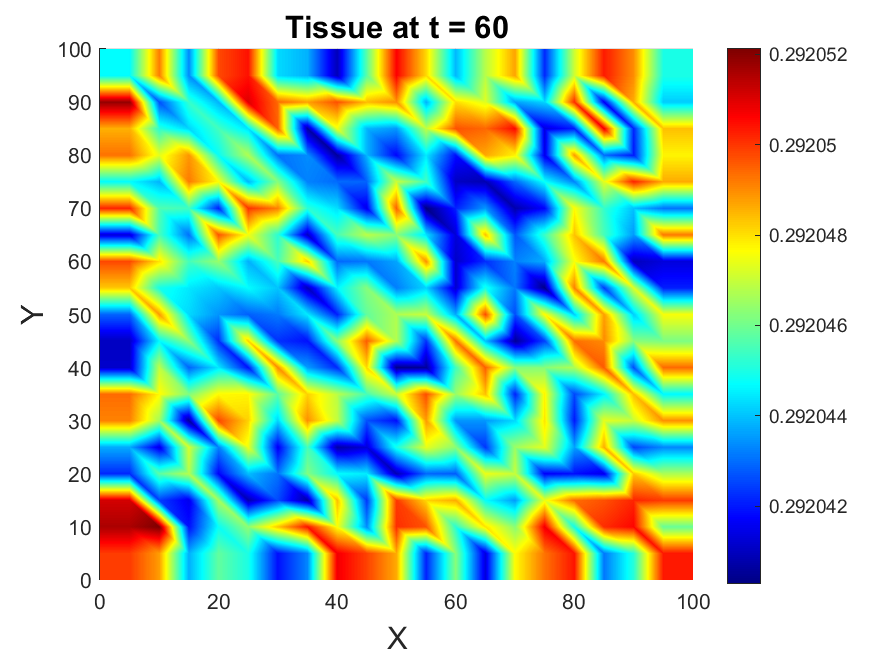}\\
	\includegraphics[width=0.24\linewidth]{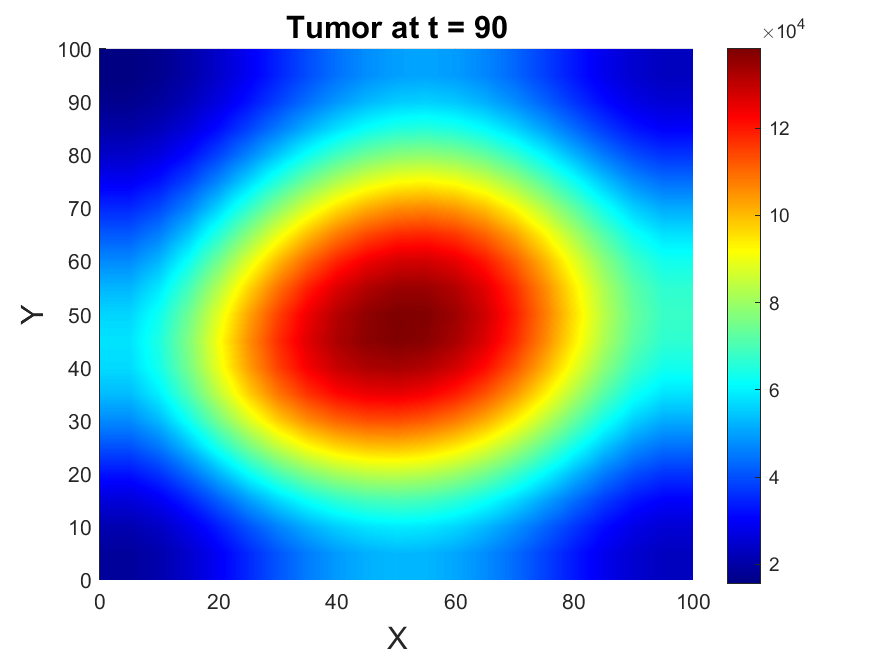}\hspace*{-0.3cm} \includegraphics[width=0.24\linewidth]{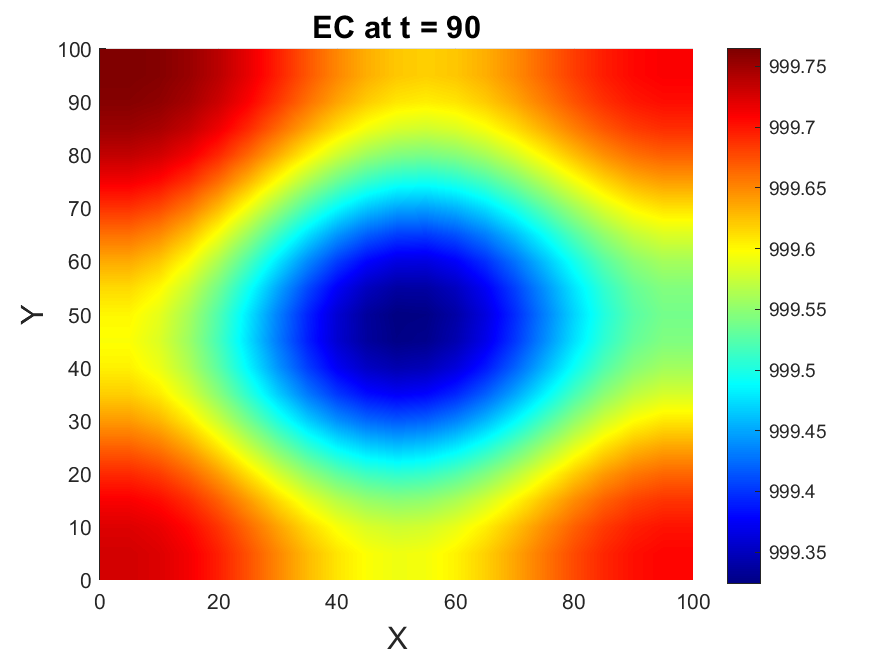}\quad \includegraphics[width=0.24\linewidth]{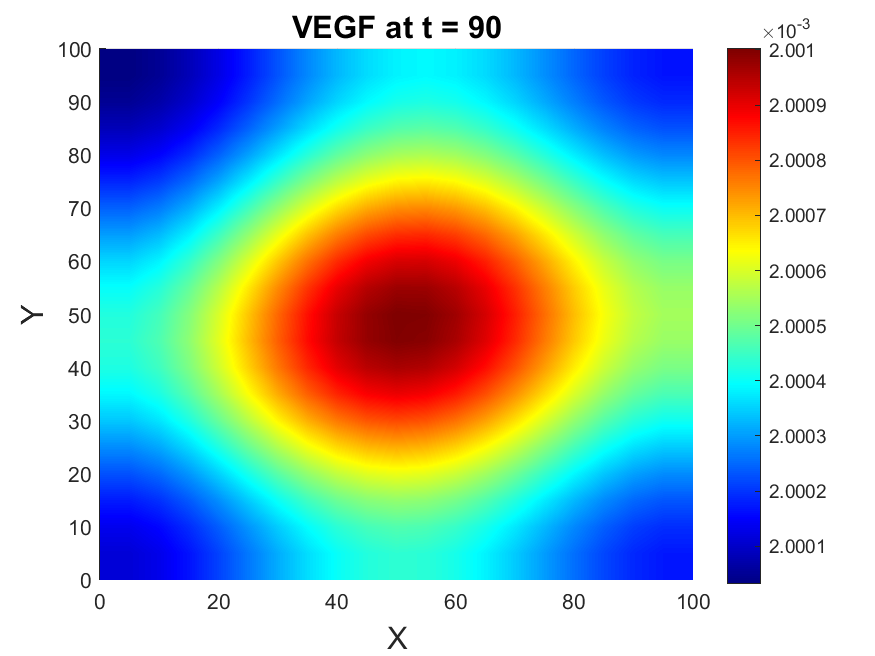}\hspace*{-0.3cm} \includegraphics[width=0.24\linewidth]{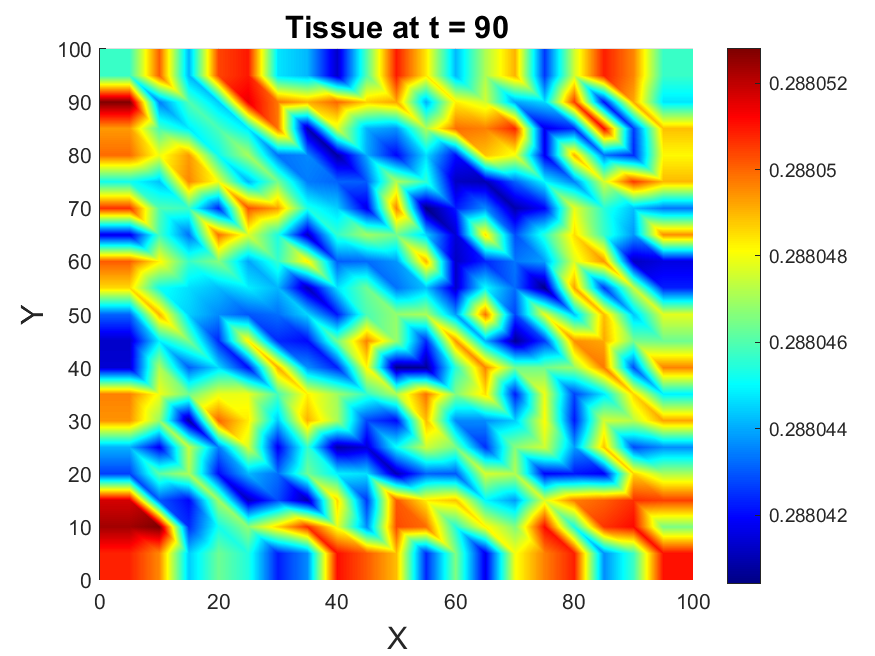}\\
	\includegraphics[width=0.24\linewidth]{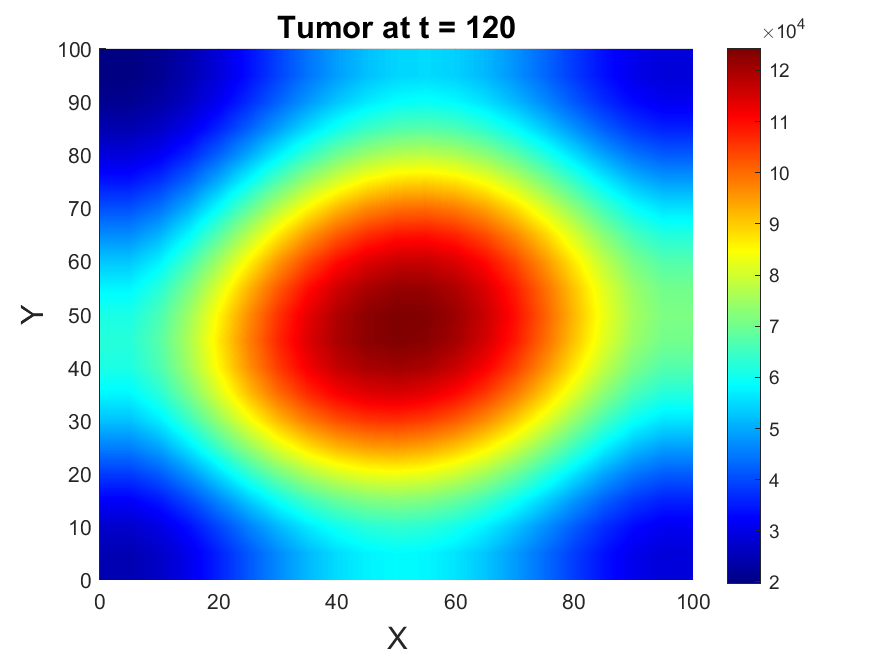}\hspace*{-0.3cm} \includegraphics[width=0.24\linewidth]{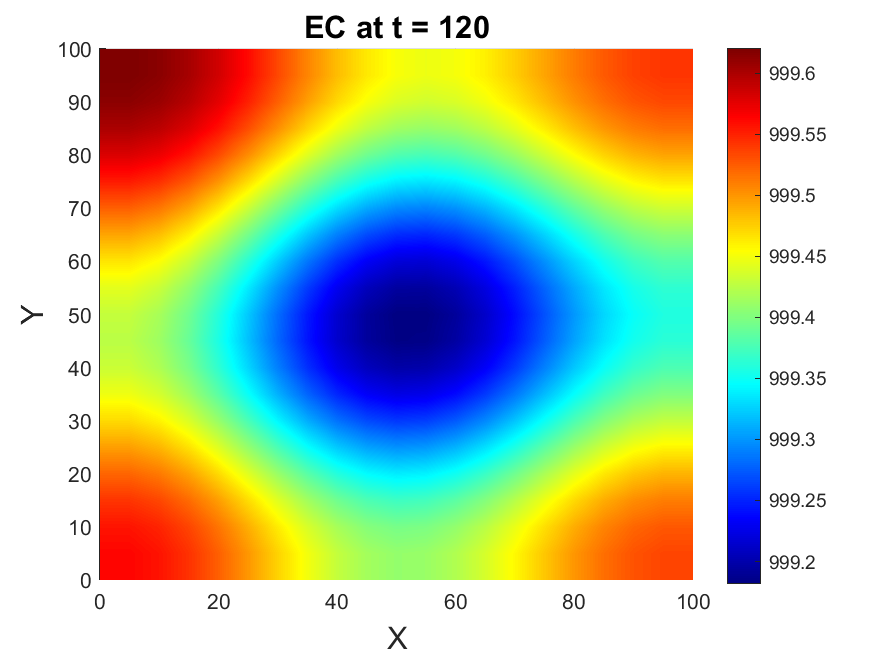}\quad \includegraphics[width=0.24\linewidth]{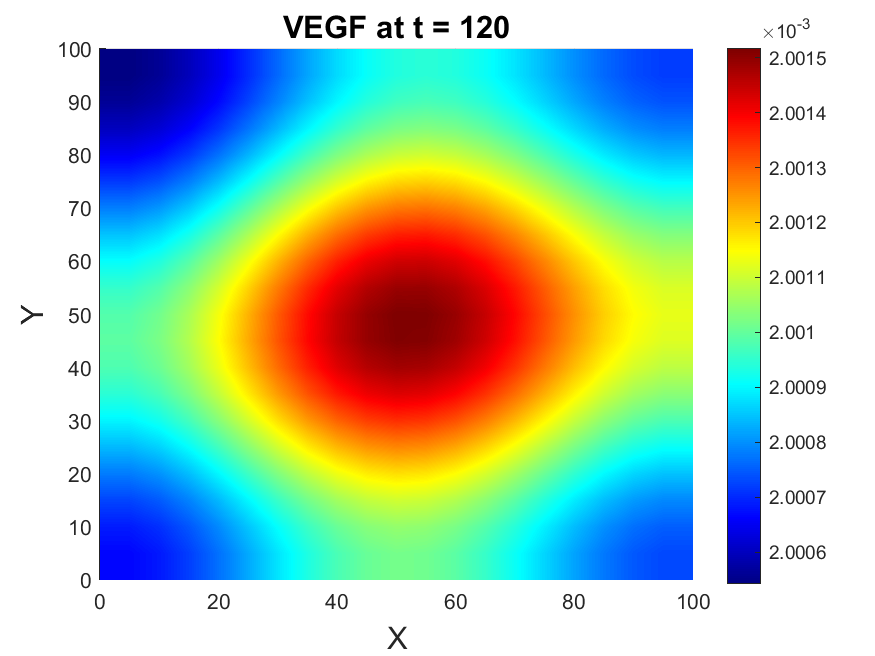}\hspace*{-0.3cm} \includegraphics[width=0.24\linewidth]{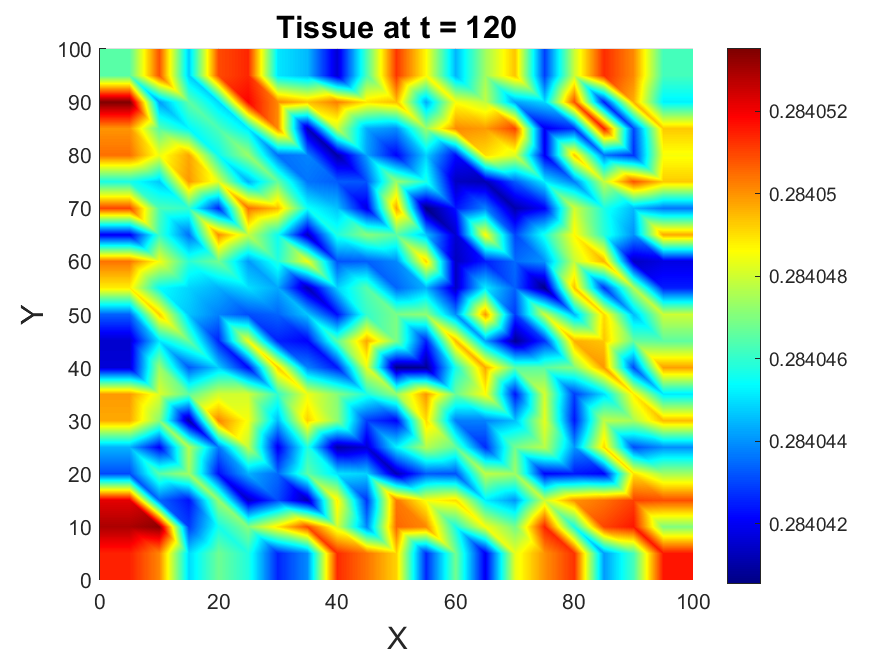}\\
	\includegraphics[width=0.24\linewidth]{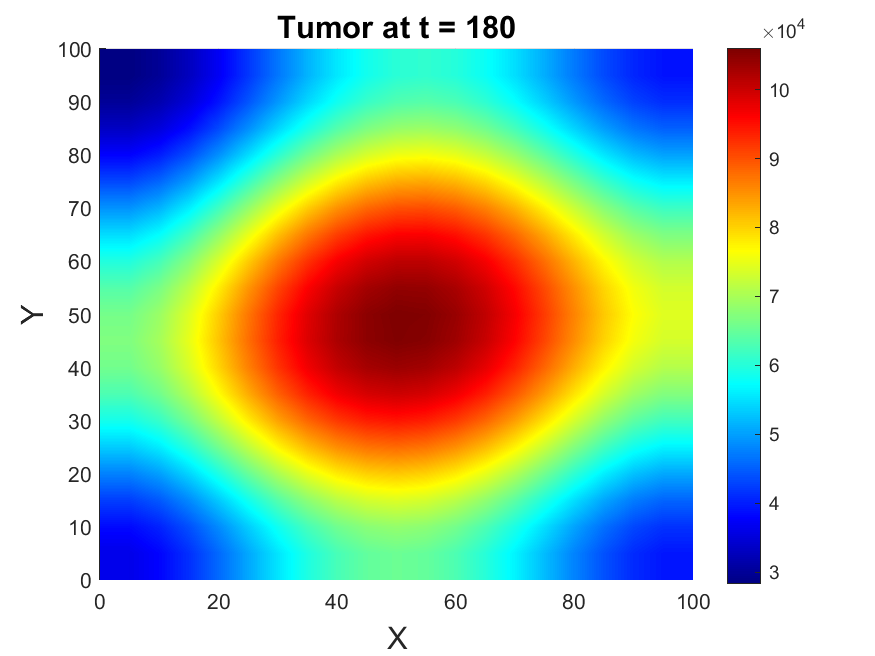}\hspace*{-0.3cm} \includegraphics[width=0.24\linewidth]{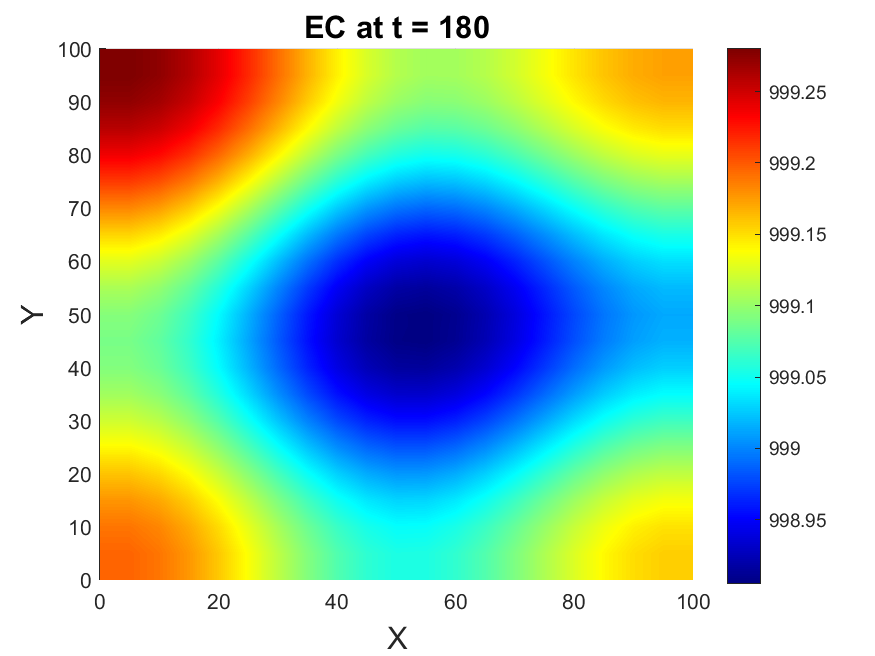}\quad \includegraphics[width=0.24\linewidth]{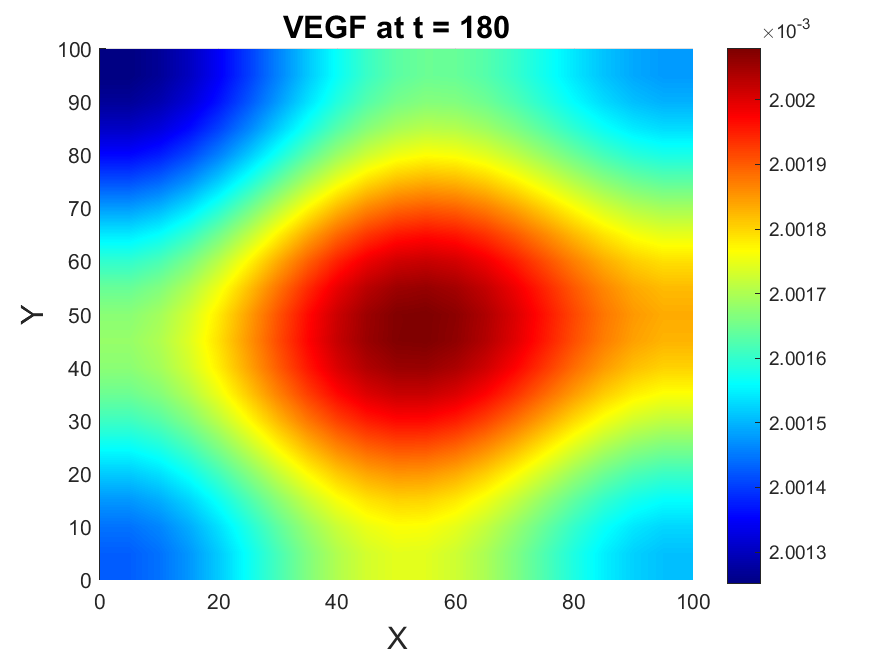}\hspace*{-0.3cm} \includegraphics[width=0.24\linewidth]{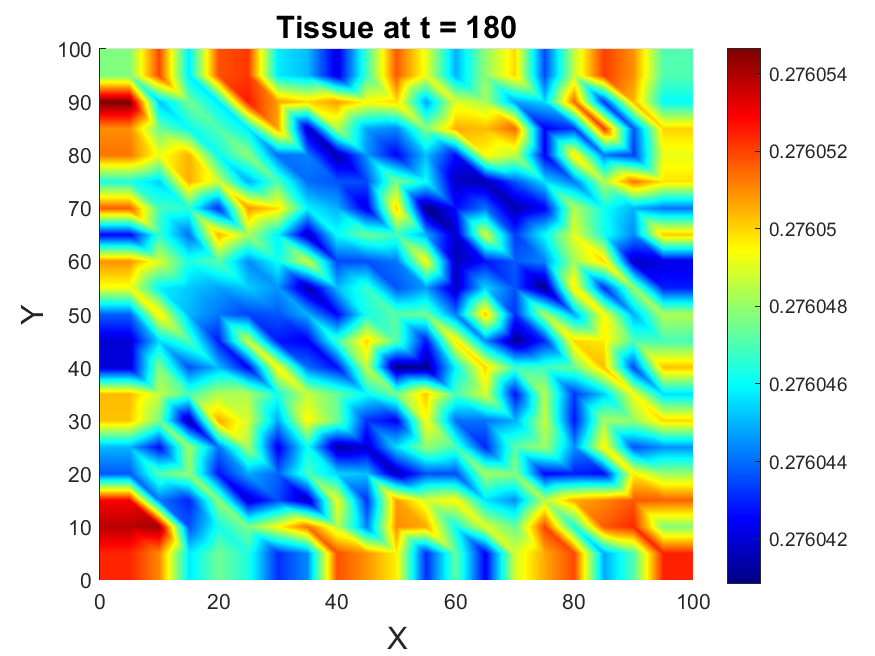}\\
	\includegraphics[width=0.24\linewidth]{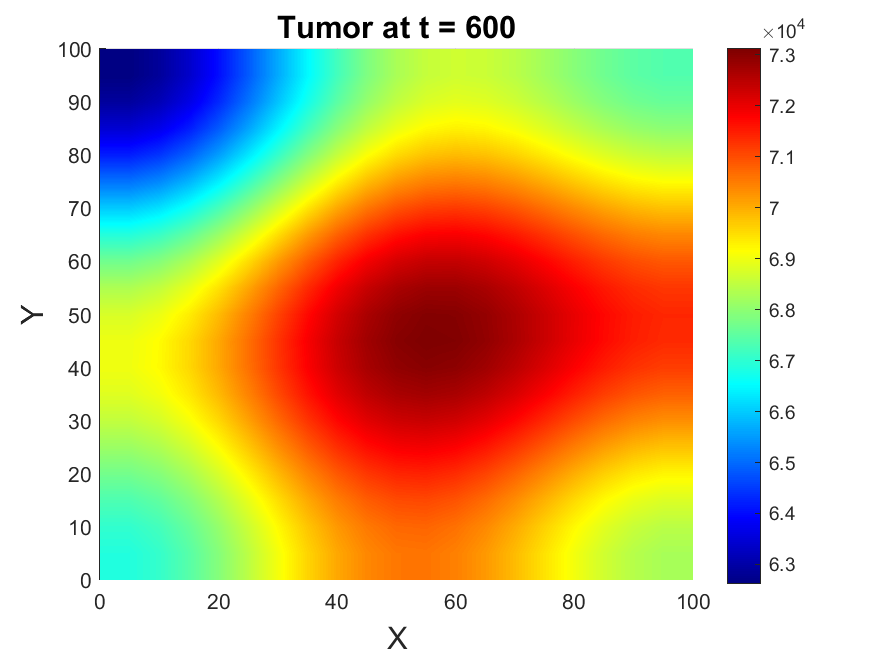}\hspace*{-0.3cm} \includegraphics[width=0.24\linewidth]{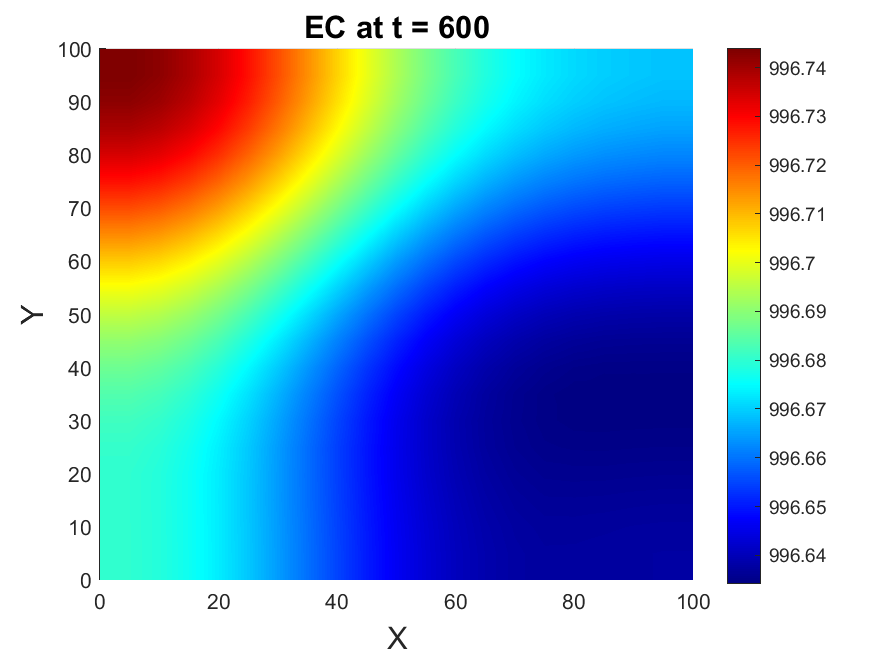}\quad \includegraphics[width=0.24\linewidth]{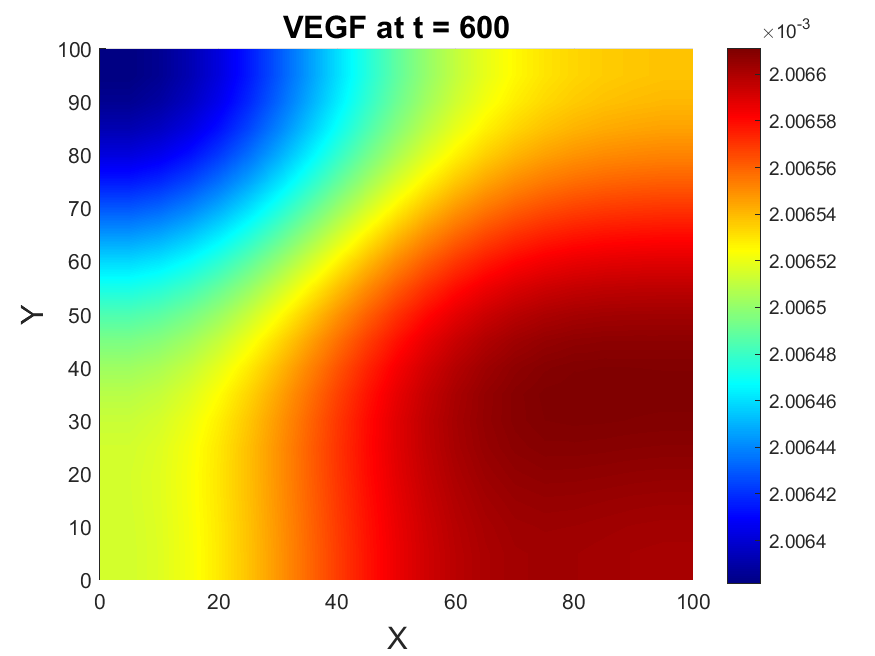}\hspace*{-0.3cm} \includegraphics[width=0.24\linewidth]{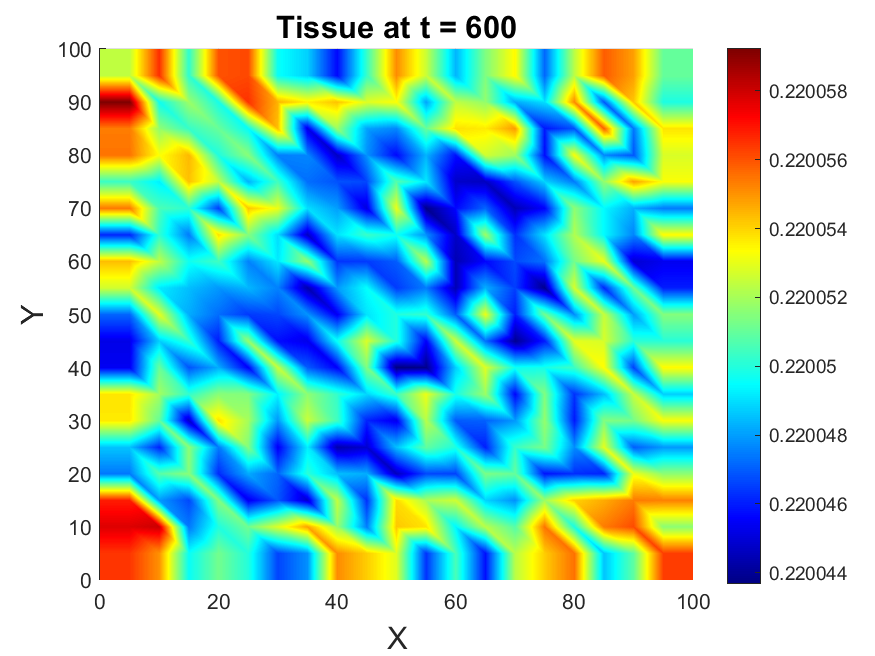}\\
	\includegraphics[width=0.24\linewidth]{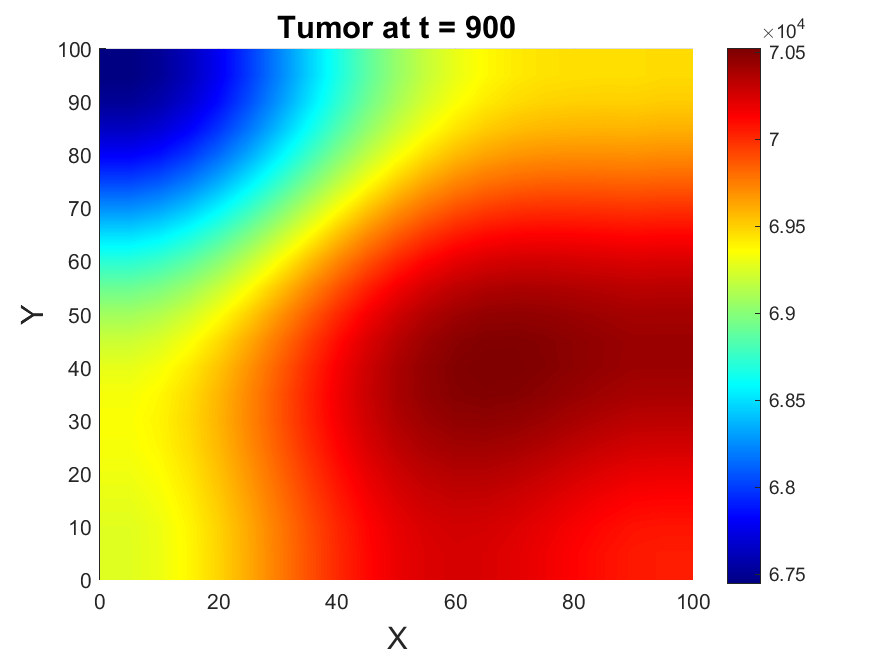}\hspace*{-0.3cm} \includegraphics[width=0.24\linewidth]{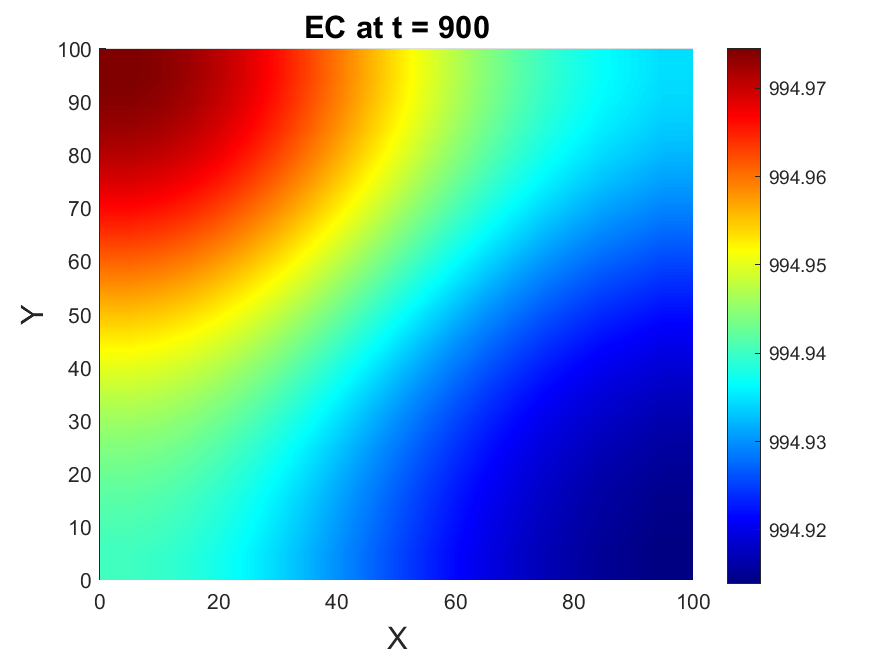}\quad \includegraphics[width=0.24\linewidth]{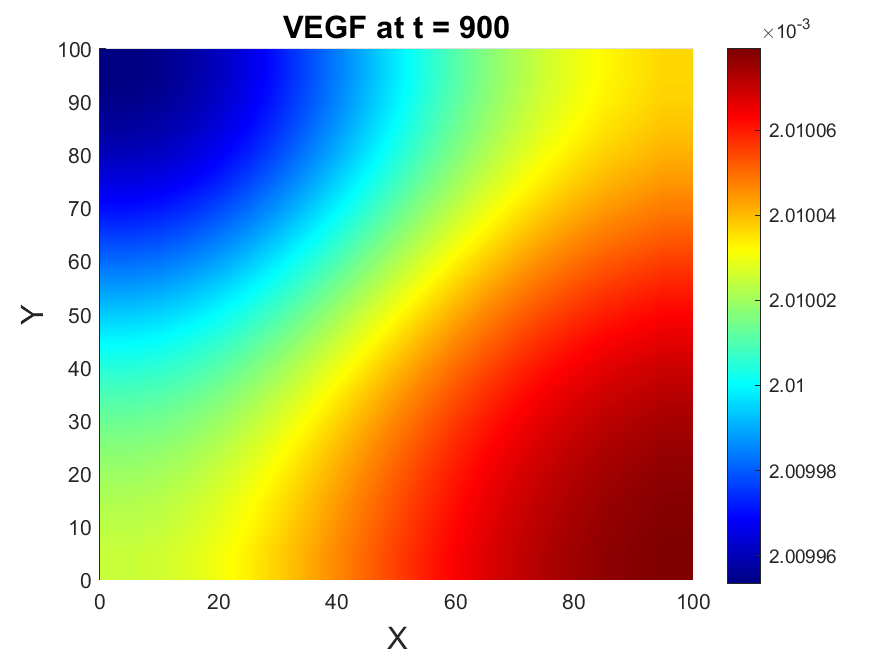}\hspace*{-0.3cm} \includegraphics[width=0.24\linewidth]{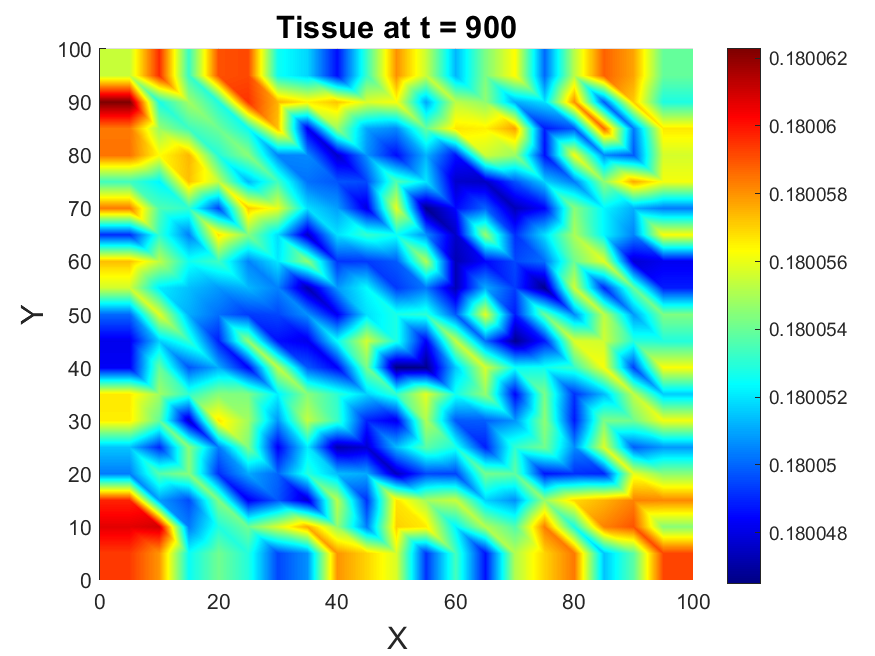}\\
	\caption{Simulations of model \eqref{q-macro-u}-\eqref{-1p2} with source terms for tumor and endothelial cells as in \eqref{eq:prolif_rates}.}\label{fig:2}
\end{figure}

\noindent
Next, we allow for proliferation of tumor cells and ECs. For the rates in \eqref{q-macro-u}, \eqref{q-macro-v} we consider 
\begin{align}
\mu_c(u,v,w):=\frac{\mu_u}{K_u}\left (1-\frac{u}{K_u}-\frac{v}{K_v}-\frac{w}{K_w}\right ),\quad \mu_e(u,v,z):=\mu_v\frac{z}{K_z}\left (1-\frac{v}{K_v}-\frac{u}{K_u}\right ),\label{eq:prolif_rates}
\end{align}
with $K_v, K_w>0$ representing the carrying capacities for ECs and tissue, respectively, $\mu_u, \mu_v>0$ are the growth rates of tumor and EC cells, respectively, and $K_z>0$ is a maximally admissible VEGF concentration. We do not account for competition of ECs with tissue, as it is negligible when compared to that with fastly growing tumor cells and intraspecific competition. The solution behavior is shown in Figure \ref{fig:2}: the tumor is evolving very similarly to the previous case; its growth is limited by all other living components of the tumor environment. The ECs start proliferating and do this more successfully in the areas where they do not infer any competition for space from tumor cells. They are also attracted by VEGF, which, however, they keep uptaking, so the highest EC densities remain in the corners of the domain, where proliferation is stronger. Eventually, the EC bulk in the upper left corner will slowly move towards the lower right corner, where VEGF and tumor cells are accumulating. This is further inhibited by the relatively large density of cancer cells present there. \\[-2ex]

\begin{figure}[h!]
	\includegraphics[width=0.27\linewidth]{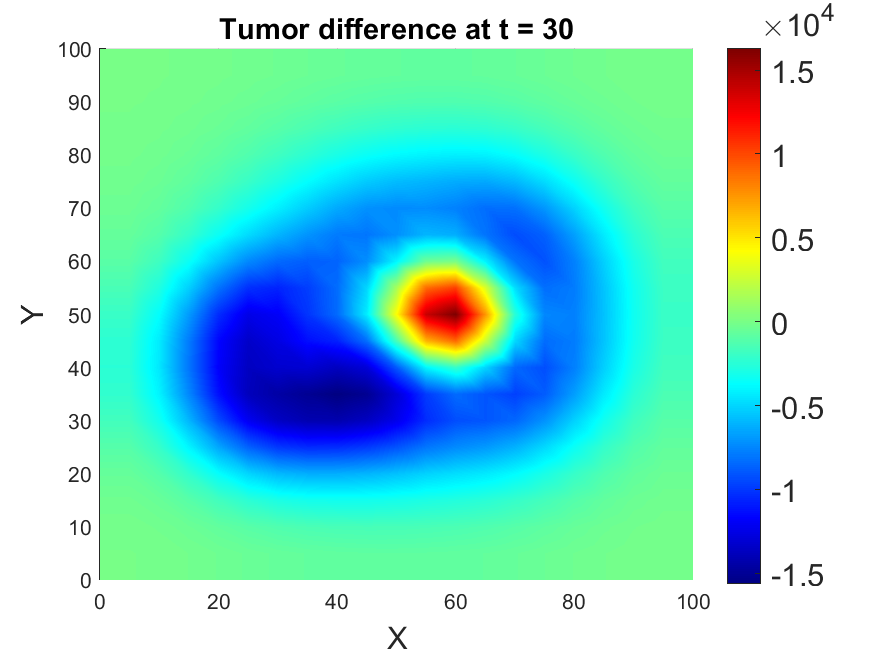}
	\includegraphics[width=0.27\linewidth]{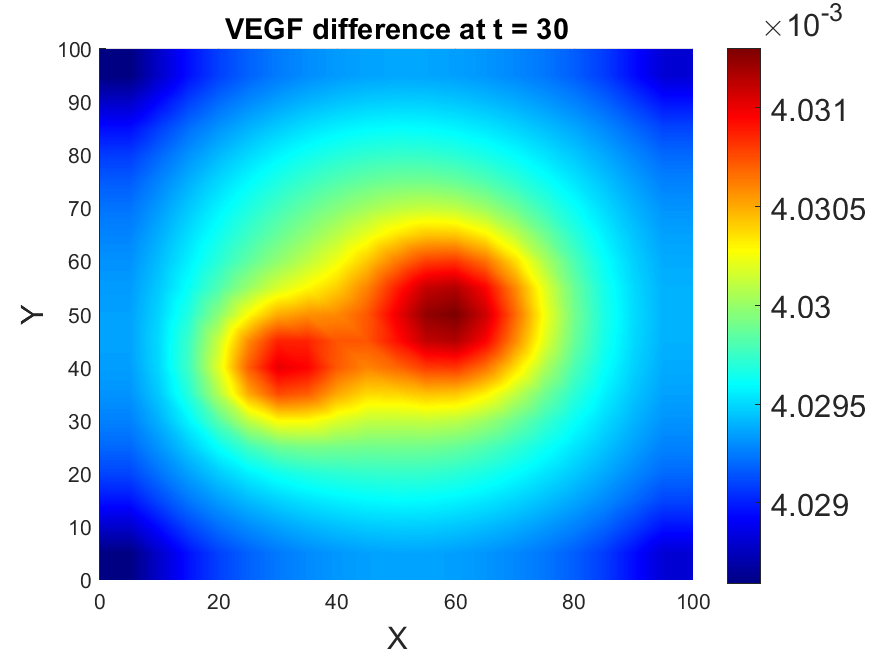}
	\includegraphics[width=0.27\linewidth]{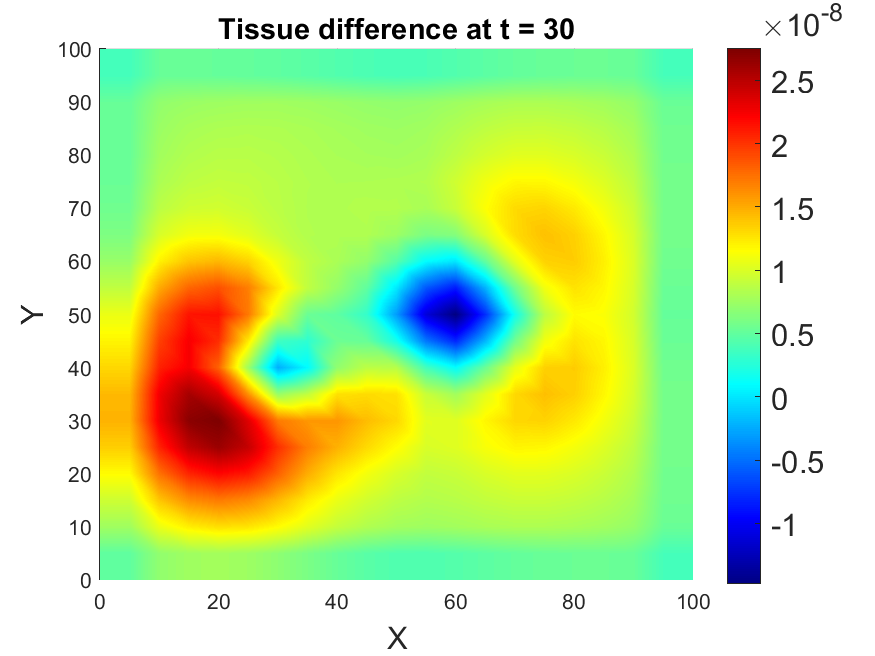}\\
	\includegraphics[width=0.27\linewidth]{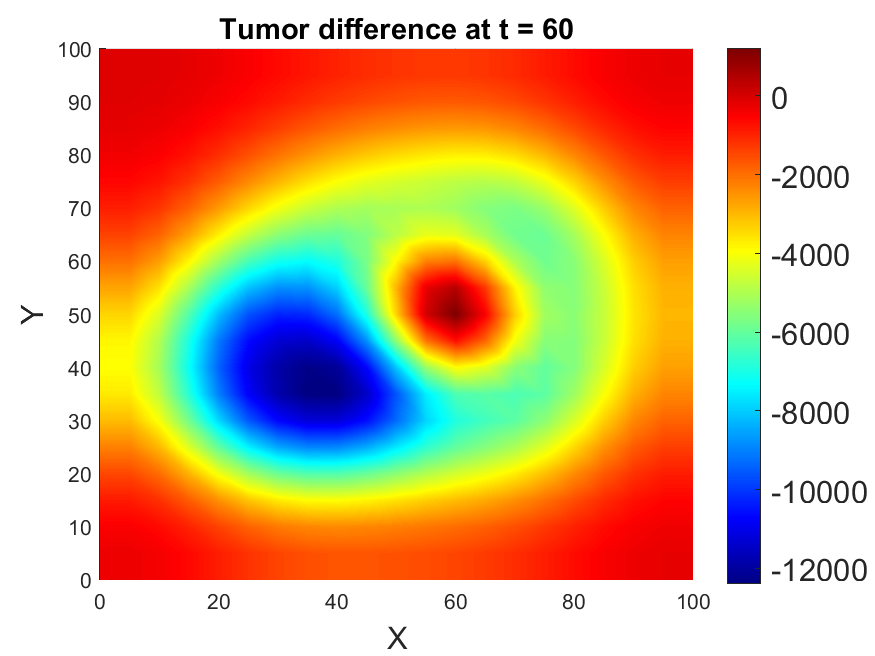}
	\includegraphics[width=0.27\linewidth]{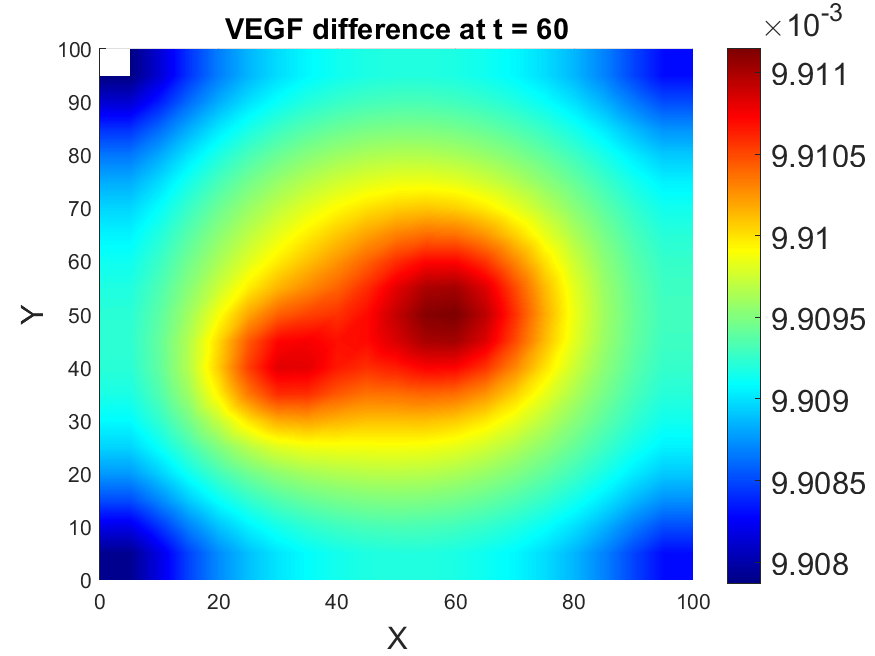}
	\includegraphics[width=0.27\linewidth]{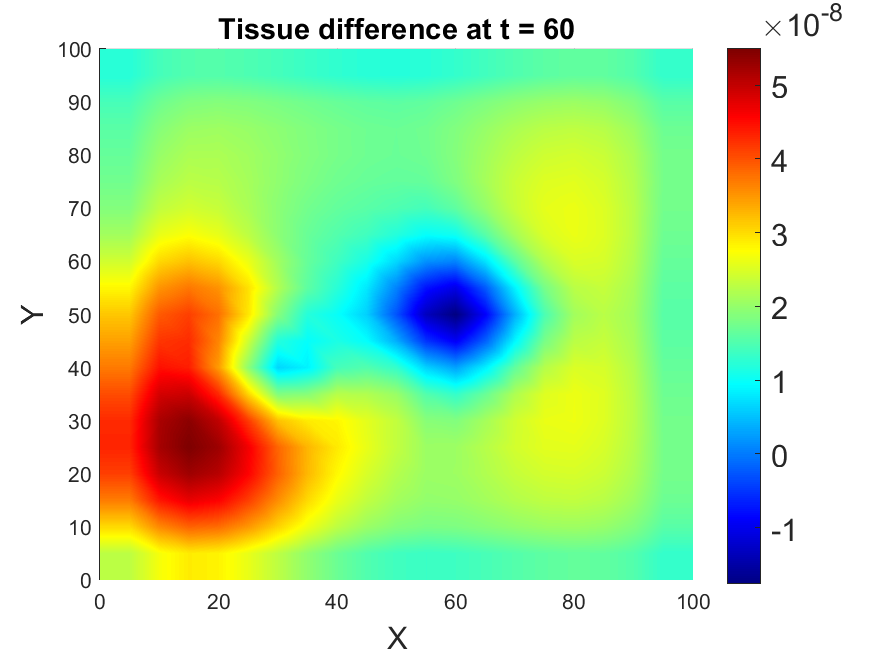}\\
	\includegraphics[width=0.27\linewidth]{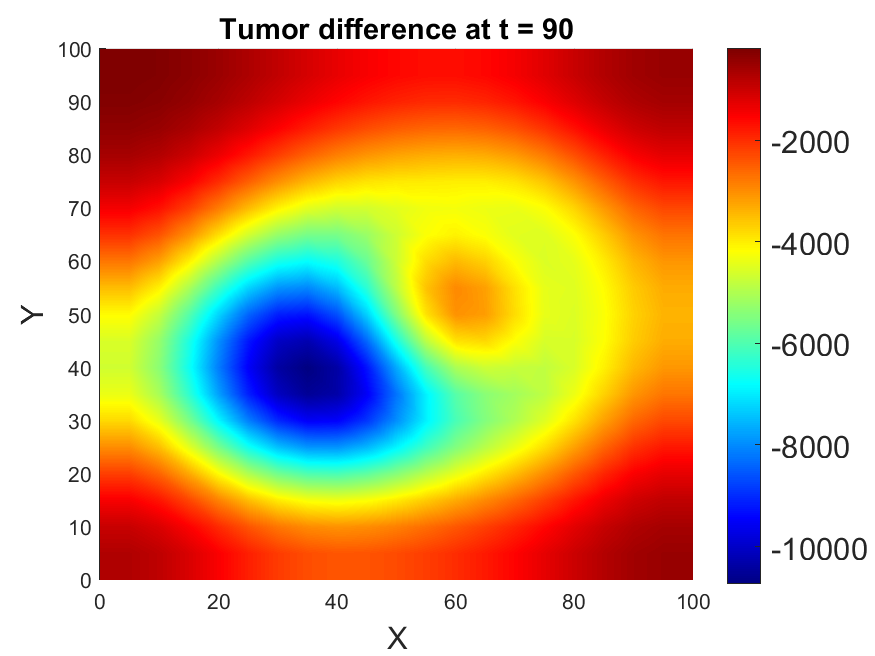}
	\includegraphics[width=0.27\linewidth]{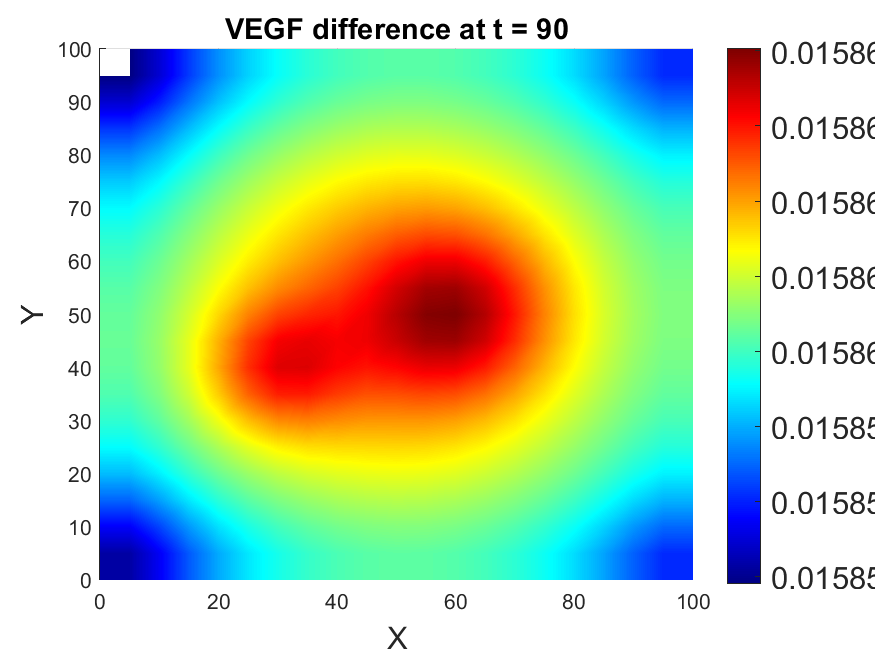}
	\includegraphics[width=0.27\linewidth]{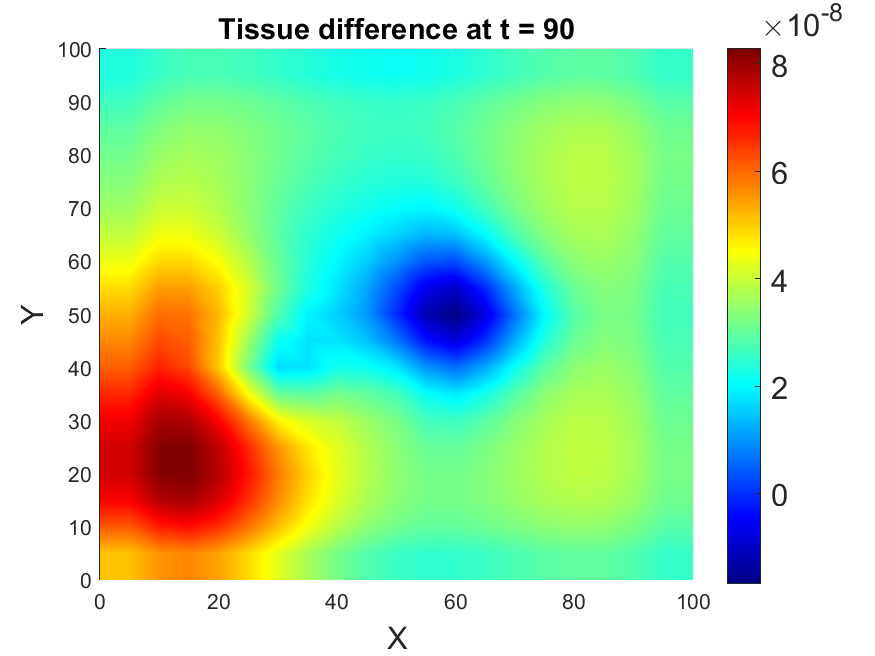}\\
	\includegraphics[width=0.27\linewidth]{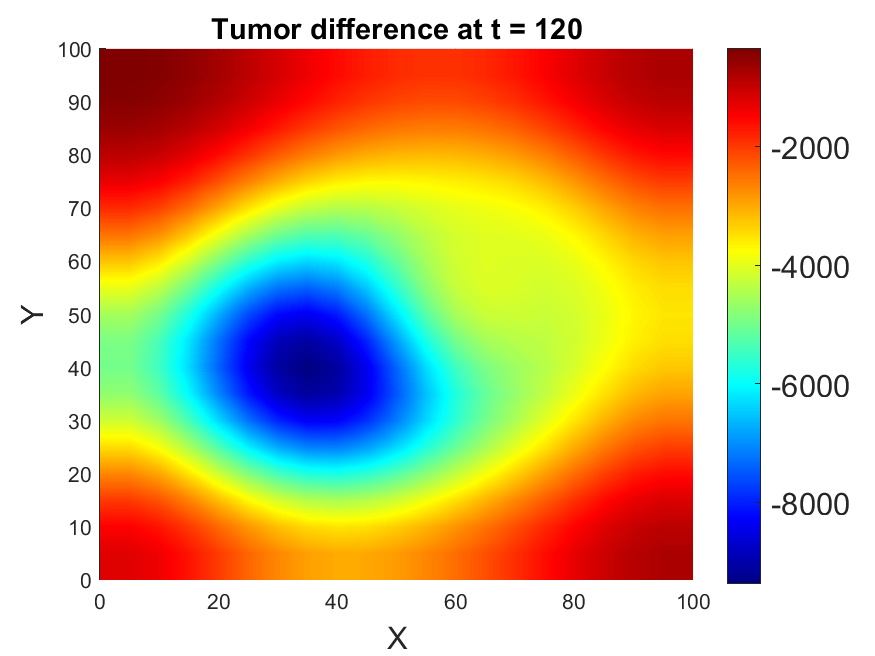}
	\includegraphics[width=0.27\linewidth]{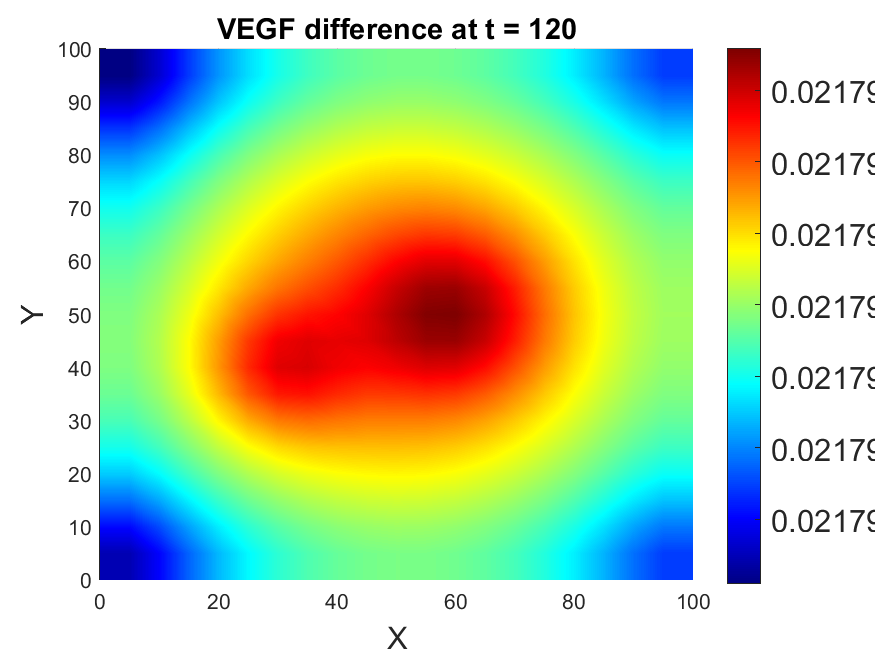}
	\includegraphics[width=0.27\linewidth]{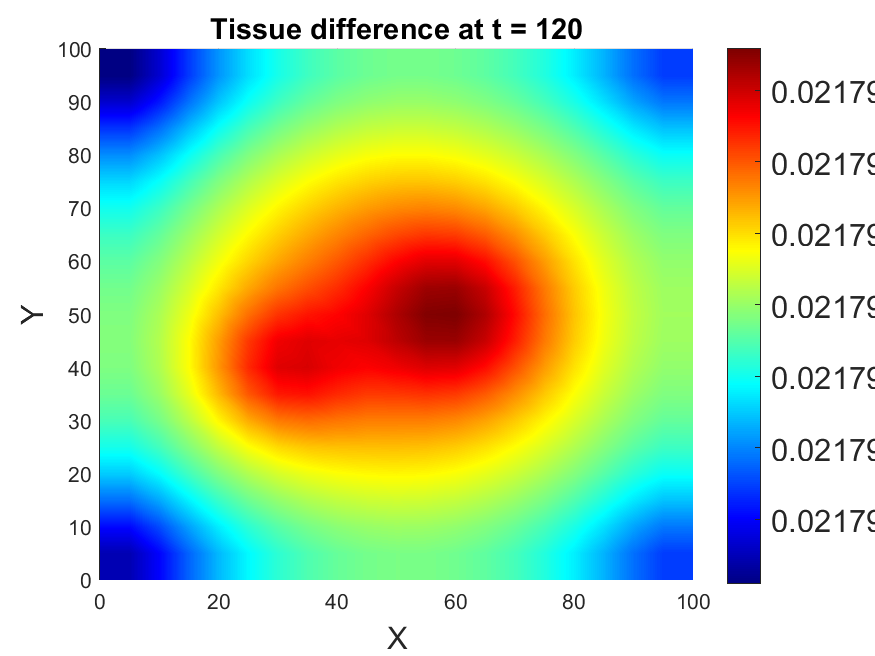}\\
	\includegraphics[width=0.27\linewidth]{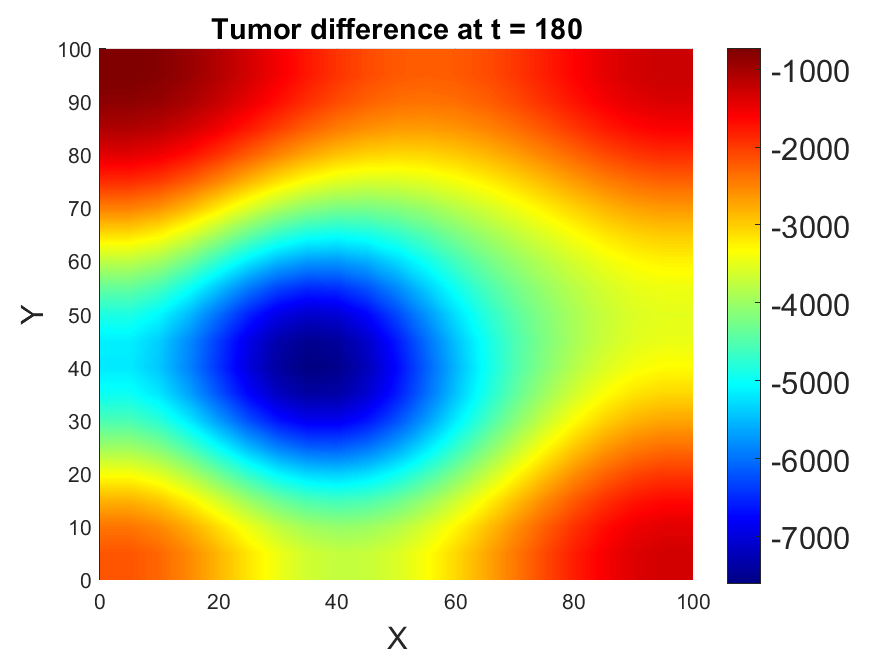}
	\includegraphics[width=0.27\linewidth]{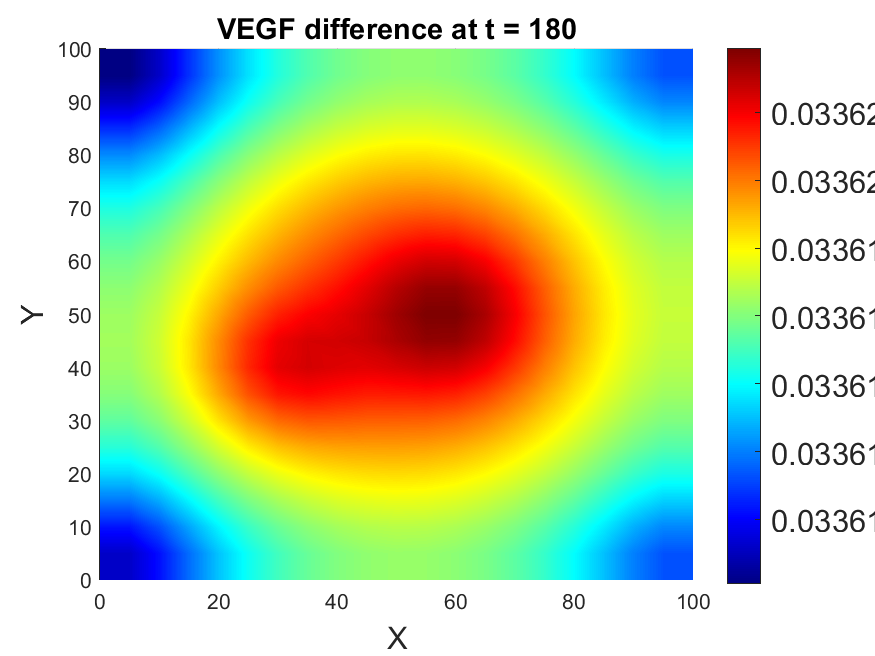}
	\includegraphics[width=0.27\linewidth]{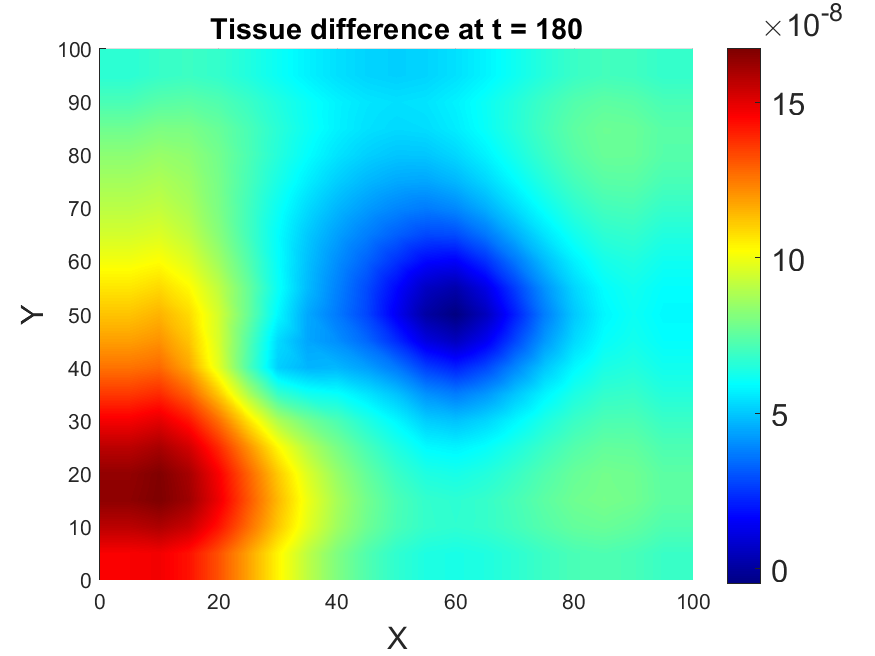}\\
	\includegraphics[width=0.27\linewidth]{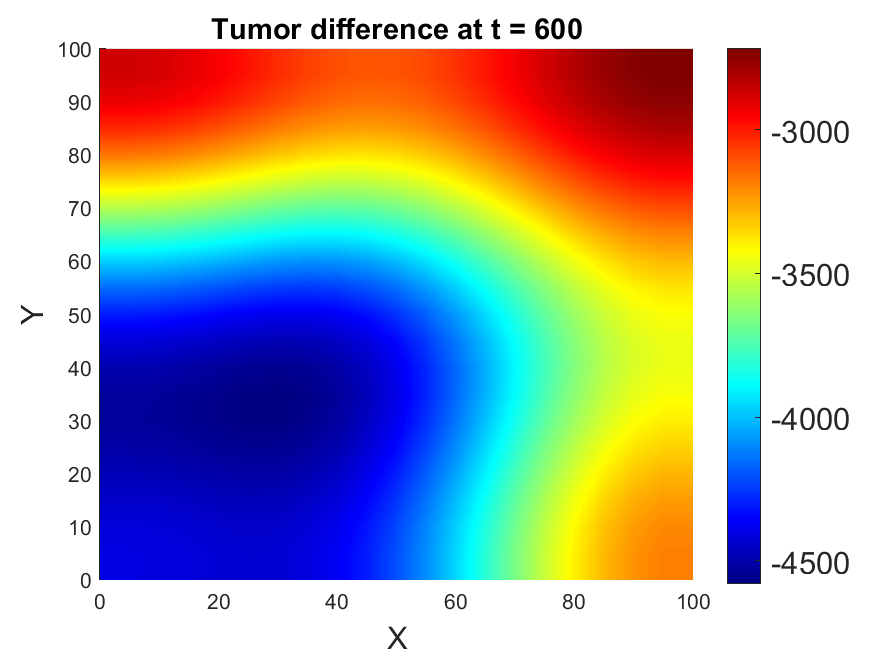}
	\includegraphics[width=0.27\linewidth]{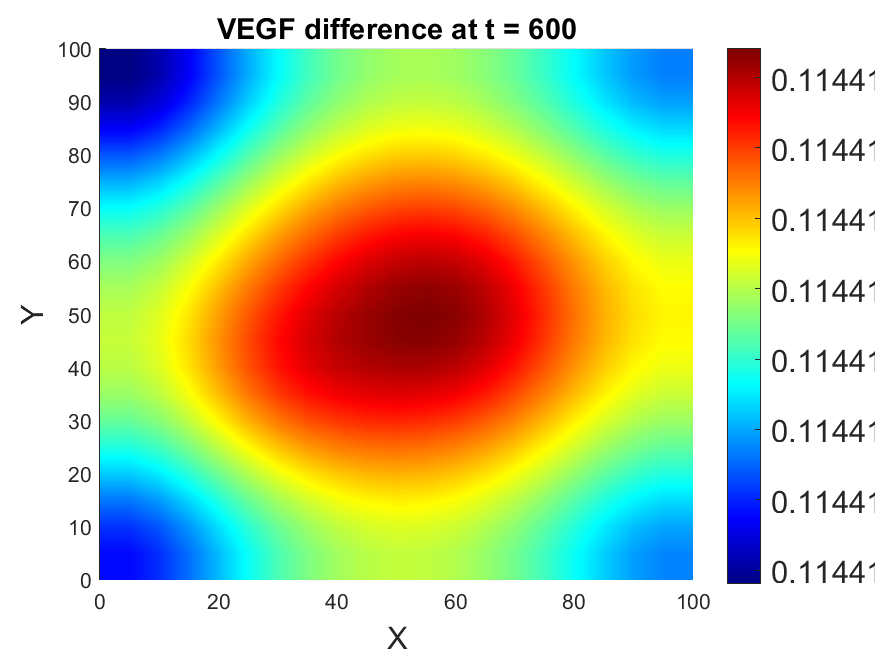}
	\includegraphics[width=0.27\linewidth]{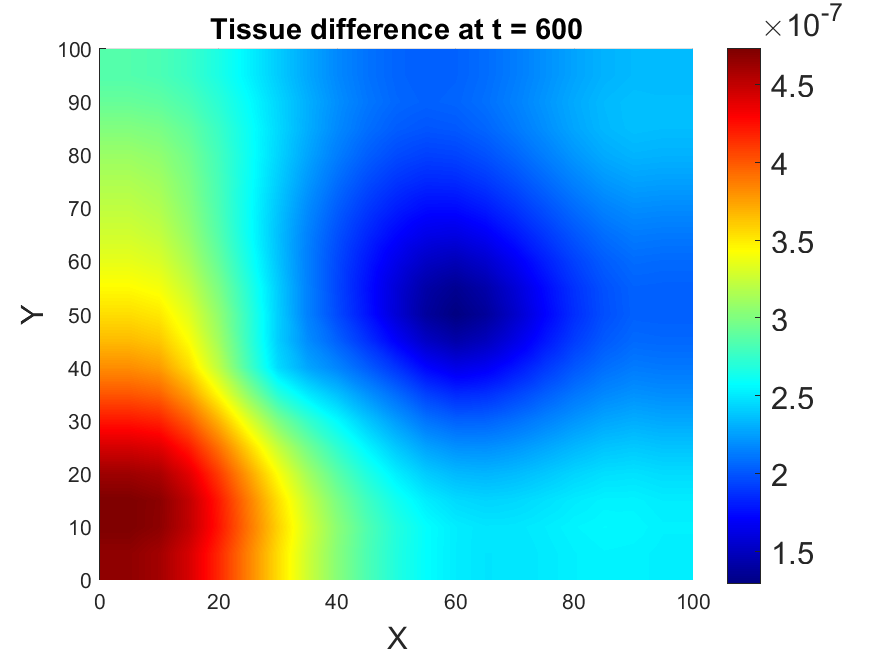}\\
	\includegraphics[width=0.27\linewidth]{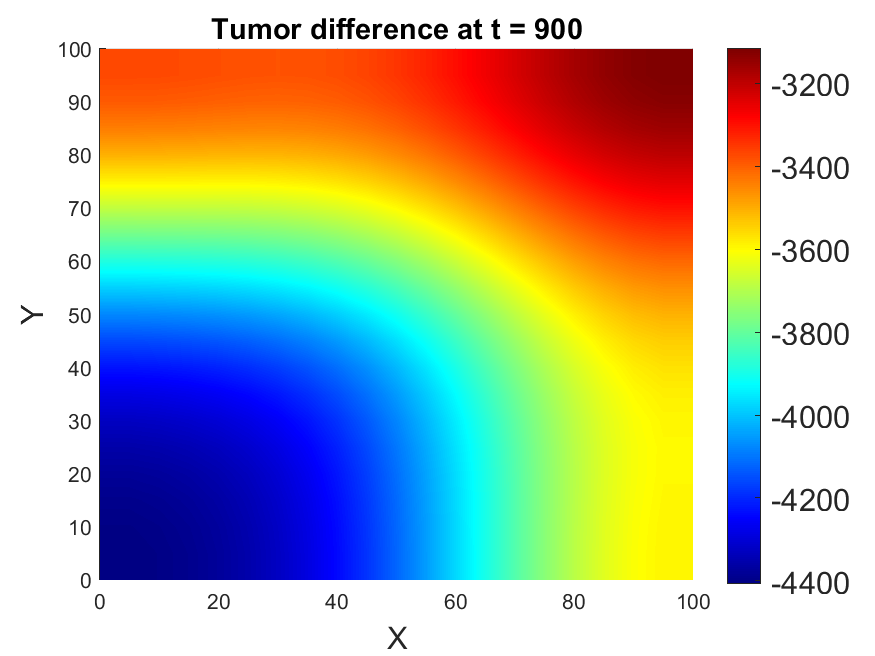}
	\includegraphics[width=0.27\linewidth]{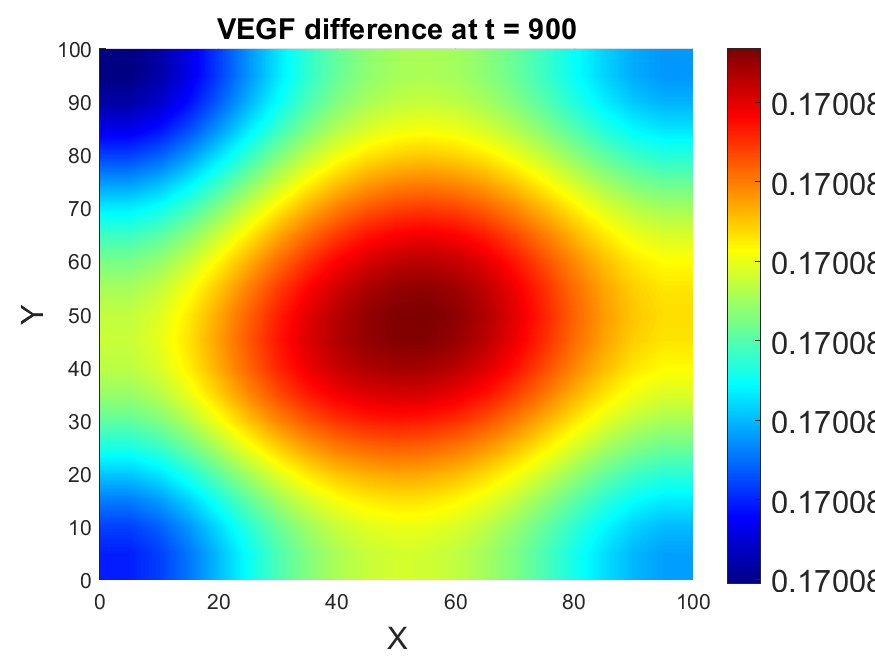}
	\includegraphics[width=0.27\linewidth]{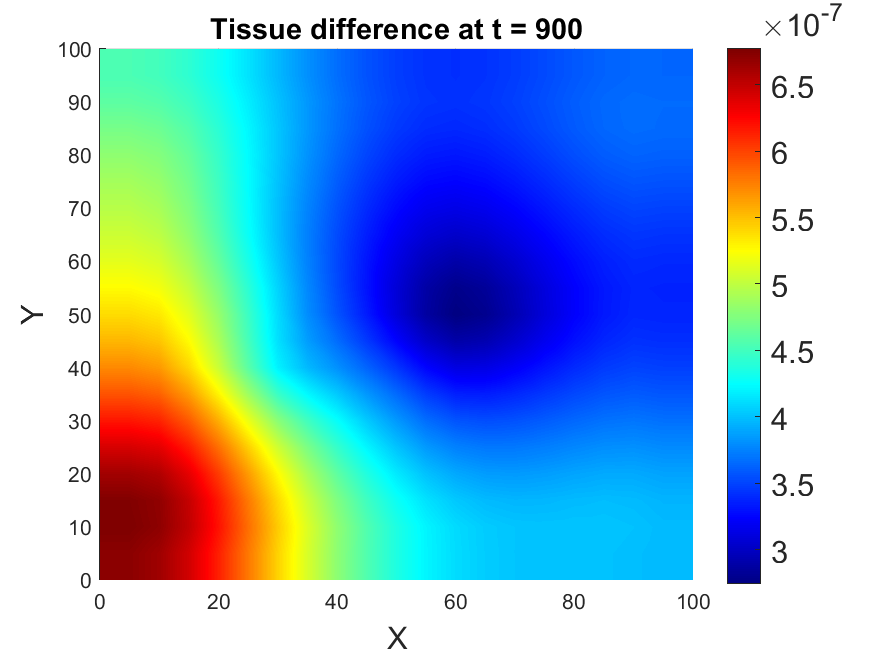}\\
	\caption{Differences between solution components of \eqref{q-macro-u}, \eqref{q-macro-v}, \eqref{-1p2} and of  \eqref{new-sys-comp}.}\label{fig:3}
\end{figure}

\noindent
Finally, we investigate the effect of the taxis cascade by comparing system \eqref{q-macro-u}, \eqref{q-macro-v}, \eqref{-1p2} with the one characterizing cancer cell invasion by performing haptotaxis up tissue gradients and chemotaxis directly following VEGF gradients. Thereby, the equation for EC density is omitted and the degradation term in the $z$-equation of \eqref{-1p2} only features natural decay. Precisely, the new system takes the form 
\be{new-sys-comp}
\left\{ \begin{array}{l}
	u_t=\Delta u - \nabla \cdot \Big( u \chi(z,w)\nabla z\Big) - \nabla \cdot \Big( u\xi(z,w)\nabla w \Big),\\
	w_t=-\psi(uw)\\
	z_t=D_z\Delta z-\mu_zz+\phi (u)
\end{array} \right.
\ee
with
\begin{align}\label{eq:chi-xi-new}
\chi(z,w)=\frac{\kappa _1}{(B(z,w))^2(1+B(z,w))}, \qquad \xi(z,w)=\frac{\kappa _2}{(B(z,w))^2(1+B(z,w))},
\end{align}
$\kappa_1, \kappa_2>0$ constants, and functions $\phi, \psi, B, \chi$ as previously. We also keep the intial conditions for $u,w,z$, and the no-flux boundary conditions for the PDEs for $u$ and $z$. To clearly put in evidence the effect of the taxis cascade we only plot the differences between solution components computed with \eqref{q-macro-u}, \eqref{q-macro-v}, \eqref{-1p2} (without source terms) and those computed with \eqref{new-sys-comp}. The results are shown in Figure \ref{fig:3}. The direct following of the chemoattractant yields a faster spread of tumor cells and accumulation at sites with higher VEGF concentration, while locations with less chemoattractant remain populated at relatively large densities and for a longer time. Overall, the model with direct taxis predicts higher cell densities all over the domain - with the correspoding accumulation tendency triggered by the tactic signal and with the largest differences at the sites whose leaving is substantially  delayed. More chemoattractant is produced by the cells at their respective locations, and more tissue is degraded by the cells available there. In absence of source terms both model versions suggest that the solution components remain bounded. The simulation outcome in Figure \ref{fig:3} also suggests that direct taxis towards VEGF can lead to a quick accumulation of tumor cells, which in presence of sufficient strong proliferation has explosion potential. Thus, the taxis cascade seems more appropriate in the analyzed context, while direct taxis has a tendency to overestimate tumor spread, accompanied by VEGF production, and slightly underestimate tissue degradation. 

\section{Discussion and outlook}\label{sec:discussion}

\noindent
The macroscopic setting \eqref{q-macro-u}, \eqref{q-macro-v}, \eqref{-1p2} deduced in this note aligns to the classification of multiple taxis models given in \cite{Kolbe2021}; more precisely, it  belongs to category (iii) therein, which has been far less addressed in literature. Such models with taxis cascade can thus make up a whole subclass of the mentioned category, but they can also be formulated in the simpler context of the population(s) performing just one taxis each. Our model addressing tumor invasion is only a paradigm; the same mathematical framework can accomodate many other applications involving different kinds of populations and their tactic signals, see e.g. \cite{Tao2019,Winkler2019,Tao2023} and references therein. The mathematical challenges coming with such models are tightly related to the concrete form of taxis and source terms. The hitherto achieved results mainly refer to various versions of forager-scrounger-nutrient taxis cascades under more or less restrictive assumptions on the data of the problem. To our knowledge the model obtained and analyzed in this note is the first one featuring a taxis cascade with chemo- and haptotaxis - with the usual mathematical issues resulting therefrom.   \\[-2ex]

\noindent
All model versions considered here required the assumptions in Subsection \ref{subsec:announce_results} to be satisfied. We also chose in Subsection \ref{subsec:model} very simple turning kernels for both cell populations. This enabled linear diffusions of tumor cells and ECs, with positive and constant diffusion coefficients. A similar meso-to-macro upscaling with more realistic turning kernels which take into account the heterogeneous, often anisotropic structure of the tissue, leads as in \cite{CEKNSSW,CKSS,Conte2023,CS20,EHKS,EHS,EKS,ZS22} to space-dependent  motility coefficients. More detailed microscopic models accounting for mechanical and/or chemical influences can lead as e.g., in \cite{DKSS20} to nonlinear (occasionally flux-limited) diffusion. Such nonconstant, possibly solution-dependent and degenerating diffusion and drift terms raise supplementary challenges in terms of well-posedness and long time behavior of solutions, see e.g., \cite{Winkler2018,Winkler2017,Heihoff2023}. 
	
\clearpage

\clearpage

\clearpage

\clearpage

\subsection*{Acknowledgement} 

\noindent
The second author acknowledges support by Deutsche
Forschungsgemeinschaft in the framework of the project  {\em
Emergence of structures and advantages in cross-diffusion systems}
(No.~411007140, GZ: WI 3707/5-1).

\phantomsection
\printbibliography

@article{CS20,
	title = {Mathematical modeling of glioma invasion: acid- and vasculature mediated go-or-grow dichotomy and the influence of tissue anisotropy},
	volume = {407},
	ISSN = {0096-3003},
	url = {http://dx.doi.org/10.1016/j.amc.2021.126305},
	DOI = {10.1016/j.amc.2021.126305},
	journal = {Applied Mathematics and Computation},
	publisher = {Elsevier BV},
	author = {Conte,  Martina and Surulescu,  Christina},
	year = {2021},
	pages = {126305}
}

@article{Conte2023,
	title = {Mathematical modeling of glioma invasion and therapy approaches via kinetic theory of active particles},
	volume = {33},
	ISSN = {1793-6314},
	url = {http://dx.doi.org/10.1142/S0218202523500227},
	DOI = {10.1142/s0218202523500227},
	number = {05},
	journal = {Mathematical Models and Methods in Applied Sciences},
	publisher = {World Scientific Pub Co Pte Ltd},
	author = {Conte,  Martina and Dzierma,  Yvonne and Knobe,  Sven and Surulescu,  Christina},
	year = {2023},
	pages = {1009-1051}
}

@article{alikakos,
	doi = {10.1080/03605307908820113},
	url = {https://doi.org/10.1080/03605307908820113},
	year = {1979},
	publisher = {Informa {UK} Limited},
	volume = {4},
	number = {8},
	pages = {827--868},
	author = {N.D. Alikakos},
	title = {{LPBounds} of solutions of reaction-diffusion equations},
	journal = {Communications in Partial Differential Equations}
}

@ARTICLE{amann,
	author = {Herbert Amann},
	title = {Compact embeddings of vector-valued Sobolev and Besov spaces},
	journal = {Glasnik Mat.},
	year = {2000},
	pages = {161--177},
	volume={35},
}

@article{BBTW,
	doi = {10.1142/s021820251550044x},
	url = {https://doi.org/10.1142/s021820251550044x},
	year = {2015},
	month = may,
	publisher = {World Scientific Pub Co Pte Lt},
	volume = {25},
	number = {09},
	pages = {1663--1763},
	author = {N. Bellomo and A. Bellouquid and Y. Tao and M. Winkler},
	title = {Toward a mathematical theory of Keller{\textendash}Segel models of pattern formation in biological tissues},
	journal = {Mathematical Models and Methods in Applied Sciences}
}

@article{BHN,
	doi = {10.1016/0362-546x(94)90101-5},
	url = {https://doi.org/10.1016/0362-546x(94)90101-5},
	year = {1994},
	publisher = {Elsevier {BV}},
	volume = {23},
	number = {9},
	pages = {1189-1209},
	author = {Piotr Biler and Waldemar Hebisch and Tadeusz Nadzieja},
	title = {The Debye system: existence and large time behavior of solutions},
	journal = {Nonlinear Analysis: Theory,  Methods {\&} Applications}
}

@article{fontelos_friedman_hu,
	doi = {10.1137/s0036141001385046},
	url = {https://doi.org/10.1137/s0036141001385046},
	year = {2002},
	publisher = {Society for Industrial {\&} Applied Mathematics ({SIAM})},
	volume = {33},
	number = {6},
	pages = {1330-1355},
	author = {Marco A. Fontelos and Avner Friedman and Bei Hu},
	title = {Mathematical Analysis of a Model for the Initiation of Angiogenesis},
	journal = {{SIAM} Journal on Mathematical Analysis}
}

@article{friedman_tello,
	doi = {10.1016/s0022-247x(02)00147-6},
	url = {https://doi.org/10.1016/s0022-247x(02)00147-6},
	year = {2002},
	publisher = {Elsevier {BV}},
	volume = {272},
	number = {1},
	pages = {138-163},
	author = {Avner Friedman and J.Ignacio Tello},
	title = {Stability of solutions of chemotaxis equations in reinforced random walks},
	journal = {Journal of Mathematical Analysis and Applications}
}

@article{FIWY,
	doi = {10.3934/dcds.2016.36.151},
	url = {https://doi.org/10.3934/dcds.2016.36.151},
	year = {2015},
	publisher = {American Institute of Mathematical Sciences ({AIMS})},
	volume = {36},
	number = {1},
	pages = {151-169},
	author = {Tomomi Yokota and Michael Winkler and Akio Ito and Kentarou Fujie},
	title = {Stabilization in a chemotaxis model for tumor invasion},
	journal = {Discrete and Continuous Dynamical Systems}
}

@article{giga_sohr,
	doi = {10.1016/0022-1236(91)90136-s},
	url = {https://doi.org/10.1016/0022-1236(91)90136-s},
	year = {1991},
	publisher = {Elsevier {BV}},
	volume = {102},
	number = {1},
	pages = {72-94},
	author = {Yoshikazu Giga and Hermann Sohr},
	title = {Abstract Lp estimates for the Cauchy problem with applications to the Navier-Stokes equations in exterior domains},
	journal = {Journal of Functional Analysis}
}

@book {GT,
	AUTHOR = {Gilbarg, D. and Trudinger, N.S.},
	TITLE = {Elliptic Partial Differential Equations of Second Order},
	PUBLISHER = {Springer-Verlag, Berlin/Heidelberg},
	YEAR = {2001},
}

@misc{LSU,
	Author = {O.A. {Ladyzhenskaya} and V.A. {Solonnikov} and N.N. {Ural'tseva}},
	Title = {{Linear and quasi-linear equations of parabolic type. Translated from the Russian by S. Smith.}},
	Year = {1968},
	HowPublished = {{Translations of Mathematical Monographs. 23. Providence, RI: American Mathematical Society (AMS). XI, 648 p.}},
	Zbl = {0174.15403},
}

@article{lions_ARMA,
	doi = {10.1007/bf00249679},
	url = {https://doi.org/10.1007/bf00249679},
	year = {1980},
	month = dec,
	publisher = {Springer Science and Business Media {LLC}},
	volume = {74},
	number = {4},
	pages = {335--353},
	author = {P. L. Lions},
	title = {R{\'{e}}solution de probl{\`{e}}mes elliptiques quasilin{\'{e}}aires},
	journal = {Archive for Rational Mechanics and Analysis}
}

@article{pang_wang_JDE,
	doi = {10.1016/j.jde.2017.03.016},
	url = {https://doi.org/10.1016/j.jde.2017.03.016},
	year = {2017},
	publisher = {Elsevier {BV}},
	volume = {263},
	number = {2},
	pages = {1269-1292},
	author = {Peter Y.H. Pang and Yifu Wang},
	title = {Global existence of a two-dimensional chemotaxis-haptotaxis model with remodeling of non-diffusible attractant},
	journal = {Journal of Differential Equations}
}

@article{taowin_subcritical,
	doi = {10.1016/j.jde.2011.08.019},
	url = {https://doi.org/10.1016/j.jde.2011.08.019},
	year = {2012},
	publisher = {Elsevier {BV}},
	volume = {252},
	number = {1},
	pages = {692-715},
	author = {Youshan Tao and Michael Winkler},
	title = {Boundedness in a quasilinear parabolic-parabolic Keller-Segel system with subcritical sensitivity},
	journal = {Journal of Differential Equations}
}

@article{taowin_consumption,
	doi = {10.1016/j.jde.2011.07.010},
	url = {https://doi.org/10.1016/j.jde.2011.07.010},
	year = {2012},
	publisher = {Elsevier {BV}},
	volume = {252},
	number = {3},
	pages = {2520-2543},
	author = {Youshan Tao and Michael Winkler},
	title = {Eventual smoothness and stabilization of large-data solutions in a three-dimensional chemotaxis system with consumption of chemoattractant},
	journal = {Journal of Differential Equations}
}

@article{taowin_JDE2014,
	doi = {10.1016/j.jde.2014.04.014},
	url = {https://doi.org/10.1016/j.jde.2014.04.014},
	year = {2014},
	publisher = {Elsevier {BV}},
	volume = {257},
	number = {3},
	pages = {784-815},
	author = {Youshan Tao and Michael Winkler},
	title = {Energy-type estimates and global solvability in a two-dimensional chemotaxis-haptotaxis model with remodeling of non-diffusible attractant},
	journal = {Journal of Differential Equations}
}

@article{win_JDE,
	doi = {10.1016/j.jde.2010.02.008},
	url = {https://doi.org/10.1016/j.jde.2010.02.008},
	year = {2010},
	publisher = {Elsevier {BV}},
	volume = {248},
	number = {12},
	pages = {2889-2905},
	author = {Michael Winkler},
	title = {Aggregation vs. global diffusive behavior in the higher-dimensional Keller-Segel model},
	journal = {Journal of Differential Equations}
}

@article{win_CPDE2012,
	doi = {10.1080/03605302.2011.591865},
	url = {https://doi.org/10.1080/03605302.2011.591865},
	year = {2012},
	publisher = {Informa {UK} Limited},
	volume = {37},
	number = {2},
	pages = {319-351},
	author = {Michael Winkler},
	title = {Global Large-Data Solutions in a Chemotaxis-(Navier-)Stokes System Modeling Cellular Swimming in Fluid Drops},
	journal = {Communications in Partial Differential Equations}
}

@article{win_NON_exp,
	doi = {10.1088/1361-6544/aa565b},
	url = {https://doi.org/10.1088/1361-6544/aa565b},
	year = {2017},
	publisher = {{IOP} Publishing},
	volume = {30},
	number = {2},
	pages = {735-764},
	author = {Michael Winkler},
	title = {Global existence and slow grow-up in a quasilinear Keller-Segel system with exponentially decaying diffusivity},
	journal = {Nonlinearity}
}

@book{bellom3,
	author = {Bellomo, Nicola and Bellouquid, Abdelghani and Gibelli, L. and Outada, N.},
	title = {A Quest Towards a Mathematical Theory of Living Systems},
	publisher={Birkh\"auser},
	year={2018},
}

@article{EHKS,
	author = {Engwer, Christian and Hillen, Thomas and Knappitsch, Markus and Surulescu, Christina},
	journal = {Journal of mathematical biology},
	number = {3},
	pages = {551-582},
	title = {Glioma follow white matter tracts: a multiscale DTI-based model},
	volume = {71},
	year = {2015},
}

@article{EHS,
	author = {Engwer, Christian and Hunt, Alexander and Surulescu, Christina},
	journal = {Mathematical medicine and biology: a journal of the IMA},
	number = {4},
	pages = {435-459},
	title = {Effective equations for anisotropic glioma spread with proliferation: a multiscale approach},
	volume = {33},
	year = {2016},
}

@article{EKS,
	author = {Engwer, Christian and Knappitsch, Markus and Surulescu, Christina},
	journal = {Mathematical Biosciences \& Engineering},
	number = {2},
	pages = {443-460},
	title = {A multiscale model for glioma spread including cell-tissue interactions and proliferation},
	volume = {13},
	year = {2016},
}

@article{OtHi,
	author = {Othmer, Hans G and Hillen, Thomas},
	journal = {SIAM Journal on Applied Mathematics},
	number = {4},
	pages = {1220-1250},
	title = {The diffusion limit of transport equations II: Chemotaxis equations},
	volume = {62},
	year = {2003},
}

@article{PeTaWa,
	author = {Perthame, Beno{\^\i}t and Tang, Min and Vauchelet, Nicolas},
	journal = {Journal of mathematical biology},
	number = {5},
	pages = {1161-1178},
	title = {Derivation of the bacterial run-and-tumble kinetic equation from a model with biochemical pathway},
	volume = {73},
	year = {2016},
}

@article{CKSS,
	author = {Corbin, Gregor and Hunt, Alexander and Klar, A and Schneider, F and Surulescu, Christina},
	journal = {Mathematical Models and Methods in Applied Sciences},
	pages = {1771-1800},
	title = {Higher-order models for glioma invasion: from a two-scale description to effective equations for mass density and momentum},
	volume = {28},
	year = {2018},
}

@article{CEKNSSW,
	author={G. Corbin and C. Engwer and A. Klar and J. Nieto and J. Soler and C. Surulescu and M. Wenske},
	title={Modeling glioma invasion with anisotropy- and hypoxia-triggered motility enhancement: from subcellular dynamics to macroscopic PDEs with multiple taxis},
	 volume = {31},
	ISSN = {1793-6314},
	url = {http://dx.doi.org/10.1142/S0218202521500056},
	DOI = {10.1142/s0218202521500056},
	number = {01},
	journal = {Mathematical Models and Methods in Applied Sciences},
	publisher = {World Scientific Pub Co Pte Ltd},
	year = {2020},
	pages = {177-222}
}

@article{DKSS20,
	title = {Multiscale Modeling of Glioma Invasion: From Receptor Binding to Flux-Limited Macroscopic PDEs},
	volume = {20},
	ISSN = {1540-3467},
	url = {http://dx.doi.org/10.1137/21M1412104},
	DOI = {10.1137/21m1412104},
	number = {2},
	journal = {Multiscale Modeling \& Simulation},
	publisher = {Society for Industrial \& Applied Mathematics (SIAM)},
	author = {Dietrich,  Anne and Kolbe,  Niklas and Sfakianakis,  Nikolaos and Surulescu,  Christina},
	year = {2022},
	pages = {685-713}
}

@article{ZS22,
	title = {A Novel Derivation of Rigorous Macroscopic Limits from a Micro-Meso Description of Signal-Triggered Cell Migration in Fibrous Environments},
	volume = {82},
	ISSN = {1095-712X},
	url = {http://dx.doi.org/10.1137/20M1365442},
	DOI = {10.1137/20m1365442},
	number = {1},
	journal = {SIAM Journal on Applied Mathematics},
	publisher = {Society for Industrial \& Applied Mathematics (SIAM)},
	author = {Zhigun,  Anna and Surulescu,  Christina},
	year = {2022},
	pages = {142-167}
}

@article{Lorenz2014,
	doi = {10.1142/s0218202514500249},
	url = {https://doi.org/10.1142/s0218202514500249},
	year = {2014},
	publisher = {World Scientific Pub Co Pte Lt},
	volume = {24},
	number = {12},
	pages = {2383-2436},
	author = {T. Lorenz and C. Surulescu},
	title = {On a class of multiscale cancer cell migration models: Well-posedness in less regular function spaces},
	journal = {Mathematical Models and Methods in Applied Sciences}
}

@article{Kumar2021,
	title = {Multiscale modeling of glioma pseudopalisades: contributions from the tumor microenvironment},
	volume = {82},
	ISSN = {1432-1416},
	url = {http://dx.doi.org/10.1007/s00285-021-01599-x},
	DOI = {10.1007/s00285-021-01599-x},
	number = {6},
	journal = {Journal of Mathematical Biology},
	publisher = {Springer Science and Business Media LLC},
	author = {Kumar,  Pawan and Li,  Jing and Surulescu,  Christina},
	year = {2021},
}

@article{Ascher1995,
	title = {Implicit-Explicit Methods for Time-Dependent Partial Differential Equations},
	volume = {32},
	ISSN = {1095-7170},
	url = {http://dx.doi.org/10.1137/0732037},
	DOI = {10.1137/0732037},
	number = {3},
	journal = {SIAM Journal on Numerical Analysis},
	publisher = {Society for Industrial \& Applied Mathematics (SIAM)},
	author = {Ascher,  Uri M. and Ruuth,  Steven J. and Wetton,  Brian T. R.},
	year = {1995},
	pages = {797-823}
}

@book{weickert1998anisotropic,
	title={Anisotropic diffusion in image processing},
	author={Weickert, J.},
	volume={1},
	year={1998},
	publisher={Teubner Stuttgart}
}

@article{Hanahan2011,
	title = {Hallmarks of Cancer: The Next Generation},
	volume = {144},
	ISSN = {0092-8674},
	url = {http://dx.doi.org/10.1016/j.cell.2011.02.013},
	DOI = {10.1016/j.cell.2011.02.013},
	number = {5},
	journal = {Cell},
	publisher = {Elsevier BV},
	author = {Hanahan,  Douglas and Weinberg,  Robert A.},
	year = {2011},
	pages = {646-674}
}

@article{Abakarova1995,
	title = {The metastatic potential of tumors depends on the pH of host tissues},
	volume = {120},
	ISSN = {1573-8221},
	url = {http://dx.doi.org/10.1007/BF02445579},
	DOI = {10.1007/bf02445579},
	number = {6},
	journal = {Bulletin of Experimental Biology and Medicine},
	publisher = {Springer Science and Business Media LLC},
	author = {Abakarova,  O. R.},
	year = {1995},
	pages = {1227-1229}
}

@article{Gatenby2004,
	title = {Why do cancers have high aerobic glycolysis?},
	volume = {4},
	ISSN = {1474-1768},
	url = {http://dx.doi.org/10.1038/nrc1478},
	DOI = {10.1038/nrc1478},
	number = {11},
	journal = {Nature Reviews Cancer},
	publisher = {Springer Science and Business Media LLC},
	author = {Gatenby,  Robert A. and Gillies,  Robert J.},
	year = {2004},
	pages = {891-899}
}

@article{Hardee2012,
	title = {Mechanisms of Glioma-Associated Neovascularization},
	volume = {181},
	ISSN = {0002-9440},
	url = {http://dx.doi.org/10.1016/j.ajpath.2012.06.030},
	DOI = {10.1016/j.ajpath.2012.06.030},
	number = {4},
	journal = {The American Journal of Pathology},
	publisher = {Elsevier BV},
	author = {Hardee,  Matthew E. and Zagzag,  David},
	year = {2012},
	pages = {1126-1141}
}

@article{Wang2019,
	title = {Association between Tumor Acidity and Hypervascularity in Human Gliomas Using pH-Weighted Amine Chemical Exchange Saturation Transfer Echo-Planar Imaging and Dynamic Susceptibility Contrast Perfusion MRI at 3T},
	volume = {40},
	ISSN = {1936-959X},
	url = {http://dx.doi.org/10.3174/ajnr.a6063},
	DOI = {10.3174/ajnr.a6063},
	number = {6},
	journal = {American Journal of Neuroradiology},
	publisher = {American Society of Neuroradiology (ASNR)},
	author = {Wang,  Y.-L. and Yao,  J. and Chakhoyan,  A. and Raymond,  C. and Salamon,  N. and Liau,  L.M. and Nghiemphu,  P.L. and Lai,  A. and Pope,  W.B. and Nguyen,  N. and Ji,  M. and Cloughesy,  T.F. and Ellingson,  B.M.},
	year = {2019},
	pages = {979-986}
}

@article{Xu2002,
	title = {Acidic Extracellular pH Induces Vascular Endothelial Growth Factor (VEGF) in Human Glioblastoma Cells via ERK1/2 MAPK Signaling Pathway},
	volume = {277},
	ISSN = {0021-9258},
	url = {http://dx.doi.org/10.1074/jbc.m108347200},
	DOI = {10.1074/jbc.m108347200},
	number = {13},
	journal = {Journal of Biological Chemistry},
	publisher = {Elsevier BV},
	author = {Xu,  Lei and Fukumura,  Dai and Jain,  Rakesh K.},
	year = {2002},
	pages = {11368-11374}
}

@article{Oudin2016,
	title = {Tumor Cell-Driven Extracellular Matrix Remodeling Drives Haptotaxis during Metastatic Progression},
	volume = {6},
	ISSN = {2159-8290},
	url = {http://dx.doi.org/10.1158/2159-8290.CD-15-1183},
	DOI = {10.1158/2159-8290.cd-15-1183},
	number = {5},
	journal = {Cancer Discovery},
	publisher = {American Association for Cancer Research (AACR)},
	author = {Oudin,  Madeleine J. and Jonas,  Oliver and Kosciuk,  Tatsiana and Broye,  Liliane C. and Guido,  Bruna C. and Wyckoff,  Jeff and Riquelme,  Daisy and Lamar,  John M. and Asokan,  Sreeja B. and Whittaker,  Charlie and Ma,  Duanduan and Langer,  Robert and Cima,  Michael J. and Wisinski,  Kari B. and Hynes,  Richard O. and Lauffenburger,  Douglas A. and Keely,  Patricia J. and Bear,  James E. and Gertler,  Frank B.},
	year = {2016},
	pages = {516-531}
}

@article{Liu2019,
	title = {Glioma Cell Migration Dynamics in Brain Tissue Assessed by Multimodal Optical Imaging},
	volume = {117},
	ISSN = {0006-3495},
	url = {http://dx.doi.org/10.1016/j.bpj.2019.08.010},
	DOI = {10.1016/j.bpj.2019.08.010},
	number = {7},
	journal = {Biophysical Journal},
	publisher = {Elsevier BV},
	author = {Liu,  Chao J. and Shamsan,  Ghaidan A. and Akkin,  Taner and Odde,  David J.},
	year = {2019},
	pages = {1179-1188}
}

@article{Giese1996,
	title = {Glioma Invasion in the Central Nervous System},
	volume = {39},
	ISSN = {1524-4040},
	url = {http://dx.doi.org/10.1097/00006123-199608000-00001},
	DOI = {10.1097/00006123-199608000-00001},
	number = {2},
	journal = {Neurosurgery},
	publisher = {Ovid Technologies (Wolters Kluwer Health)},
	author = {Giese, Alf and Westphal, Manfred},
	year = {1996},
	pages = {235-252}
}

@article{Hormuth2021,
	title = {Biologically-Based Mathematical Modeling of Tumor Vasculature and Angiogenesis via Time-Resolved Imaging Data},
	volume = {13},
	ISSN = {2072-6694},
	url = {http://dx.doi.org/10.3390/cancers13123008},
	DOI = {10.3390/cancers13123008},
	number = {12},
	journal = {Cancers},
	publisher = {MDPI AG},
	author = {Hormuth,  David A. and Phillips,  Caleb M. and Wu,  Chengyue and Lima,  Ernesto A. B. F. and Lorenzo,  Guillermo and Jha,  Prashant K. and Jarrett,  Angela M. and Oden,  J. Tinsley and Yankeelov,  Thomas E.},
	year = {2021},
	pages = {3008}
}

@article{Winkler2018,
	title = {Singular structure formation in a degenerate haptotaxis model involving myopic diffusion},
	volume = {112},
	ISSN = {0021-7824},
	url = {http://dx.doi.org/10.1016/j.matpur.2017.11.002},
	DOI = {10.1016/j.matpur.2017.11.002},
	journal = {Journal de Math\'ematiques Pures et Appliqu\'ees},
	publisher = {Elsevier BV},
	author = {Winkler,  Michael},
	year = {2018},
	pages = {118-169}
}

@article{Winkler2017,
	title = {Global weak solutions to a strongly degenerate haptotaxis model},
	volume = {15},
	ISSN = {1945-0796},
	url = {http://dx.doi.org/10.4310/CMS.2017.v15.n6.a5},
	DOI = {10.4310/cms.2017.v15.n6.a5},
	number = {6},
	journal = {Communications in Mathematical Sciences},
	publisher = {International Press of Boston},
	author = {Winkler,  Michael and Surulescu,  Christina},
	year = {2017},
	pages = {1581-1616}
}

@article{Heihoff2023,
	title = {Global solutions to a haptotaxis system with a potentially degenerate diffusion tensor in two and three dimensions},
	volume = {36},
	ISSN = {1361-6544},
	url = {http://dx.doi.org/10.1088/1361-6544/acadcb},
	DOI = {10.1088/1361-6544/acadcb},
	number = {2},
	journal = {Nonlinearity},
	publisher = {IOP Publishing},
	author = {Heihoff,  Frederic},
	year = {2023},
	pages = {1245-1278}
}

@article{Kolbe2021,
	title = {Modeling multiple taxis: Tumor invasion with phenotypic heterogeneity,  haptotaxis,  and unilateral interspecies repellence},
	volume = {26},
	ISSN = {1553-524X},
	url = {http://dx.doi.org/10.3934/dcdsb.2020284},
	DOI = {10.3934/dcdsb.2020284},
	number = {1},
	journal = {Discrete \& Continuous Dynamical Systems - B},
	publisher = {American Institute of Mathematical Sciences (AIMS)},
	author = {Kolbe,  Niklas and Sfakianakis,  Nikolaos and Stinner,  Christian and Surulescu,  Christina and Lenz,  Jonas},
	year = {2021},
	pages = {443-481}
}

@article{Tao2019,
	title = {Large time behavior in a forager-exploiter model with different taxis strategies for two groups in search of food},
	volume = {29},
	ISSN = {1793-6314},
	url = {http://dx.doi.org/10.1142/S021820251950043X},
	DOI = {10.1142/s021820251950043x},
	number = {11},
	journal = {Mathematical Models and Methods in Applied Sciences},
	publisher = {World Scientific Pub Co Pte Ltd},
	author = {Tao,  Youshan and Winkler,  Michael},
	year = {2019},
	pages = {2151-2182}
}

@article{Winkler2019,
	title = {Global generalized solutions to a multi-dimensional doubly tactic resource consumption model accounting for social interactions},
	volume = {29},
	ISSN = {1793-6314},
	url = {http://dx.doi.org/10.1142/S021820251950012X},
	DOI = {10.1142/s021820251950012x},
	number = {03},
	journal = {Mathematical Models and Methods in Applied Sciences},
	publisher = {World Scientific Pub Co Pte Ltd},
	author = {Winkler,  Michael},
	year = {2019},
	pages = {373-418}
}

@article{Tao2023,
	title = {Small-signal solutions to a nonlocal cross-diffusion model for interaction of scroungers with rapidly diffusing foragers},
	volume = {33},
	ISSN = {1793-6314},
	url = {http://dx.doi.org/10.1142/S0218202523500045},
	DOI = {10.1142/s0218202523500045},
	number = {01},
	journal = {Mathematical Models and Methods in Applied Sciences},
	publisher = {World Scientific Pub Co Pte Ltd},
	author = {Tao,  Youshan and Winkler,  Michael},
	year = {2023},
	pages = {103-138}
}

@article{Xu2020,
	title = {Asymptotic Behavior of a Tumor Angiogenesis Model with Haptotaxis},
	volume = {8},
	ISSN = {2227-7390},
	url = {http://dx.doi.org/10.3390/math8050664},
	DOI = {10.3390/math8050664},
	number = {5},
	journal = {Mathematics},
	publisher = {MDPI AG},
	author = {Xu,  Chi and Wang,  Yifu},
	year = {2020},
	pages = {664}
}

@article{Anderson1998,
	title = {Continuous and Discrete Mathematical Models of Tumor-induced Angiogenesis},
	volume = {60},
	ISSN = {0092-8240},
	url = {http://dx.doi.org/10.1006/bulm.1998.0042},
	DOI = {10.1006/bulm.1998.0042},
	number = {5},
	journal = {Bulletin of Mathematical Biology},
	publisher = {Springer Science and Business Media LLC},
	author = {Anderson,  A},
	year = {1998},
	pages = {857-899}
}

@article{MORALESRODRIGO2013,
	title = {Global existence and asymptotic behavior of a tumor angiogenesis model with chemotaxis and haptotaxis},
	volume = {24},
	ISSN = {1793-6314},
	url = {http://dx.doi.org/10.1142/S0218202513500553},
	DOI = {10.1142/s0218202513500553},
	number = {03},
	journal = {Mathematical Models and Methods in Applied Sciences},
	publisher = {World Scientific Pub Co Pte Ltd},
	author = {Morales-Rodrigo,  Cristian and Tello, J. Ignacio},
	year = {2013},
	pages = {427-464}
}

@inbook{Chen2017,
	title = {Modelling of Tumour-Induced Angiogenesis Influenced by Haptotaxis},
	ISBN = {9789811039577},
	ISSN = {2196-887X},
	url = {http://dx.doi.org/10.1007/978-981-10-3957-7_9},
	DOI = {10.1007/978-981-10-3957-7_9},
	booktitle = {Emerging Trends in Neuro Engineering and Neural Computation},
	publisher = {Springer Singapore},
	author = {Chen,  Wei and Zhang,  Li and Liu,  Chengyu and Hossain,  Alamgir},
	year = {2017},
	pages = {173-191}
}

@article{Pang2019,
	title = {Asymptotic behavior of solutions to a tumor angiogenesis model with chemotaxis-haptotaxis},
	volume = {29},
	ISSN = {1793-6314},
	url = {http://dx.doi.org/10.1142/S0218202519500246},
	DOI = {10.1142/s0218202519500246},
	number = {07},
	journal = {Mathematical Models and Methods in Applied Sciences},
	publisher = {World Scientific Pub Co Pte Ltd},
	author = {Pang,  Peter Y. H. and Wang,  Yifu},
	year = {2019},
	pages = {1387-1412}
}

\end{document}